\documentclass[a4paper,11pt]{amsart}

\usepackage{amscd,amsfonts,amsmath,amssymb}

\input xy
\xyoption{all}
\usepackage[all]{xy}

\newcommand{\cal}{\mathcal}

\newcounter{subsub}[subsection]

\def\hfl#1#2{\smash{\mathop{\hbox to 12mm{\rightarrowfill}}
\limits^{\scriptstyle#1}_{\scriptstyle#2}}}

\long\def\InsertFig#1 #2 #3 #4\EndFig{
\hbox{\hskip #1 mm$\vbox to #2 mm{\vfil\includegraphics{#3}}#4$}}
\long\def\LabelTeX#1 #2 #3\ELTX{\rlap{\kern#1mm\raise#2mm\hbox{#3}}}

\newcommand\bn {{\bigskip \noindent}}

\newtheorem{theorem}{Theorem}[section]
\newtheorem{lemma}[theorem]{Lemma}
\newtheorem{example}[theorem]{Example}
\newtheorem{corollary}[theorem]{Corollary}
\newtheorem{proposition}[theorem]{Proposition}
\newtheorem{question}[theorem]{Question}
\newtheorem{remark}[theorem]{Remark}

\newtheorem{definition}[theorem]{Definition}
\newtheorem{conjecture}[theorem]{Conjecture}

\newtheorem{examples}[theorem]{Examples}
\newtheorem{abundance conjecture}[theorem]{Abundance Conjecture}
\newtheorem{simplified lemma}[theorem]{Simplified Lemma}

\newtheorem{notation}[theorem]{Notation}
\newtheorem{main questions}[theorem]{Main Questions}

\newcommand\ra{{\rightarrow}}
\newcommand\mero{{\dashrightarrow}}

\begin{document}
\title{SPECIAL VARIETIES AND CLASSIFICATION
THEORY}
\author{Fr\'ed\'eric CAMPANA}
\maketitle

\begin{center}
\begin{minipage}{130mm}
\scriptsize

{\bf Abstract}: We describe two structure theorems in the birational geometry of complex projective 
manifolds, analogues of the structure theorems for Lie algebras,  
reducing these first to semi-simple and solvable ones, the solvable being iterated extensions of abelian ones.
The roles of semi-simple, solvable and abelian are respectively played here by the orbifolds of general type, special manifolds, 
and third: manifolds either rationally connected or with $\kappa=0$. 

A complex projective manifolds $X$, is said to be 
{\bf special} if it has no map onto an orbifold of general type,
 the orbifold structure on the base being given by the divisor of multiple fibres (suitably defined).
We show that  manifolds 
either rationally connected or with $\kappa=0$ are special. 
We also show that $X$ is special if and only if it has no Bogomolov sheaf, which 
are rank one coherent subsheaves of maximal Kodaira dimension $p$ of $\Omega_X^p$, for some $p>0$. 

An orbifold version of Iitaka's $C_{n,m}$ additivity conjecture implies that special manifolds 
have a canonical decomposition as towers of fibrations with 
generic (orbifold) fibres either rationally connected or with zero Kodaira dimension. This decomposition is obtained by alternatively 
iterating orbifold Iitaka fibrations and suitably defined rational quotients. 

For any complex projective manifold $X$, we next construct a functorial fibration $c_X:X\ra C(X)$ called its {\bf core}$^1$, 
such that its general fibre is special, and contains any special subvariety of $X$ meeting it. 
The core is thus the constant (resp. identity) map iff $X$ is special (resp. of general type). 

As a consequence of the special case of the orbifold version of Iitaka conjecture $C_{n,m}$ when the base orbifold is of general type, 
we show that the core has an orbifold base of general type if $X$ is not special.

The core thus canonically decomposes any $X$ into its ``components" in the dichotomy special (the fibres) versus general type 
(in the orbifold sense: the base). Special manifolds being (conjecturally) orbifold combinations of the 
first two types of the fundamental trichotomy (rationally connected, $\kappa=0$, general type) of classification theory.

This decomposition is deeply linked with other aspects of classification theory: we indeed conjecture, and show in some few cases, 
 that special manifolds have  an almost abelian fundamental group, and are exactly the ones having zero 
Kobayashi metric, or a potentially dense set of rational points over any field of definition finitely generated over 
$\mathbb Q$. This immediately leads to a natural extension of Lang's conjectures to arbitrary $X$'s (and even to orbifolds).

This decomposition gives a simple synthetic view of the structure of arbitrary $X$'s, and 
indicates that the natural frame of classification theory is the category of orbifolds, to which our observations should be 
extended. 

The present observations actually work for any compact K\" ahler $X$, which are considered here.
Most of them can also be applied to other algebraically closed fields 
(with the usual problems of separability in positive characteristics).\\

$^1$ From the french ``coeur''. A case to which A.  Dumas's statement applies well: 
``Comme l'anglais n'est que du fran\c cais mal prononc\' e..." (in ``Vingt ans apr\`es", Chap LXVII: ``Londres")

\end{minipage}
\end{center}

 \tableofcontents
\section*{\bf Introduction} 

The structure of arbitrary complex Lie algebras is described to a large extent by the following two theorems: 
the first one says that any such Lie algebra is 
in a unique way an extension of a semi-simple Lie algebra by a solvable one. Next, 
this extension can further be decomposed, also in a canonical functorial way, as a composition of extensions with abelian kernels. 

The objective of this paper is to establish two similar structure results for arbitrary complex projective varieties $X$, in birational geometry. 

The first theorem asserts the existence and uniqueness of a fibration $c_X:X\mero C(X)$, called the core of $X$, such that its general fibres are 
special, while its orbifold base is of general type. The second result (only conjectural) asserts that the core 
can be obtained in a functorial canonical way  as a composition of orbifold rational quotients and Iitaka fibrations. This statement depends on an 
orbifold version of Iitaka's conjecture $C_{n,m}$.

One thus gets a canonical decomposition of $X$, the elementary components of which are the varieties of pure type in the fundamental trichotomy 
of algebraic geometry (rationally connected, $\kappa=0$, general type). The orbifold structures seem to be unavoidable, and to lie at the heart of 
classification theory, which should thus deal with this category of objects.\\
 
 We now describe the content of the different sections.\\

In \S 1, we introduce the notion of orbifold base of a holomorphic fibration $f:X\ra Y$ 
($f$ is thus surjective, and its fibres are connected) between 
compact K\" ahler manifolds. This is simply (\ref{multdiv}) a $\mathbb Q$-divisor on $Y$, suitably 
defined by the multiple fibres of $f$. This orbifold structure carries naturally a canonical bundle, a Kodaira dimension and 
a fundamental group. 

The Kodaira dimension $\kappa(Y,f)$ is the minimum of the Kodaira dimensions of the orbifold bases, when $f'$ runs 
over all models of $f$, which are the fibrations bimeromorphically equivalent to $f$. 
This minimum is achieved by $f$ itself when $\kappa(Y)\geq 0$ (by \ref{k>=0}), but 
may be smaller in general. 

We next show (\ref{behKod}) that $\kappa(Y,f)$ behaves as in the ``classical" case (where $f=id_X$) under composition with 
generically finite maps $u:X'\mero X$ and Stein factorisation.

In the next section, we define canonically a rank one subsheaf $F_f$ of $\Omega _X^p$ by saturating the inverse image of $K_Y$. 
This sheaf is an intrinsic invariant of the equivalence class of fibration. Its Kodaira dimension is $\kappa(Y,f)$. 
As an aplication, the orbifold base of $f$ has Kodaira dimension $\kappa(Y,f)$ if $f$ is what we call neat (see \ref{neat}).

We next show (\ref{semicont}) the countable upper semi-continuity of $\kappa(Y,f)$ in families.

Finally, we show how to compute, on suitable models, the orbifold structure 
on the base of a composition $g\circ f$ of two fibrations, from 
the orbifold bases of $g$ and $f$. This computation plays an essential role 
in the proof of the shown case of the orbifold additivity conjecture in \S 4. It is also the clue to the definition 
of the orbifold base of a fibration between orbifolds. See section \ref{orbcat}.\\

In \S 2, we first define general type fibrations $f:X\mero Y$ as the ones having $\kappa(Y,f)=dim(Y)>0$, 
special manifolds $X$ are defined as the ones having no $f$ of general type, and special fibrations 
as the ones having special general fibres. We list without proof some examples of special manifolds. The two most 
important ones are the rationally connected manifolds, and manifolds with $\kappa=0$. 

Notice however 
that being special does not give 
any restriction on the Kodaira dimension, except for the top one. For example, any elliptic surface with base 
elliptic or rational is special if it has a section. Somewhat unexpectedly maybe, the moduli of fibrations do not play any role in our 
considerations.

We next show various geometric properties of special and general type fibrations. The most important one (\ref{domin'}) is that
any special fibration $f:X\mero Y$ dominates any general type fibration $g:X\mero Z$, in the sense that there exists 
$\phi:Y\mero Z$ such that $g=\phi\circ f$. 

From which one concludes that on any $X$ there exists at most one fibration both special and of general type. 
One of the main results of the present paper is the {\bf existence} of such a fibration. We get it by the two possible approaches:
 either from ``above", as the ``lowest special" fibration on $X$, obtained by geometric means (see section \ref{core}), 
or from ``below", as the ``highest general type 
fibration" on $X$ (see \ref{core'}). The second approach is less geometric and hides the motivations 
and origins of the constructions made, but has technical advantages similar to those of the axiomatic 
method, turning theorems into definitions, once the foundations have been laid. (See footnote 2 below)
\footnote{B.Russell once compared these advantages to the ones of stealing over honest work}

The next subsection shows the bijective correspondance between Bogomolov sheaves and general type fibrations (\ref{gtBog}). From 
which we conclude that special manifolds are characterised by the absence of such sheaves. Applying a result from [C95], we 
obtain a first simple proof that $X$is special if either rationally connected or with $c_1(X)=0$.

We next show (\ref{specZR}), among other things, that the general fibre $X_s$ of a 
fibration $h:X\ra S$ is special if it is special 
for any $s$ in a set $E$ not contained in a countable union of Zariski closed subsets with empty interior of $Z$. 
This property, which we call Zariski regularity (for specialness) is, together with \ref{domin'}, one of the main ingredients 
in the 
proof that the core of a manifold is a special fibration (\ref{c_X}).

This section ends with a brief sketch of the extension of the considerations of the present paper to the orbifold category. It seems that most of them 
extend without any difficulty to the orbifold context, including the notions of orbifold modification and of differential 
form on an orbifold $(Y/D)$ with smooth $Y$ and 
the support of $D$ of normal crossings. In particular the notion of Bogomolov sheaf and the sheaf $F_f$ 
associated to a fibration between orbifolds make sense. The geometric aspect seems more delicate to handle (see section \ref{orbrc} for 
some subtle problems on surface orbifolds which are classical in the non-orbifold context). 

The notion of orbifold differential form we introduce interpolates between the usual one when $D$ is empty, and the classical notion of 
$\Omega_Y(log (D))$ when $D$ is reduced, which corresponds to the limiting case of infinite multiplicities. It seems that this topic 
desrves by itself further developments.\\

In \S 3, we construct and start to study the core $c_X:X\mero C(X)$ of a manifold. 
We first show (\ref{c_X}) that its general fibres are 
special (this fails for general singular varieties, and does not follow from the original definition, 
obtained applying the construction of meromorphic quotients 
recalled in section \ref{geomquot}). 
From this result, we immediately get (\ref{RCspec}) that rationally connected manifolds 
are special (simply because $\mathbb P^1$ is special). The notions around rational connectedness are recalled. 
We next deduce in (\ref{CR=C}) from [G-H-S] that if $R(X)$ is the rational quotient of $X$ 
(see \ref{ratquot'} for this notion), then $C(R(X))=C(X)$. 

We next describe (\ref{coresurf} and \ref{corenonspec3f}) the core and list the special manifolds in dimensions $2$ and $3$, after 
having introduced (see \S 3.6) the ``higher Kodaira dimensions" of a (compact complex) manifold. From the description 
so obtained, we deduce that in these cases the core is a fibration of general type (when $X$ is not special). The 
``decomposition theorem" (\ref{decconj}) asserts that this is true in any dimension. \\

The next \S 4 states (see \ref{$C^{orb}_{n,m}$}) the orbifold version of Iitaka conjecture $C_{n,m}$, 
and shows (\ref{cnmtheor}) the special 
case when its orbifold base is of general type. 
This result is one of the main technical tools of the present paper. 
The proof consists in extending the classical weak positivity 
results of Fujita, Kawamata and Wiehweg for direct image sheaves of 
pluricanonical forms to the orbifold situation, by suitably introducing the orbifold 
divisor into the proofs.\\

\S 5 contains the geometric consequences of the orbifold additivity result (theorem \ref{cnmtheor}): 
Manifolds with $\kappa=0$ are special (\ref{k=0}), 
the Albanese map of a special manifold is surjective, connected, and has no multiple 
fibres in codimension one (\ref{albsurj}). 
We also show (\ref{ess=n} (resp. \ref{ess=n-1})), that the essential dimension $ess(X):=dim(C(X))$ of $X$ is equal to 
its dimension $dim(X)=n$ (resp. to $(n-1)$) iff $X$ is of general type (resp. iff $X$ has a fibration of general type with generic 
fibre a curve either rational or elliptic). 

These are two special cases of the ``decomposition theorem" already stated above. This theorem is established 
in (\ref{decconj}). We also show that the decomposition theorem implies the invariance of the essential 
dimension under finite \' etale covers, in particular that such covers of a special manifold are still special.
 We finally observe (see subsection \ref{essbog}) that 
the decomposition theorem asserts that any nonspecial $X$ has a 
(unique) fibration both special and of general type, and a unique maximum Bogomolov sheaf.

We then show that the Stein factorised product of two fibrations of general type is itself of general type \footnote{Maybe surprisingly, it was not no 
noticed, since [Bo78] that one could so naturally define a fibre product of two Bogomolov sheaves, although the 
techniques used in the present paper are available since two decades. This is because the question is obvious if one starts from ``above": the geometry 
of orbifolds and general type fibrations, but very hidden and unnatural if one starts from the sheaves.} This gives, together with 
another application of \ref{cnmtheor}, the second construction of the core ``from below", as the ``highest general type" fibration on a given $X$.\\

The \S 6 presents various conjectures related to the core. We first start with defining and comparing some orbifold variants 
of the notion of rational connectedness. One of them (the weakest, possibly equivalent to the strongest of these notions) turns out 
to give a good substitute of the notion of rational quotient in the orbifold context, provided the orbifold version 
of conjecture $C_{n,m}$ holds.   This implies (see section 
\ref{c=rj^n}) that the core can be canonically and functorially decomposed into a succession of orbifold Iitaka-Moishezon
\footnote{Usually called Iitaka fibration. The construction obtained 
independently by B.Moishezon [Mois] remained apparently unnoticed} fibrations 
and rational quotients (defined in section \ref{orbRJ}). See \S 6 for more details. This construction thus ``only" depends 
on $C_{n,m}$ in its orbifold version.

So that the core actually gives directly the final step of this decomposition, without producing the intermediate steps. 

In subsection \ref{definv}, we conjecture 
the deformation invariance of all the invariants introduced (in particular: $ess(X)$ and the higher Kodaira dimensions. 
So that in particular, deformations of special manifolds should be special). The next subsections deal with two other similar theories built 
on related notions of specialness and general type fibrations? They seem to both differ from the one presented here, but presently no example is known 
that distinguishes them.\\ 

In \S 7, we consider the fundamental group, and conjecture (``abelianity conjecture" \ref{abconj}) that  a 
special manifold should have an almost abelian fundamental group. This conjecture is supported and motivated by the fact that
rationally connected manifolds are simply connected, that this conjecture is standard for the case $\kappa=0$, and the preceding 
conjectural orbifold decomposition of any special manifold as a tower of fibrations with fibres of these two types. This conjecture seems to 
sum up all conjectures with the same conclusion.

We show  that this conjecture is true for linear and torsionfree solvable representations of the fundamental group, as an 
immediate consequence of previous result by various authors. As usual, we extend this cojecture to the orbifold case, anyway necessary 
to solve the non-orbifold one in higher dimensions.\\

In \S 8, we establish and apply an orbifold version (see (\ref{ko})) of the famous extension theorem of Kobayashi-Ochiai, asserting 
roughly that meromorphic nondegenerate (ie: somewhere submersive) maps $\psi$ from a dense 
Zariski open subset $U$ of a complex manifold $V$
to a variety of general type $Y$ do not have essential singularities along the boundary divisor $D:=(V-U)$. Our version essentially says this holds true 
if $\psi$ factorises as $\psi=f\circ \phi$, with $f:X\mero Y$ is a fibration of general type, and $\phi:U\mero X$ meromorphic. 

The proof is an orbifold version of the proof of Kobayashi-Ochiai. This result implies among other more general results that a 
manifold $X$ is special if there is a meromorphic nondegenerate map $\phi:\mathbb C^n\mero X$. This is a special instance 
where the Kobayashi pseudometric of $X$ vanishes identically if $\phi$ has dense image in the metric topology. 

The next \S 9 discusses two conjectures concerning the Kobayashi pseudometric of a compact K\" ahler $X$. 
The first one (\ref{conj3h}) states that special manifolds are exactly the ones having zero Kobyashi pseudometric $d_X$. 
The second one (\ref{splitd}) states that for any $X$, $d_X$ is the lift by the core $c_X:X\mero C(X)$ 
of the orbifold Kobayashi pseudometric of the base orbifold of $c_X$, defined in \ref{orbd}, and that this orbifold pseudometric 
is a metric outside some algebraic subset $S$ of $C(X)$. Recall that the orbifold base of $c_X$ is of general type (if not a point), 
so this second part of the conjecture \ref{splitd} is the orbifold version of Lang's hyperbolicity conjecture. 

The next \S 10 (resp. \S 11) is the exact analog of \S 9 for arithmetic geometry (resp. the function field version of \S 10)). 
We do not give a more detailed description, here, because our discussion rests only on known results by many authors. 

This paper ends with the \S 12, which is a technical appendix exposing (with proofs) the construction 
of meromorphic quotients of equivalence relations of geometric origin (\ref{quot}) identifying points connected by chains 
with irreducible components in a given family of subvarieties.

This construction appeared it seems for the first time in [C81]. We follow here the original argument, 
in a slightly simplified form. (This construction was exposed also in [Ko 96] by a different approach, 
corrected and exposed in a detailed way in [De]).

We complete this construction by introducing the notions of Zariski regularity (section \ref{Zreg}) and 
stability (section \ref{stab}), of constant use in this context. In particular, the notion of stability 
permits to give a criterion (\ref{quotstab}) for the quotient to have fibres in the class as subvarieties 
used to construct the chains providing the equivalence relation.\\

The present version is the sixth one. It differs from the original one in two aspects: first by 
the definition of multiple fibres (based on inf, no longer gcd. See \ref{multdiv}), 
second, and more importantly, by the proof that the base orbifold of the core is of general type, 
which was only conjectured in the preceding versions.
The first point was also observed, independently of me, by S.Lu, after he read the first versions of the present paper: 
one can replace the gcd-multiplicities 
by the present ones without changing any statement or proof (except for a small 
refinement in the proof of the former version of proposition \ref{injsheav}).
The new definition is better behaved with respect to hyperbolicity questions and Bogomolov sheaves (but seemingly not with respect to the 
fundamental group, for which $gcd$-multiplicities work very well).\\

{\bf Aknowlegements:} I would like to thank F. Bogomolov, J.P. Demailly, and C. Voisin for 
discussions which permitted me to improve sections 7 and 10.  P.Eyssidieux asked me about the function 
field version of \S 10, which lead to \S 11. F.Bogomolov noticed that general type fibrations 
might be linked with what is called here Bogomolov sheaves; this link is established in \ref{gtBog} (One direction was in fact shown
 in the first version in the transcendantal context of section \ref{orbko}). Moreover, after a first 
version of this work was written, I got further stimulating comments and references 
included here from D.Abramovitch, J.L. Colliot-Th\'el\`ene and B.Totaro 
(who independently suggested me to include Conjecture IV$_A$).

I thank heartily thank L.Bonavero, S.Druel, J.P.Demailly, L.Manivel, C.Mourougane,M.Paun and E.Peyre for 
their carefull reading of {\S 8,9,12,14,12,12 and 7} of the previous version, (which are now \S 2,4,8,12,8,8 and 10), 
which they exposed in Grenoble in 
July of $2002$. And also for pointing many mistakes, gaps and inaccuracies. They also suggested me to reorganise 
the exposition according to the technical order of proofs. The present version thus differs considerably from the 
preceding ones in this respect. Further, L.Bonavero and S.Druel provided me with considerable help in the 
manipulation of Latex.

Special thanks are due to E.Viehweg, who gave me a decisive hint in the proof of the crucial Lemma \ref{viehw}.

\section{\bf Orbifold Base of a Fibration.}

\subsection{\bf Fibrations}

Before giving definitions, let us start with

\subsubsection{\bf A motivating example}

\begin{example}\label{amotexamp}

Let $X_0:= E\times \mathbb P_{1}(\Bbb
C)$, where $E$ is an elliptic curve. Let $X:=\tilde {X}/j$, where
$\tilde{X}:=E\times C$, where $C$ is a hyperelliptic curve and
$j=t\times h$ is the involution acting diagonally on $\tilde X$ by
a translation $t$ of order 2 on $E$, and by the hyperelliptic
involution $h$ on $C$. Then both $X_0$ and $X$ have natural fibrations
on $\mathbb P_{1}(\mathbb C)$ with generic fibre $E$. 

\end{example}

They cannot be
distinguished by this information, although they differ fundamentally
at the
levels of Kodaira dimension, fundamental group, Kobayashi pseudo
metric and arithmetics (for appropriate choice of $E$). If, however,
we
take into account the multiple fibres of the fibrations onto $\Bbb
P_{1}(\mathbb C)$, we see that in the second case (of $X$), the base is
not really
$\mathbb P_{1}(\mathbb C)$, but rather the orbifold of general type $C/h$.

This we shall now generalise.
\subsubsection{\bf Fibrations} \label{fib}

A {\bf fibration} $f:X\mero Y$ is a 
surjective ({\it i.e.} dominant) meromorphic
map with connected fibres between normal
irreducible compact complex analytic spaces $X$ and $Y$ (see below for 
the precise definition of the fibres of 
a meromorphic map). 

Of course, this fibration is said to be holomorphic (or regular) if so is the map $f$. 

Another
fibration $f':X'\mero Y'$ is said to be {\bf equivalent} to $f$ if there
exists bimeromorphic maps $u:X\mero X'$ and $v:Y\mero Y'$ such that
$f'\circ u=v\circ f$. We denote by $F$ or $X_y$ the generic fibre of
$f$. We also say that $f'$ is a {\bf model} or a {\bf representative} of $f$
(we shall not distinguish between $f$ and its equivalence class).

A fibration $f:X\mero Y$ canonically defines (see \cite{C85}) a meromorphic
map $\phi_f:Y\mero \mathcal C (X)$,
with $\mathcal C (X)$ the Chow-Scheme of $X$, by sending the generic 
$y\in Y$ to the point of
$\mathcal C (X)$ parametrising the reduced fibre $X_y$ of $f$ over $y$. 
Let $\Phi_f\subset \mathcal C (X)$
be the image of $Y$ by $\phi_f$: it is a compact irreducible analytic 
subset of $\mathcal C (X)$ bimeromorphic
to $Y$ such that its incidence graph is bimeromorphic to $Y$.
The above correspondance induces a bijective map between equivalence 
classes of fibrations and
compact irreducible analytic subsets of $\mathcal C (X)$ with incidence 
graph bimeromorphic to $X$.\\

\begin{definition} \label{neat} Assume that $f:X\ra Y$ is a holomorphic fibration. An irreducible Weil divisor $D$ on $X$ is said to be 
{\bf $f$-exceptional} if $f(D)$ has codimension $2$ or more in $Y$. We say that $f:X\ra Y$ is {\bf neat} 
if $f$ is holomorphic, $X,Y$ are smooth, and if there moreover exists a 
bimeromorphic holomorphic map $u:X\ra X'$ with $X'$ smooth such that each $f$-exceptional divisor of $X$ 
is also $u$-exceptional.
\end{definition}

By Raynaud and Hironaka flattening theorems, any fibration $f_0:X_0\mero Y_0$ has a 
neat model, in which $X'$ may be choosen 
arbitrarily (bimeromorphic to the domain $X_0$ of the initial fibration $f_0$). 
Indeed: first blow-up $X'$ in such a way that 
$f_0\circ b:X'\mero Y_0$ is holomorphic, where $b:X'\mero X_0$ is bimeromorphic; 
then flatten $(f_0\circ b)$ by blowing-up $Y_0$ to get a smooth $Y$. Finally, 
take a smooth model $X$ of $X'\times _{Y_0}Y$. \\

Assume the fibration $f$ is holomorphic. We say that the holomorphic fibration $f':X'\ra Y'$ 
with $X',Y'$ smooth {\bf dominates $f$} 
if there exists a commutative diagram in which $u,v$ are bimeromorphic; obviously, $f'$ is then equivalent to $f$:

\centerline{
\xymatrix{ X' \ar[r]^{u}\ar[d]_{f'} & X\ar[d]^{f}\\  
Y'  \ar[r]^{v} & Y\\ 
}
}

Observe that any $f'$ dominating a neat fibration is itself neat.

More generally, we say that the fibration $f:X'\mero Y'$ dominates the fibration $f:X\mero Y$ if there exists 
a diagram as above with $u, v$ meromorphic and bimeromorphic. This notion is well-defined on equivalence classes.\\

Let $f:X\ra Y$ be a holomorphic fibration with $X,Y$ smooth. We say (as in \cite{Vi82}) that $f$ is 
{\bf prepared} if the locus $Y^*\subset Y$ of points $y$ with a smooth $f$-fibre $X_y$ has a complement 
contained in a normal crossing divisor $D\subset Y$ such that, moreover, $f^{-1}(D)$ is also a 
divisor of normal crossings in $X$.

Any fibration is dominated by a prepared model, again by an immediate application of Hironaka's results. 
The models of fibrations considered can thus always be assumed to be prepared.

\subsubsection{\bf Multiplicities} \label{mult}

Assume now $X,Y$ to be smooth and $f$ 
to be a holomorphic fibration. Let $|\Delta|
\subset Y$ be the union of all codimension one
irreducible components of the locus of $y$'s $\in Y$ such that the
scheme-theoretic
fibre of $f$ over $y$ is not smooth. For each component $\Delta_i$ of
$ |\Delta |$, let $D_i:=\sum_{j\in J}m_{i,j} D_{i,j}$ be the union
of all components (counted with their scheme-theoretic
multiplicities) of $f^*(\Delta_i)$ which are  mapped surjectively
onto $\Delta_i$ by $f$.
Define $m_i=m(f,\Delta_i)=:\inf (m_{i,j},j\in J)$: it is called the {\bf multiplicity of
$f$ along $\Delta_i$}.

\subsubsection{\bf Multiplicity Divisor of a Fibration} \label{multdiv}

The $\mathbb Q$-divisor $\Delta(f)=\sum_{i\in
I}(1-1/m_i)\Delta_i$ on $Y$ will be called the {\bf multiplicity divisor
of $f$}. The motivation for its 
introduction comes from the examples above
and in Example \ref{galois} below.

Observe we can also write: $\Delta(f)=\sum_{D\subset Y
}(1-1/m(f,D)).D$ , where $D$ ranges over the set of all irreducible divisors of $Y$, and 
$m(f,D)$ is the {\bf multiplicity of $f$ along $D$}, defined as $m(f,D):=\inf (m_{j},j\in J(f,D))$, where 
$f^*(D):=\sum_{j\in J(f,D)}m_{j} D_{j}+R$ , and $J(f,D)$ is the set of all components of $f^*(D)$ which are  mapped surjectively
onto $D$ by $f$, while $R$ is $f$-exceptional. Indeed, $m(f,D)=1$ if $D$ is not a component of $\mid \Delta\mid$. 

\subsection{\bf Orbifolds.}

\subsubsection{\bf Notion of Orbifold} \label{orb}

An {\bf orbifold $(Y/\Delta)$} is a pair
consisting of a compact irreducible complex space $Y$ together with a Weil
$\Bbb Q$-divisor of the form:
$\Delta=:\sum_{i\in I}(1-1/m_i)\Delta_i$
for distinct prime divisors $\Delta_i$ of $Y$, and positive integers $m_i$. 
We also say that $\Delta$ is an orbifold structure on $Y$. We write $\mid \Delta\mid$ for the support of $\Delta$, in which 
the coefficient of each $\Delta_i$ is one. Or equivalently: each $m_i=+\infty$.

If $Y$ is smooth, and if the support
of $\Delta$ is a simple normal crossings divisor, 
such pairs $(Y/\Delta)$ were introduced 
and used by V.Shokurov in \cite{Sk} under the 
name of ``standard pairs''.
We shall say that the orbifold $(Y/\Delta)$ is a {\bf klt-orbifold}
if the pair $(Y,\Delta)$ is klt (this is an abbreviation for 
``Kawamata-Log-terminal"). This seems to be the right
category to consider, morphisms being the obvious ones.

The next example shows why such pairs are rather called orbifolds here. We also note that 
this term is commonly used in similar situations either in differential geometry, or when $Y$ is a curve 
(in [F-M] or [L], for example).

\begin{example}\label{galois} 
Let $f:X\ra Y$ be a finite (ramified) Galois cover
between manifolds. Let its ramification divisor be:
$\Delta=:\sum_{i\in I}(1-1/m_i)\Delta_{i}$,
the order of ramification above the generic point 
of $\Delta_i$ being $m_i$.
\end{example}

This preceding example occurs in the construction of fibrations with multiple fibres, as in \ref {amotexamp} for example.

\subsubsection{\bf Orbifold Base of a Fibration}

 \begin{definition} \label{orbbase}
Let $f$ be a holomorphic fibration as in \ref{fib}, and
$\Delta=:\Delta(f)$ be the multiplicity divisor of $f$. We call
$(Y/\Delta)$
the {\bf orbifold base} of $f$.

\end{definition}

\subsubsection{\bf Orbifold invariants and fibrations}

\begin{definition}\label{cankodfib}
Let $(Y/\Delta)$ be an orbifold.
We define its {\bf canonical bundle} as the $\Bbb Q$-divisor 
$K_{Y/\Delta}=:K_{Y}+\Delta$ on $Y$ and 
its {\bf Kodaira dimension} as $\kappa(Y/\Delta)=:\kappa(Y,K_{Y}+\Delta)$.

\end{definition}

One can also associate to an orbifold its fundamental group, considered below.

The term ``orbifold" arises in this context from the following reason: first from the special case 
when $f$ is the quotient of a fibration $f':X'\ra Y'$ by a finite group $G$ acting on some Galois cover $X'$ of 
$X$, the action of $G$ preserving the fibration $f'$, and such that $f'$ has no multiple fibre in codimension $1$. 
In this case, $(Y/\Delta(f))$ is precisely the orbifold $(Y'/G)$. We then say that the orbifold $Y/\Delta(f)$ 
has the {\bf unfolding} $(Y',G)$. When $Y$ is a curve, such an unfolding exists, except when $Y=\mathbb P^1$ and $\Delta$ 
consists of one or two punctures with different multiplicities. In the general case, such unfoldings exist locally on $Y$ 
(but not globally in general. M.Kapovitch explained me a very beautiful construction in dimension two for orbifolds of 
general type).

\subsection{\bf The Kodaira Dimension of a Fibration}

\subsubsection{\bf Kodaira Dimension} 

Define now in general, for $f:X\mero Y$ a 
fibration between 
irreducible compact complex spaces $X$ and $Y$:
$$\kappa(Y,f):=\inf\lbrace\kappa(\bar Y / \Delta(\bar f))\rbrace,$$
where $\bar f:\bar X \ra \bar Y$ 
ranges over all holomorphic fibrations $\bar f$ between 
manifolds $\bar X$ and $\bar Y$ which are 
equivalent to $f$.

We call $\kappa (Y,f)$ {\bf the Kodaira dimension of the fibration 
$f$.} Notice indeed that this notion
depends only on the equivalence class of $f$.

\begin{definition} \label{adm}

We shall say that $f:X\ra Y$ is
{\bf admissible} if it is holomorphic, with $X,Y$ smooth, and if $\kappa(Y,f)=\kappa(Y,K_{Y}+\Delta(f))$.
\end{definition}

We shall see in section \ref{diffib} that neat fibrations are admissible. In  \ref{k>=0}, we shall see that 
any $f:X\ra Y$ is admissible if $\kappa(Y)\geq 0$.

Notice that it follows from \ref{behKod} right below that 
if $f'$ dominates $f$, and if $f$ is admissible, so is $f'$.

\subsubsection{\bf Generically Finite Maps}

We consider a commutative diagram  of holomorphic surjective maps
between compact irreducible normal complex spaces :

\centerline{
\xymatrix{ X' \ar[r]^{u}\ar[d]_{f'} & X\ar[d]^{f}\\  
Y'  \ar[r]^{v} & Y\\ 
}
}

We assume moreover that $f,f'$ are holomorphic fibrations 
and $u,v$ generically finite (and such that $f\circ u=v\circ f'$).

\begin{theorem}\label{behKod} Let $f:X\ra Y$ and
 $f':X'\ra Y'$ be two holomorphic fibrations
 and let holomorphic maps $u:X'\ra X$ and $v:Y'\ra Y$ 
 such that $f\circ u=v\circ f'$ be given.
 \begin{enumerate}
  \item [1.] Assume also that $u,v$ are bimeromorphic, then
   \begin{enumerate}
    \item [a.] $\kappa(Y/\Delta(f\circ u))=\kappa(Y/\Delta(f))$
     and $\kappa(Y'/\Delta(f'))\leq\kappa(Y/\Delta(f))$,
    \item [b.] If, moreover, $\kappa(Y)\geq 0$, then equality holds, and
     $\kappa(Y/\Delta(f))=\kappa(Y,f)$.
   \end{enumerate}
  \item [2.] Assume that $u$ and $v$ are generically finite and
   surjective. Then:
   \begin{enumerate}
    \item [a.] $\kappa(Y',f')\geq\kappa(Y,f)$,
    \item [b.] $\kappa(Y',f')=\kappa(Y,f)$ if $u$ is \'etale, and $X,X'$ are smooth.
   \end{enumerate}
 \end{enumerate}
\end{theorem}

We shall prove Theorem \ref{behKod} as a consequence of
several lemmas, some of independent interest. 

A different, shorter, proof of \ref{behKod} will be given in section \ref{diffib} below.

\begin{lemma}\label{behmod}
Assume $u,v$ bimeromorphic. Then:
\begin{enumerate}
\item [a.] $\Delta(f\circ u)=\Delta(f)$ and so 
$\kappa(Y/\Delta(f\circ u))=\kappa(Y/\Delta(f))$,
\item [b.] $\Delta(f')=v^*(\Delta(f))+E$, for some $\mathbb Q$-divisor
$E$ supported on the exceptional locus of $v$. Therefore
$v_*(\Delta(f'))=\Delta(f)$ and 
$\kappa(Y'/\Delta(f')) \leq \kappa(Y/\Delta(f))$.
\end{enumerate}
\end{lemma}

{\bf Proof: }

{\bf (a)} For any $D\subset Y$, using notations of \ref{multdiv}, we have: 

$(f\circ u)^*(D)=u^*(f^*(D))=u^*(\sum_{J(D,f)}m_j.D_j+R)=\sum _{j\in J(D,f)} u^*(m_j.(\overline{D_j}+R_j))+u^*(R)$,
 
where $\overline{D_j}$ is the strict transform by $u$ of $D_j$, and $R_j$ is the $u$-exceptional part of $u^*(D_j)$. 

Obviously, $u^*(R)$ is $(f\circ u)$-exceptional. For each component $R_{jk}$ of $R_j$, its multiplicity in $(f\circ u)^*(D)$ 
is divisible by $m_j$, by the above equalities (and the factoriality of the smooth $X$), Thus $m(D,f\circ u)=m(D,f)$, and 
$\Delta(f\circ u)=\Delta(f)$, as claimed. (Observe that $R_{jk}$ does not need to be $(f\circ u)$-exceptional, in general).

{\bf (b)} By (a), we can and shall assume that $X=X'$, to ease notations. Let $D'\subset Y'$ be 
an irreducible divisor, and let $D:=v(D)$. We just need to show that $m(D,f)=m(D',f')$ if $D'$ is 
not $u$-exceptional, that is: if $D$ is a divisor of $Y$. Then: $v^*(D)=D'+E$, with $E$ an effective $v$-exceptional divisor of $Y'$. 
Thus $\sum_{J(D,f)}m_j.D_j+R=f^*(D)=(f')^*(D')+(f')^*(E)$. Now observe that $(f')^*(E)$ is $f$-exceptional, and that each $D_j$ is mapped 
surjectively onto $D$ by $f$, and so must be surjectively mapped to $D'$ by $f'$. Thus: $\sum_{J(D,f)}m_j.D_j=\sum_{J(D',f')}m_j.D'_j$, and we 
get the claim.

We thus get the first assertion of (b). The others easily follow from it:
write $E=E^{+}-E^-$, with $E^+$ and $E^-$ effective
and $v$-exceptional. We thus get:
$\kappa(Y'/\Delta(f'))\leq\kappa(Y',K_{Y'}+v^*(\Delta)+E^+)=\kappa(Y'/v^*(\Delta))=\kappa(Y/\Delta)$, 
and so the conclusion.

\begin{remark} Modifications of $X$ thus don't alter
$\kappa(Y/\Delta(f))$,
but modifications of the base may let it
decrease. See the example below,which shows that 
strict inequality can actually occur in 
Theorem \ref{behKod} 

\end{remark}

\begin{example}

Probably the simplest example is when $Y=\Bbb P^2$ and $\Delta_{red}$
is a union of $2k\geq 6$ distinct lines meeting at one point $a\in \Bbb P^2$.
Corresponding fibrations are easily constructed, as follows:
Let $p:Y^+\ra Y$ be the double cover branched exactly along $\Delta_{red}$, and
$h$ the involution of $Y^+$ exchanging the sheets of $p$.
Let $E$ be an elliptic curve and $t$ a translation of order $2$ of $E$.
Let $X^+:=E\times Y^+$, $j:=t\times h$ the diagonal involution of $X^+$,
and $X_0:=X^+/G$, where $G$ is the group of order $2$ generated by $j$.
Let $f_0:X_0\ra Y$ and $f^+:X^+\ra Y^+$ be the induced holomorphic 
fibrations. Let $d:X\ra X_0$
be a desingularisation (induced by a desingularisation of $Y^+$, for example),
and $f:=d\circ f_0:X\ra Y$.

\end{example}
Then $\Delta(f)=(1/2)\Delta_{red}$, and so 
$\kappa(Y/\Delta(f))=-\infty$ (resp. $0$; resp. $2$) if
$k\leq 2$ (resp. $k=3$; resp. $k\geq 4$).
On the other hand, let $v:Y'\ra Y$ be the blow-up of $a\in \Bbb P^2$, let
$E$ be its exceptional divisor, and $f':X'\ra Y'$
be the lifting of the meromorphic map $v^{-1}\circ f:X\ra Y'$
to a suitable modification $X'$ of $X$.
An easy local computation (for example) shows that 
$\Delta(f')=\overline {\Delta}$,
the strict transform of $\Delta(f)$
by $v$. (This is just because the normalisation
of $(Y^+\times _{Y} Y')$ does not ramify over the generic point of $E$).
We thus conclude that $\kappa(Y'/\Delta(f'))=-\infty$,
whatever $k\geq 1$. (Because ($K_{Y'}+\Delta(f'))=(k-3)(v^*(H)-E)-2E$, if
$H$ is the hyperplane line bundle on $Y$, and $v^*(H)-E$ defines 
the unique ruling of $Y'$).
In particular, $\kappa(Y'/\Delta(f')) <\kappa(Y/\Delta(f))$ if $k\geq 3$.
Finally, we observe that the same construction does not lead to this 
last (strict)
inequality when $\Delta(f)$ is a normal crossing divisor (but such other examples should exist).

\begin{corollary}\label{existadmiss} 
Let $f:X\mero Y$ a 
fibration between
irreducible compact complex spaces $X$ and $Y$.
\begin{enumerate}
\item [1.] there exists $\bar f:\bar X \ra \bar Y$ an {\bf admissible}
holomorphic model of $f$
between manifolds $\bar X$ and $\bar Y$, and bimeromorphic maps $u:\bar X\ra X$ and $v:\bar Y\ra Y$
such that $f\circ u=v\circ \bar f$,
\item [2.] for any $X'$ bimeromorphic to $X$, there exists
$u : X'' \ra X'$ bimeromorphic and $f'': X'' \ra Y''$, a 
holomorphic {\bf admissible and prepared} 
model of $f$ between manifolds $X''$ and $Y''$
such that any effective divisor $B'' \subset X''$ 
which is $f''$-exceptional, is also $u$-exceptional.     
\end{enumerate}
\end{corollary} 

{\bf Proof: } It is an easy application of Hironaka smoothing and Hironaka-Raynaud Flattening theorems.

(1) First take an admissible model $\bar f:\bar X\ra \bar Y$ of $f$. Modify $\bar Y$ to $Y'$, smooth and dominating $Y$. 
Making base change by $Y'$ over $Y$, and smoothing $X\times_{\bar Y}Y'$, we get the fibration $f'$, with the desired properties. 

(2) Start with the previous fibration $f'$. We can first modify $X'$ so that the new $X'$ dominates any given bimeromorphic model of $X$. 
Then flatten the new $f'$ by modifying $Y'$, to get $Y"$, smooth and the non-smooth locus of $f'$ 
contained in a divisor of normal crossings. As before, make base change of $Y"$ over $Y'$, and smooth $X'\times_{Y'} Y"$ 
to obtain a fibration $f":X"\ra Y"$ enjoying the claimed properties $\square$

\begin{proposition} \label{k>=0}
Assume that $\kappa(Y) \geq 0$ in Lemma 
\ref{behmod}
above. Then:
$$\kappa(Y/\Delta(f))
=\kappa(Y'/\Delta(f'))= 
\kappa(Y,f).$$

(In other words: any holomorphic model of $f$ is then admissible).
\end{proposition}

{\bf Proof:} Because $Y$ is smooth, using the notations of
the proof of Lemma \ref{behmod}, we have:
$ v^*(\Delta(f)) - bE\leq \overline {\Delta(f)} \leq \Delta(f')$, 
for some nonnegative rational number $b$,
with $\overline {\Delta(f)} $ denoting the strict transform of 
$\Delta(f)$ by $v$. Also, 
here $E$ denotes the reduced exceptional divisor of the map 
$v$, and the inequality:
$D\leq D'$ between two $\Bbb Q$-divisors $D,D'$ on $Y'$ means that 
their difference
$(D'-D)$ is an effective $\Bbb Q$-divisor.

Moreover, we have:
$K_{Y'}\geq v^*(K_Y)+aE$, for
some strictly positive rational number $a$ (here we use the smoothness of $Y$,
but $Y$ having just terminal singularities would be sufficient).

We can thus write, as $\mathbb Q$-divisors:

$K_{Y'}+\Delta(f')\geq $
$v^*(K_Y)+aE+v^*(\Delta(f))-bE=
v^*(K_Y+\Delta(f))+(a-b)E$.

We are thus finished if $(a-b)\geq 0$. So we now assume the contrary: 
$(a/b)<1$.

From: $bE\geq v^*(\Delta(f))-\overline {\Delta(f)}$, we get first: 
$aE= (a/b)bE\geq (a/b)(v^*(\Delta(f))-\overline {\Delta(f)})$, and then:

$\overline{\Delta(f)}+a.E\geq (a/b).v^*(\Delta(f))+(1-(a/b)).\overline{\Delta(f))}$.

From which we deduce:
$(K_{Y'}+\Delta(f'))\geq v^*(K_Y)+\overline {\Delta(f)}+aE\geq D$, 
with:

$D:=(a/b)v^*(K_Y+\Delta(f))+(1-a/b)(v^*(K_Y))+\overline 
{\Delta(f)} \geq (a/b)v^*(K_Y+\Delta(f))$.
Where the last inequality follows from our assumption that $\kappa(Y)\geq 0$ $\square$

\begin{lemma}\label{genfin} Let the situation be as in Theorem
\ref{behKod} above. 
Assume
$u$ and $v$ are generically finite and surjective. Then: $\kappa(Y',f')\geq\kappa(Y,f)$. 
\end{lemma}

{\bf Proof:}
We need only to show that for any admissible $f'$ as above, we can
find a holomorphic representative $f_0:X_0\ra Y_0$ of $f$ with: 
$\kappa(Y_0/\Delta(f_0))\leq\kappa(Y'/\Delta(f'))$, because then:

$\kappa(Y,f)\leq\kappa(Y_0/\Delta(f_0))\leq\kappa(Y'/\Delta(f'))=\kappa(Y',f')$.

Modifying $X',X,Y$ and $Y'$, we can and shall thus assume $f'$ to be admissible.

Let $v=v"\circ v'$ be the Stein factorisation of $v$, with $v':Y'\ra Y"$
connected (hence bimeromorphic), $v":Y"\ra Y$ finite, and $Y"$ normal.

We also denote here
with
$\Delta:= \Delta(f)$, $\Delta':=\Delta(f')$, $\Delta
":=\Delta(f")$, $\Delta_0:=\Delta(f_0)$,... the relevant multiplicity
divisors for the corresponding fibrations $f$, $f'$, $f":=v'\circ 
f'$, $f_0$,....

Let now $w:Y"\ra Y'_o$ be a modification, and define
$f'_0:=w\circ f":=w\circ v'\circ f'$, $w_o:=w\circ v':Y'\ra Y'_0$.

The relevant diagram is the following:

\centerline{
\xymatrix{
X'\ar@/_4pc/[ddr]_-{f'_{0}}\ar[rr]^-{u}\ar[d]_-{f'}\ar[dr]^-{f''} & & 
X\ar[d]^-{f}\\
Y'\ar[r]^-{v'}\ar@/_1.8pc/[rr]_(.3){v} \ar[rd]_{w_0}& Y''\ar[r]^-{v''}\ar[d]|\hole^(.3){w} 
& Y\\
& Y_{0}'
}
}

 \begin{lemma} We can and shall further assume,
modifying $Y,X$ and $X',Y'$ if needed, that such a $w$ exists, with 
$f'_o$ admissible equivalent to $f'$.

\end{lemma}

{\bf Proof} Indeed: this existence follows from the next lemma, applied to our 
initial $Y'_0,Y_0$ in place of $V,U$, where $f'_0$ dominates $f_0$ (
which means: there exists $g,h$ such that $f_0\circ h=g\circ f'_0$, and
$f_0:X_0\ra Y_0$ is equivalent to $f$, with $g:Y'_0\ra Y_0$ and 
$h:X'_0\ra X_0$ generically finite):

\begin{lemma} Let $r:V\ra U$ be generically finite surjective, 
between irreducible normal compact complex spaces. There exist modifications
$n:V'\ra V$ and $m:U'\ra U$ such that:

(1) $s:=m^{-1}\circ r\circ n: V'\ra U'$ is holomorphic

(2) if $s=r"\circ r'$ is the Stein factorisation of $s$, then there 
exists a (holomorphic) factorisation $n":V"\ra V$ of
$m\circ r":V"\ra U $ through $r$ (ie: $r\circ n"=m\circ r"$). Here 
$r':V'\ra V"$ is connected and $r":V"\ra U'$ is finite $\square$

We can assume $V'$ and $U'$ to be smooth.

\end{lemma}

The relevant diagram is the following:

\centerline{
\xymatrix{ 
V' \ar@/^2pc/[rr]^-{n}\ar[r]^-{r'} \ar[dr]_-{s} & V"
\ar[d]^-{r"}\ar[r]^-{n"} & V \ar[d]^-{r}\\
 & U' \ar[r]^-{m} & U
}
}

{\bf Proof:} Just flatten the map $r$ by suitably modifying $U$ $\square$
 


To complete the proof of \ref {genfin}, we shall next need the following two properties {\bf (a), (b)}:

{\bf (a)} We then have (because $v"$ is finite):

(*a) $K_{Y"}+\Delta "=(v")^*(K_Y+\Delta)+R"$, an equality
between Weil Divisors, and $R"$ effective.

To see this,  observe we may first assume that $f'$ is deduced from 
$f$ by base change
by $v$ and smoothing of a main component (this is simply because 
$\Delta(f\circ u)\geq \Delta(f)$ if $Y'=Y$, and $u$ is generically finite. The argument is the same as in 
\ref{behmod}).

Next, a simple computation in local analytic coordinates over the generic  
point $y"$ of $Y"$ shows that:
if $r$ is the ramification order of $v"$ at $y'$ and $m$ is the 
multiplicity of the fibre of $f$ over $y:=v"(y")$, then the
multiplicity of the $f"$-fibre over $y"$ is $m':=(m/d)$, if 
$d:=gcd(r,m)$ and $r':=r/d$ 

For this, we need to choose $y"$ lying
outside the codimension two analytic subset of $Y"$ consisting
 of points $y"$ which are 
either singular on $Y"$, or
mapped by $v"$ to a point $y:=v"(y")$ either singular on 
$\Delta(f)\cup S$, or with $f$-fibre
of nongeneric multiplicity, where $S:=v"(S")$, with $S"$ the Weil 
divisor of points of $Y"$ at which  $v"$ ramifies.

 From this we see that near $y"$ on $Y"$, we have:

$K_{Y"}+\Delta"=(v")^*(K_Y+(1-1/r).R)+(1-1/m').(1/r).D$,
with $R:=v"(R^+)$ ($S"=R^+$ the
reduced ramification divisor of $v"$ near $y"$), and $D$ the (unique) reduced
component of $\Delta$ near $y$. If $R^+$ or $D$ is empty, our 
equality (*a) is obvious. Otherwise, $R=D$, by our choice of $y"$, and
(*a) follows from the inequality: $(1-1/m)\leq 
(1-1/m.r')=(1-1/m'.r)=(1-1/r)+(1-1/m').1/r$.\\

{\bf (b)} Moreover, we also have (quite generally): $K_{Y"}+\Delta 
"=(v')_{*}(K_{Y'}+\Delta ')$. This is an immediate check.

From the initial constructions made, we thus know that $f'_0$ is 
admissible, and dominates $f_0$, admissible (in the sense (**) just 
before the lemma).

We then get: $K_{Y"}+\Delta "$=
   $(v')_{*}(K_{Y'}+\Delta ')$=
   $(v")^{*}(K_Y+\Delta)+R"$.

Hence: $K_{Y'_0}+\Delta '_0$=
   $(w)_{*}(K_{Y"}+\Delta")$=
    $(w)_{*}((v")^*(K_{Y}+\Delta)+R"))$=
$(w_0)_{*}(v^*(K_Y+\Delta)+R'))$.

And so: $\kappa(Y',f')$=$\kappa(Y'_0,K_{Y'_0}+\Delta '_0)\geq$
$\kappa(Y'_0,(w_0)_{*}(v^*(K_Y+\Delta)))\geq$
$\kappa(Y',v^*(K_Y+\Delta))$=
$\kappa(Y,K_Y+\Delta)\geq\kappa(Y,f)$, as claimed in 
\ref{genfin}, the proof of which is thus complete $\square$\\

The last property asserted in Theorem \ref{behKod} is:

\begin{lemma}\label{etale}
Let the situation be as in Theorem
\ref{behKod} above. 
Assume now that
$u:X'\ra X$ is \'etale (ie: unramified). Then:
$\kappa(Y',f')=\kappa(Y,f)$.
\end{lemma}

{\bf Proof:} By Lemma \ref{genfin}, 
we need only to show that
$\kappa(Y',f')\leq \kappa(Y,f)$ for any $f,f'$ as above.

We can assume $f$ to be admissible, and still $u$ to be
\' etale: indeed, just replace $X'$ by $X'_1:=X'\times_{X}X_1$ , if 
$m:X_1\ra X$ is a modification having
an admissible fibration $f_1:X_1\ra Y_1$ equivalent to $f$. Then 
$X'_1$ is smooth, and $u_1:X'_1\ra X_1$ is \' etale.

The Stein factorisation $f'_1:X'_1\ra Y'_1$ of $f_1\circ u_1:X'_1\ra 
Y_1$ does not need to have $Y_1$ smooth. But a
modification $m':X'_2\ra X'_1$ of $X'_1$
allows to assume this to be true (accordingly modifying $Y_1$). Now 
the birational
invariance of the fundamental group for complex manifolds shows that
$X'_2$ may be assumed (after further modification) to be of the form 
$X'\times _{X} X_2$ for some modification $X_2$ of $X$.
Let us check this.

We can indeed assume that $u:X'_1:=X'\ra X$ is Galois, of group $G$, and 
that $f'$ is the (connected part of the) Stein factorisation of $(f\circ u)$.
The map $v:Y'\ra Y$ is thus $G$-equivariant as well.
But then, the modification $m':X'_2\ra X'_1$ can be assumed to be also 
$G$-equivariant, which shows that $X'_2$ is obtained from $u$ 
by base-change  over a modification of $X$. This establishes our 
preliminary claim (that $f$ can be chosen admissible).

Next, by Stein factorisation of $v=v"\circ v'$, 
we see as in the above proof of Lemma \ref{behKod}, 
that it is sufficient to show that:

(*) $K_{Y"}+\Delta "=(v")^*(K_Y+\Delta)$ 
(that is: $R"$ is empty).

Because then $K_{Y"}+\Delta "=(v')_*(K_{Y'}+\Delta ')$, and so
$K_{Y'}+\Delta '=(v)^*(K_Y+\Delta)+E'$, for some $v$-exceptional (not
necessarily effective)
$\mathbb Q$-divisor $E'$ of $Y'$, which implies the desired reverse
inequality.\\
By the finiteness of the map $v"$, the equality (*) has only to be
shown near any point $y\in Y$ lying
outside a codimension two subset $S$ of $Y$. We can thus assume that
$y$ is a generic point on some component of the support of
$\Delta(f)$.
Let us now cut Y by divisors in general position through $y$: we are
reduced to the case when $Y$ is a curve (the argument being local on
$Y$ in
the analytic topology, we don't need any algebraicity assumption of
$Y$). But then an easy local computation (see \cite{C98}, for example)
shows that for any
$y"\in Y"$, the order of ramification $r$ of $v"$ at $y"$ divides the
multiplicity $m$ of the fibre $X_y$ of $f$ over $y:=v"(y")$.
The multiplicity $m"$ of the fibre $X'_{y"}$ of $f'$ over $y"$ is
thus $m':=m/r$, since $u$ is \'etale.\\
We now compute $K_{Y"}+\Delta "$ near $y"$:
   $K_{Y"}+\Delta "=(v")^*(K_Y+(1-1/r)[y])+(1-1/m')[y"]$=
$(v")^*(K_Y+((1-1/r)+(1/r).(1-1/m))[y])$=
$(v")^*(K_Y+(1-1/m)[y])$, as claimed. (Here [y] is the reduced
divisor on the curve $Y$ supported by the point $y$) $\square$

\subsection{\bf The Sheaf of Differential Forms determined by a fibration}\label{diffib}

\begin{definition}\label{diffsheav} Let $X\in \mathcal C$ be smooth, and let $f:X\mero Y$ be a meromorphic fibration, with $Y$ 
reduced but not necessarily smooth . The rank one coherent
subsheaf $F_f$ of $\Omega ^p_X, p:=dim(Y)$ is defined as the saturation in $\Omega_X^p$ of 
$f^*(K_{Y_0}$, if $Y_0$ is the smoot locus of $Y$.
\end{definition}

Let us remark that the subsheaf $F_f$ so defined is a bimeromorphic invariant: 
it is preserved not only under modifications of $X$, but also under modifications of $Y$. 
In other words: $F_f$ depends only on the equivalence class of $f$.

We define by $\kappa(f)$ its Kodaira dimension.

Let us make more precise  what is understood by $\kappa(X,F_f)$: for 
$m>0$, define $H^0(X,F^m)$ to be the complex vector space of sections of
the subsheaf of $(\Omega_X^p)^{\otimes m}$ which coincides with 
$F_f^{\otimes m}$ over the Zariski open subset of $X$
(with codimension two or more complement) over which $F_f$ is locally 
free. Then define $\kappa$ in the usual way.

To see the usual property of $\kappa$ (being an integer, or 
$-\infty$), notice that the data are bimeromorphically invariant on $X$.
So that $F_f$ can be considered as the injective image of some locally 
free rank one sheaf $L$ on $X$ (after some suitable modification).
  The claimed property is then obvious, since it holds for $L$. We 
can then always implicitely
assume the existence of $L$ in the sequel, and also the holomorphic 
character of any meromorphic map defined on $X$ (such as, for example,
the ones defined by linear systems $L^{\otimes m}$).

We shall now describe $F_f$ in more detail.

\begin{definition} In the preceding situation, define
$F(f):=f^*(K_X)\otimes \mathcal O$$_X(\lceil f^*(\Delta(f))\rceil )$, 
where the symbol
used is the usual round-up (defined as $\lceil *\rceil :=-[-(*)]$,
applied to the coefficients of the irreducible components of the effective
$\Bbb Q$-divisor under consideration. Here [*] is the integral part, 
of course).
\end{definition}

\begin{definition}\label{partsup} Let $f:X\ra Y$ be a holomorphic fibration, and $S$ an effective 
divisor on $X$. We say that $S$ is {\bf partially 
suported on the fibres of f} if:
$f(S)\neq Y$, and if $T$ is any irreducible component of 
$f(S)$ of codimension one in $Y$, then $f^{-1}(T)$ contains an
irreducible component mapped on $T$ by $f$, but not contained in the support of $S$.
\end{definition}

Observe that if $S$ is partially supported on the fibres of $f$, so are its positive multiples. 

The introduction of this notion is due to the following:

\begin{lemma}\label{sectpartsupp} Let $f:X\ra Y$ be a fibration, and $S$ a divisor 
of $X$ partially supported on the fibres of $f$. Let $L$ be a line bundle on $Y$. 
The natural injection of sheaves $L\subset f_*(f^*(L)+S)$ is an isomorphism. 
\end{lemma}

{\bf Proof:} The assertion is of local nature on $Y$. So we can assume that $L$ is trivial. We then just need to show 
that $f_*(\cal O$$_X(S))\cong \cal O$$_Y$. We assume that $S\subset f^*(\cal O$$_Y(T))$, for some effective divisor $T$ 
on $Y$. Local sections of the the sheaf on the right are of the form $f^*(u/t)$, where $u$ is holomorphic on $Y$, while 
$t$ is a local equation of $T$. The sections of the sheaf on the left are meromorphic functions on $X$ of the same form, but 
with poles contained in $S$. Because $S$ is partially supported on the fibres of $f$, we get the claim ($t$ divides $u$) $\square$

We apply this to the following situation:

\begin{proposition}\label{F(f)} Let $f:X\ra Y$ be a holomorphic fibration between manifolds. 
There exists a Zariski closed subset $A\subset Y$ of codimension at least $2$ such that $F(f)+S$ and 
$F_f$ are naturally isomorphic over $(X-B):=f^{-1}(Y-A)$, where $S\subset X$ is an effective divisor partially supported 
on the fibres of $f$.
\end{proposition}

{\bf Proof:} Let $A$ be the union of the singular set of the support of $\Delta(f)$ and of 
images of all $f$-exceptional divisors on $X$. Let us remark that the above natural isomorphism 
is immediate outside of $\Delta$, because if $T\subset Y$ is a one-codimensional component of 
the locus of non-smooth fibres of $f$, then $f^{-1}(T)$ contains a reduced component at 
the generic point of which $f$ is smooth.

 So we consider the situation near a smooth point of some $\Delta_i$ not lying in $A$. 

In suitable local coordinates at the generic point of 
$D_{ij}$, in the notations of the lines preceding \ref{multdiv},
we have:

$(x)=(x_1,...,x_n)$, and $(y)=(y_1,...,y_p)$, with:
$f(x)=(y_1:=x_1^{m_{ij}},y_2:=x_2,...,y_p:=x_p)$.

And so: $\lceil f^*(K_Y+\Delta(f))\rceil$ is generated by
$x_1^{\lceil (m_{ij}/m_i)-1\rceil }.d(x')$, with 
$d(x'):=dx_1\wedge...\wedge dx_p$. A simple check shows that this is exactly the 
claim. (One may even observe that the divisor $S$ has the same description as the one 
given to define $F(f)$, by adding to $\Delta$ 
the one-codimensional components of the locus of non-smooth fibres of $f$) $\square$

\begin{corollary}\label{compsect} Let $f:X\ra Y$ be a fibration as in \ref{F(f)} above. Let $m>0$ 
be a sufficiently divisible integer. Then: 

(1) The natural isomorphism between $F(f)+S$ and $F_f$ over $(X-B)$ 
extends to a natural injection of $H^0(X,F_f^{\otimes m})$ into 
$H^0(X,m.(F(f)+S)) \cong H^0(Y,m.(K_Y+\Delta(f)))$.

(2) If $f$ is neat, this injection is bijective.
\end{corollary}

{\bf Proof:} (1) We start by observing that, by \ref{sectpartsupp}, the bijection:
 
$H^0(X,m.(F(f)+S)) \cong H^0(Y,m.(K_Y+\Delta(f)))$ actually holds. 

The natural map at the level of sections of $m$-th powers induces an isomorphism over a codimension 
two subset of $Y$. Because $m.(K_Y+\Delta(f))$ is locally free on $Y$, the said isomorphism thus extends as an injection, 
by Hartog's theorem. 

(2) This is because $B$ is mapped to a codimension two or more Zariski closed subset of $X'$ if $u:X\ra X'$ 
is a modification with $X'$ smooth and sending the $f$-exceptional divisors of $X$ in codimension $2$ or more in $X'$. 
Then the sections of $m.(F(f)+S)$ over $(X-B)$ extend to sections of $F_f^{\otimes m}$, as claimed $\square$

We now have the following easy but important consequence: 

\begin{proposition}\label{Fvsf} Let $f:X\ra Y$ be a holomorphic fibration, with 
$X,Y$ smooth. Then: 

(1) $\kappa(f)\leq \kappa(Y/\Delta(f)$

(2) $\kappa(f)=\kappa(Y/\Delta(f)$ if $f$ is neat (see \ref{neat} for this notion). 

(3) $\kappa(f)=\kappa(Y,f)$
\end{proposition}

{\bf Proof:} (1) and (2) are simply restatements of \ref{compsect}. We deduce (3) by choosing a neat
 admissible model $f'$ of $f$. Then: $\kappa(f)=\kappa(Y'/\Delta(f'))=\kappa(Y,f)$ $\square$

This allows us to give a short proof of the basic properties shown in Theorem \ref{behKod}:

\begin{corollary} The statements of Theorem \ref{behKod} hold.
\end{corollary}

{\bf Proof:} Using the notations there, we have a natural inclusion: $u^*(F_f)\subset F_{f'}$, 
which is an equality if $u$ is \' etale. the conclusions follow from the standard properties of 
the Kodaira dimension $\square$

Remark that one can also use the sheaves $F_f$ to simplify some of the geometric proofs given in section \ref{specdomgt}

\subsection{\bf Semi-continuity of the Kodaira Dimension}

\begin{proposition} \label{semicont}Let $f:X\ra Y$ and $g:Y\ra Z$
be holomorphic fibrations, with $X$ a connected compact 
complex manifold, and
$h:=g\circ f$. Let $Z^*$ be the
Zariski open subset of $Z$ over which $g$ and $h$ are smooth. Let, for
$z\in Z$, denote by $f_{z}:X_{z}\ra Y_{z}$ the restriction of $f$ to
the $z$-fiber $X_z$ of $h$, mapped
by $f$  to the $z$-fibre $Y_z$ of $g$. Then:
\begin{enumerate}

\item [1.] There exist modifications $\mu:X'\ra X$ and $\nu:Y'\ra Y$ such
that $f'_{z}$ is admissible, for $z$ general in $Z$, with $f':=\nu
^{-1}\circ f\circ \mu$ holomorphic.

\item [2.] Let $d$:=inf$(\lbrace \kappa(Y_{z},f_{z}), z\in Z^*\rbrace)$. (So
that $d\in \lbrace -\infty,0,1,\ldots,dim(Y)-dim(Z)\rbrace$). Then:
$d=\kappa(Y_{z},f_{z})$, for
$z$ general in $Z$.

\item[3] Let $A$ be the set of points $z$ of $Z^*$ such that $\kappa(Y_{z},f_{z})=dim(Y)-dim(Z)$. 
Then: either $A$ contains the general point of $Z$, or $A$ is contained in a countable union of closed proper 
analytic subsets of $Z$.

\end{enumerate}
\end{proposition}

{\bf Proof:} We don't mention modifications of $X$, since they don't
change the Kodaira dimensions of fibrations with $X$ as domain.
  For any modification $\nu:Y'\ra Y$ and
$d\in \lbrace -\infty,0,1,\ldots,dim(Y)-dim(Z)\rbrace$, let
$S^{*}_{d}(\nu)$:=$\lbrace z\in Z^{*}$ such
that: $\kappa(Y_{z},\Delta(f'_{z}))\geq d\rbrace$.
From \cite{Gr60}, and as in 5L-S], we deduce that
$S^{*}_{d}(\nu)$=$S_{d}(\nu)\cap Z^*$, where $S_{d}(\nu)$ is a
countable union of Zariski closed subsets of $Z$.

Obviously, $S_{d+1}(\nu)\subset S_{d}(\nu)$ for any $\nu$ and $d$
(this on $Z^*$ at least, which is sufficient for our purposes).

If $\nu':Y"\ra Y$ dominates $\nu$ in the sense that there exists a
$\nu":Y"\ra Y'$ with $\nu'=\nu"\circ \nu$, then obviously 
(by Theorem \ref{behKod}):
$S_{d}(\nu')\subset S_{d}(\nu)$.

We can thus define, for any $d$, $S_{d}\subset Z$ as the intersection
of all $S_{d}(\nu)$'s: it is again a countable intersection of Zariski
closed subsets of $Z$.

Define now $d$:=max$\lbrace d'$ such that $S_{d'}=Z\rbrace$. There
thus exists some $\nu$ such that $S_{d+1}(\nu)\neq Z$. Both claims then
follow immediately from the constructions just made.

The third assertion is an immediate consequence of the second.

\subsection{\bf Composition of Fibrations}\label{compfib}

This section will not be used before section \ref{orbadd}. 

Assume now $X,Y,Z$ to be smooth and $f:X\ra Y$ and $g:Y\ra Z$ 
to be holomorphic fibrations.

Our aim is to define, if $H$ is 
an orbifold structure on $Y$ (ie: an effective $\mathbb Q$-divisor on $Y$ 
with components having multiplicities of the form $(1-1/m)$, for $m$ integer), 
an orbifold structure $\Delta(g,H)$ on $Z$ in such a way that we have 
the equality: $\Delta(g,H)=\Delta(g\circ f)$ when $H=\Delta(f)$, if $f:X\ra Y$ 
is sufficiently ``high", in a sense defined in section \ref{orbadd} below.\\

Let us recall how we defined multiplicities: for any irreducible 
divisor $D\subset Z$, the multiplicity $m(g,D)$ of $g$ along $D$ is: 
$m(g,D)=:inf \lbrace m_{j},j\in J(g,D)\rbrace$.
Where $g^{+}(D):=(\sum_{j\in J(g,D)} m_{j}.D_{j})$ is the union
of all components (counted with their scheme-theoretic
multiplicities)  of $g^*(D)$ which are  mapped surjectively
onto $D$ by $g$, and $J(g,D)$ is the set of these components.\\

Let, if $D'\subset Y$ is an irreducible divisor of $Y$ mapped by $g$ to 
a divisor $D=g(D')\subset Z$, $m_*(g,D')$ be the multiplicity of $D'$ in $g^*(D)$.
Thus $m(g,D)=(\sum_{j\in J(g,D)}m_*(g,D_j).D_j)$.\\

We shall now define the notion of 
{\bf orbifold base of a fibration} $g: (Y/H)\ra Z$, when
the domain of the fibration is itself an orbifold $(Y/H)$. But $g:Y\ra 
Z$ is just a usual fibration, with $Z$ smooth.

Writing $H:=\sum_{i\in I} (1-1/m_i).H_i$, define, for any 
irreducible reduced divisor $D'\subset Y$
its muliplicity $m(H; D')$ in $H$ as being $m_i$ if $D'=H_i$,
and being $1$ otherwise (ie: if $D'$ is not a component of the 
support of $H$).

 For any irreducible divisor $D\subset Z$, with $Z$ 
smooth, define next the {\bf multiplicity $m(g,H;D)$ of $g: (Y/H)\ra Z$ along $D$} by:
$m(g,H; D):= inf_{j\in J(g,D)} \lbrace m_*(g,D_j).m(H;D_j)\rbrace$.

And, finally:

\begin{definition} \label{orbase} 

Let $g:Y\ra Z$ be a fibration, with $Z$ smooth.
Let $H$ be an orbifold structure on $Y$. We define the {\bf orbifold base
$\Delta(g,H)$ of the fibration $g:(Y/H)\ra Z$ } by:

$\Delta(g,H):=\sum_{D\subset Z} (1-1/m(g,H;D)).D$.
\end{definition}

In general, it is not true that $\Delta(g,H)=\Delta(g\circ f)$ if $H=\Delta(f)$. 
But the following results at least are available:

\begin{proposition} Let $f:X\ra Y$ and $g:Y\ra Z$ be two 
holomorphic fibrations, with $X,Y,Z$ smooth.
Then: $\Delta(g\circ f)\leq\Delta(g,\Delta(f))$.

(Recall that, for two $\mathbb Q$-divisors $A,B$ on a variety $Z$, we 
write $A\leq B$ if $(B-A)$ is effective.)
\end{proposition}

{\bf Proof:} We easily check that we have, for any prime divisor $D\subset Z$:

$m(g,\Delta (f);D)=inf_{j\in J(g,D)} \lbrace m_*(g,D_j).m(\Delta(f);D_j)\rbrace$, by definition; and:

$m(g\circ f,D)=inf_{k\in J(g\circ f,D)} \lbrace m(g\circ f,D_k)\rbrace$, also by definition. 

Now, if $f(D_k)=D_j$ is a divisor, then: $m_*(g\circ f,D_k)=m_*(g,D_j).m(\Delta(f);D_j)$, by an easy check.

Observe that the minimum of these values, taken over $J(g,D)$, is precisely $m(g,\Delta(f);D)$.

Thus:  $m(g\circ f,D)=inf \lbrace m(g,\Delta(f);D), m(f,g)\rbrace$, where: 
$m(f,g):=inf_{k\in J'(g\circ f,D)} \lbrace m_*(g\circ f,D_k)\rbrace$, 
and where $J'(g\circ f,D)$ is the set of irreducible components of
 $(g\circ f)^*(D)$ which are surjectively mapped to 
$D$ by $(g\circ f)$, but are $f$-exceptional. From this the claim follows.\\

Observe that strict inequality may thus occur. However:

\begin{proposition} 
\label{compofib}
Let $f,g$ be two fibrations, as above.
Let $f':X'\ra Y'$ be a modification of $f$ (with modifications $u:X'\ra X$ and $v:Y'\ra Y$ such that: 
$f\circ u=v\circ f'$). Then: 
\begin{enumerate}
\item[1.] $\Delta(g'\circ f')=\Delta(g\circ f)$, with $g':=g\circ v$.
\item[2.] $\Delta(g',\Delta(f'))\leq \Delta(g,\Delta(f))$
\item[3.] There exists a modification $f'$ of $f$ such that: 
$\Delta(g',\Delta(f'))\leq \Delta(g",\Delta(f"))$, for any modification $f"$ of $f'$.
\item[4.] If $f'$ is a modification as in (3) above, then: $\Delta(g'\circ f')=\Delta(g',\Delta(f'))$.
\item[5.] In (3) above, we can choose $f'$ such that, moreover, $f',g'$ and $(g'\circ f')$ are prepared, admissible and high.
\end{enumerate}
\end{proposition}

{\bf Proof:}

 (1) $g'\circ f'=(g\circ f)\circ u$. So, the claim follows from \ref{behKod}.

(2) We can and shall thus assume that $X'=X$ and $f=v\circ f$. Let $D\subset Z$ be a prime divisor. Then $J(g,D)\subset J(g',D)$, the 
difference consisting of the $v$-exceptional components of $(g')^*(D)$. Moreover, for each $j\in J(g,D)$, we have: 
$m(g,D_j)=m(g',\overline {D_j})$, with $\overline {D_j}$ the strict transform of $D_j$ by $v$. Finally:
 $m(f,D_j)=m(f',\overline{D_j})$, by 
\ref{behKod}. This implies the claim.

(3) This is simply because there are only finitely many orbifold divisors on $Z$ lying between $\Delta(g\circ f)$ and 
$\Delta(g,\Delta(f))$, and because of (1) and (2) above, which show that the first (resp. second) term is invariant (resp. decreases) 
under a modification.

(4) Assume indeed that $f$ (rather than $f'$ to ease notations) 
enjoys the property stated in (3). Assume we have $\Delta(g\circ f)<\Delta(g,\Delta(f))$. 
This means that there exists an $f$-exceptional prime divisor $D'\subset X$ such that $g(f(D')):=D\subset Z$ is a divisor 
in $Z$, and that the multiplicity $m'$ of $D'$ in $(g\circ f)^*(D)$ is equal to $m(g\circ f,D)$, and is so 
strictly less than $m(g,\Delta(f));D)$. 

Take a modification $f'$ of $f$ such that the strict transform of $D'$ in $X'$ is no longer $f'$-exceptional. The multiplicity 
of $D"$, the strict transform of $D'$ in $X'$ is then $m'$ (by an easy check). Thus:
 $m(g', \Delta(f');D)\leq m'<m(g,\Delta(f);D)$. But this contradicts the property (3), supposed to hold for $f$. Contradiction; the claim follows. 

(5) By modifying $Z$, we can assume that $(g\circ f)$ and $g$ are admissible, high, and moreover that the non-smooth loci of these 
two fibrations are contained in a normal crossings divisor. By modifying next $Y$, we can assume that $g$ is prepared, $f$ admissible and high, 
and that the non-smooth locus of $f$ is contained in a normal crossings divisor of $Y$. 
Finally modify $X$ to get the remaining stated properties $\square$

\section{\bf Special Fibrations and General Type Fibrations}

\subsection{\bf Special or General Type Fibrations}

\begin{definition} Let $f:X\mero Y$ be a fibration, with $X,Y$ compact irreducible.

\begin{enumerate}
\item[1.] The fibration $f:X\mero Y$ is said to be of
{\bf general type} if $\kappa(Y,f)=dim(Y)>0$.

\item[2.] The variety $X$ is said
to be {\bf special} if $X\in \mathcal C$ and if there is no meromorphic 
fibration $f:X\mero Y$
of general type, for any $Y$.

\item[3.] The fibration $f:X\mero Y$ is said to be 
{\bf special} if $X\in \mathcal C$, and if its general fibre is special.
\end{enumerate}

\end{definition}

Recall that here $\mathcal C$ is the class of compact complex spaces $X$
  which are bimeromorphic to
(or, equivalently: dominated by ) some compact K\"ahler manifold $X'$ 
(depending on $X$).
This class was introduced by A.Fujiki.

Recall (see \cite{C81}) also that a point of a complex space $Y$ is said to be {\bf general} if it lies 
outside of a countable union of closed analytic subsets of $Y$, none of which containing any irreducible component of $Y$.
Similarly, if $f:X\ra Y$ is a fibration, one of its fibres $X_y$ is general if it lies above a general point $y$ of $Y$.

\begin{example}\label{specex}
\end{example}

We list (proofs need tools developed below, and so are given later) some examples of special manifolds.

\begin{enumerate}
\item[0.] A variety of general type (and positive dimension) is {\bf not} special. 
(Consider its identity map: it is a fibration of general type). 
\item[1.] A {\bf curve} is special iff its genus is 0 or
1, iff its Kodaira dimension is at most zero, iff its fundamental group
is abelian, iff it is not hyperbolic. This because a curve has only 
the two trivial fibrations (constant, and identity).

The two fondamental examples of special Manifolds are direct generalisations:

\item[2.] A manifold which is {\bf rationally connected} is special. See Theorem 
\ref{RCspec} for a geometric proof and definition 
of the notions involved. Another shorter (but more abstract) proof of the 
specialness of rationally connected manifolds is given
 \ref{Ktriv} below.

\item[3.] A manifold $X$ with {\bf vanishing Kodaira dimension} 
(ie: $\kappa(X)=0$) is special. See Theorem \ref{k=0} for the proof. 
Another proof is given in \ref{Ktriv} in the special case where $c_1(X)=0$.

\item[4.] More generally, special manifolds are conjectured to be build up 
from manifolds either rationally connected, or 
with Kodaira Dimension zero by suitable compositions of fibrations with 
fibres of these two types. See section \ref{c=rj^n} for a 
precise formulation.

\item[5.] For any $d>0$ and $k\in\lbrace -\infty,0,...,d-1\rbrace$, there exists 
projective manifolds of dimension $d$ and 
Kodaira dimension $k$. See \ref{d,k,spec} for such examples.

\item[6.] A manifold $X\in \mathcal C$ is special if there exists a nondegenerate meromorphic map 
from $\mathbb C^n$ to $X$, where nondegenerate means: submersive at some point where it is holomorphic. 
For example, a complex torus, or a projective space are special 
(this follows also from [2,3] above as well).
See theorem \ref{ko} for the proof of a more general version.

\item[7] K\" ahler manifolds with nef anticanonical bundle are conjectured to be special. 
This conjecture implies 
most usual conjectures concerning these manifolds. See [D-P-S],[Zh],[Pa].

\item[8.] A manifold X of {\bf algebraic dimension zero} (denoted $a(X)=0$, to mean
that all meromorphic functions on $X$
are constant, so that meromorphic maps from $X$ to projective 
varieties are constant) is also special (simply because any meromorphic map from $X$ 
onto a projective manifold is constant).

 \item[9.] More generally (see [U],Chap. 12 for
the notions used):

\begin{theorem}\label{fibalgred}
 Let $a_{X}:X\ra Alg(X)$ be the
algebraic reduction of $X\in \mathcal C$. The generic fibre of $a_X$ is
special.

\end{theorem}

Recall that the algebraic dimension of $X$, denoted $a(X)$, is the
dimension of $Alg(X)$, and also the transcendance degree over
$\mathbb C$ of the
field of meromorphic functions on $X$. One says that $X$ is {\bf
Moishezon} if $a(X)=dim(X)$. This also means that $X$ has a
modification which is projective.

The proof of the preceding result is given in \ref{fibAlgred} below.

\end{enumerate}

\begin{question} \label{specquest} Two important stability properties of 
the class of special manifolds are expected to hold, but 
are not proved in the present paper: the stability under deformation and specialisation. 
Do these hold?
\end{question}

\subsection{\bf Special Fibrations Dominate General Type Fibrations}\label{specdomgt}

The geometric study of special manifolds is based on Theorem \ref{domin}, 
stated and proved below. Its proof rests on several preliminary 
resuts of independent interest that we now give.

\begin{lemma}\label{specsurj} Let $g:X'\ra X$ be meromorphic
surjective (ie: dominant). Assume that $X'$ is special. Then $X$ is 
special, too.
\end{lemma}

{\bf Proof:} Assume first that $g$ is connected (ie: a
fibration). Let, if any, $f:X\ra Y$ be a fibration of general type. We
can assume that
$f\circ g:X'\ra Y$ is admissible. Then obviously,
$\Delta(f\circ g)=\Delta(f)+E$, for some effective divisor $E\subset 
Y$ (because $m(f\circ g, D)\geq m(f,D)$, for any irreducible divisor $D \subset Y$).
Thus $f\circ g$ is of general type, too.
   A contradiction. No such $f$ does exist, which is what was claimed.

In the general case, Stein factorise $g$ and use the first part to
reduce to the case where $g$ is generically finite. If $f$ as in the
first part
exists, then we deduce from \ref{behKod} that (the fibration part of) the
Stein factorisation of $f\circ g$ also is of general
type. Hence again a contradiction $\square$

\begin{proposition}\label{gtsect'} Let $f:X\ra Y$ be a fibration of
general type. Let $j:Z\ra X$ be holomorphic such that $f\circ j:Z\ra Y$ is
surjective.
Let $f\circ j=g\circ h$ be the Stein factorisation of $f\circ j$, 
with $h:Z\ra Y'$
connected and $g:Y'\ra Y$ finite. Then: $h$ is a fibration of general type.
\end{proposition}

 In particular, if $dim(Z)=dim(Y)$, we get:

\begin{corollary}\label{gtsect} Let  $f:X\ra Y$ be a fibration of general type.
Let $j:Z\ra X$ be meromorphic, and such that  $f\circ j:Z\ra Y$ surjective.
Then $Z$ is a variety of general type.
\end{corollary}

{\bf Proof (of \ref{gtsect'} and \ref{gtsect}):} Assume first that $f\circ j:Z\ra Y$ is connected and
admissible (as we can, then). For any component $\Delta_i$ of
$\Delta:=\Delta(f)$, we have (restricting to components surjectively 
mapped onto
$\Delta_i$ by $(f\circ j)$):
$(f\circ 
j)^{*}(\Delta_i)$=$j^{*}(f^{*}(\Delta_i))$=$j^{*}(m_i.D_i)$=$m_i.j^{*}(D_i)$.

Thus: $\Delta(f\circ j)=\Delta(f)+E$,
for some effective $\mathbb Q$-divisor $E$
of $Y$. And $(f\circ j)$ is thus of general type in this case.

We now consider the general case: define $f': X':=
\widetilde{(X\times_{Y} Y')}\ra Y'$ be deduced from the base change
by $g$ and smoothing of the
fibre product. The map $j$ lifts meromorphically to $j':Z\ra Y'$, by
construction, because $f\circ j$ is surjective. But now $f'\circ j':Z\ra Y'$ is
a fibration.
Applying the first part, we get the claim.

For \ref{gtsect}, notice that in this situation, $h$ is bimeromorphic and of
general type by \ref{gtsect'}. Thus $Z$ is itself of general type, as claimed $\square$

\begin{example} We can now give two elementary examples of special varieties:

\label{spec1}

\begin{enumerate}

\item[1.]  $\mathbb P_n(\mathbb C)$ is special.

\item[2.] A product of special varieties is special.

\end{enumerate}
\end{example}

{\bf Proof:} Indeed (for (1)): let $f:\mathbb P_n(\mathbb C)\mero Y$ be any general
type fibration, if any. Let $m:=dim(Y)>0$.
   Choose $j:Z=\mathbb P_m(\mathbb C)\subset \mathbb P_n(\mathbb C)$
such that $f\circ j$ is surjective to contradict \ref{gtsect'}.

The proof of (2) is similar. (We shall prove more general results in section \ref{core}).

\begin{proposition}
\label{restgt} 

Let $f:X\ra Y$ and $k:Y\ra W$ be fibrations.
Assume that $f$ is of general type.
Then $f_w:X_w\ra Y_w$ is also of general type, for $w\in W$ general.
\end{proposition}

{\bf Proof:}  Let $w\in W$ be general, and recall that
$f_w:X_w\ra Y_w$  is nothing, but the restriction of $f$ to $X_w$.
But then,
   $\Delta(f_w)=\Delta(f)_{\mid Y_{w}}+E_{w}$, with $E_w$ effective and
empty for generic $w$ in $W$. Moreover, $K_{Y_{w}}= K_{Y\mid Y_{w}}$, 
by adjunction.

   Thus $K_{Y_{w}}+\Delta(f_w)= (K_{Y}+\Delta(f))_{\mid Y_{w}}$ for general $w$.
Now $(K_{Y}+\Delta(f))$ is big. Thus so is its restriction to
$Y_{w}$. By modifying adequately
$X$ and $Y$, we can assume that $f_{w}$ is admissible by the following 
lemma. We thus get the claimed property $\square$

\begin{lemma} Let $f:X\ra Y$ and $g:Y\ra Z$ be
fibrations. There exists representatives of $f$ and $g$ (also denoted 
$f$ and $g$) such
that $f_z$ is admissible, for $z\in Z$ general.
\end{lemma}

{\bf Proof:} This is a special case of \ref{semicont} $\square$

The main result in this section is:

\begin{theorem}\label{domin}
 Let $h:V\mero Z$ and $f:X\mero Y$ be
fibrations with $f$ of general type and $h$ having general fibres
which are special.
Let $g:V\mero X$ be meromorphic surjective. Then, there exists $k:Z\ra
Y$ such that $f\circ g=k\circ h$.

\end{theorem}

The situation is described in the following commutative diagram:

\centerline{
\xymatrix{ 
V \ar[r]^{g}\ar[d]_-{h=special} & X\ar[d]^-{f=gen.type}\\  
Z  \ar[r]^-{k} & Y\\ 
}
}

The special case $V=X$ deserves special mention:

\medskip

\begin{theorem}\label{domin'} Let $h:X\mero Z$ and $f:X\ra Y$ be
fibrations with $f$ of general type and $h$ {\bf special}. Then, there exists $k:Z\ra
Y$ such that $f\circ g=k\circ h$ (We say that $h$ {\bf dominates} $f$).

\end{theorem}

The corresponding diagram is:

\centerline{
\xymatrix{ 
X \ar[dr]^{f=gen.type}\ar[d]_-{h=special} \\  
Z  \ar[r]^-{k} & Y\\ 
}
}

\begin{remark} The special case where $X=Y$ is of general type is obvious, by the easy addition theorem, 
because the covering family of $Y$ by the subvarieties $h(X_z)$ has a generic member of general type. These 
subvarieties must be points.

\end{remark}

{\bf Proof (of \ref{domin}):} This is a direct consequence of \ref{domin+} below. Indeed, if
such a map $k$ does not exist, then $f\circ g (V_z)$ is
positive-dimensional, for generic $z\in Z$. But the Stein
factorisation of $f\circ g_z$ is then of general type, by \ref{domin+}.
This contradicts the
assumption that $g$ is special $\square$ 

\begin{proposition}
\label{domin+} 
Let $f:X\ra Y$ and $h:V\ra Z$ be
fibrations. Assume $f$ is of general type. Let $g:V\ra X$ be a 
surjective meromorphic map.
  Let $g_z:V_z\ra Y_z$ , be
the
   restriction of $f\circ g$ to $V_z$, with $Y_z:=h_z(V_z)$.

Assume $Y_z$ is positive dimensional. Then: the Stein factorisation
of $g_z$ is of general type for $z$ general in $Z$.

\end{proposition}

{\bf Proof:} By \ref{restgt}, we can then replace $X$ by $V$
   and $f$ by $(f\circ g)$ without
losing the hypothesis that $(g\circ g)$ is of general type both when $g$ is
generically finite or is a fibration. We thus see that
the (connected part of the) Stein factorisation $f'$ of $f\circ g$ is a fibration of general
type. Replace $X$ by $V$ and $f$ by $f'$, so that we are
reduced
to the case where $X=V$, and $f=(f\circ g)$, which we now treat.

We then notice that we can replace $Z$ by any subvariety $Z'\subset
Z$ going through a general point of $Z$, 
and $X$ by $X'\subset X$, defined by $X':=g^{-1}(Z')$, provided $f(X')=Y$.
This is because of \ref{restgt}, which shows that the (Stein factorisation
of the) restriction of $f$ to $X'$ is still of general type.

We shall then construct an appropriate $Z'\subset Z$.

Let $c:Z\ra \mathcal C$$(Y)$ be the meromorphic map sending a generic
$z\in Z$ to the reduced cycle of $Y$ supported on $Y_z$. Here $\mathcal
C$$(Y)$
denotes the Chow Scheme of $Y$. Observe next that $f$ being of
general type, $Y$ is Moishezon. By modifying suitably $Y$, we shall
assume that
$Y$ is projective. Let $W\in \mathcal C$$(X)$ be the image of $Z$ by $c$.
Thus $W$ is projective, too.

   We next choose $W'\subset W$ to be an intersection of generic
members of any very ample linear system on $W$, in such a way that
the incidence graph
$Y'\subset W'\times Y$ of the algebraic family of cycles of $Y$
parametrised by $W'$ is generically finite over $Y$. This means that
if $p:Y'\ra Y$
and $q:Y'\ra W'$ are induced by the natural projections, then $p$ is
generically finite surjective. More concretely, the generic point of
$Y$ is contained in
only finitely many of the $Y_w$'s, for $w\in W'$. Define now
$Z':=c^{-1}(W')$. (Remark that when $Z$ is Moishezon,
we don't need to consider $c$, and can just take intersections of ample
divisors of a projective modification of $Z$ directly).

We now replace $X,Z,g,f$ respectively by $X',Z'$, and their
restrictions to $X'$.

\begin{lemma} $p:Y'\ra Y$ is bimeromorphic.

\end{lemma}

{\bf Proof:} Let $c'$ be the restriction of the above map $c$ to $Z'$.
Then $Y'\subset (W'\times Y)$ is the image of $(c'\circ h)\times 
f:X'\ra W'\times Y$.
Thus $f:X'\ra Y$ lifts to $f':X'\ra Y'$ such that $p\circ f'=f$. 
Because $p$ is generically
finite and $f$ connected, we see that $p$ is bimeromorphic.
$\square$

We have, by construction, $q\circ f=c'\circ h:X\ra W'$. We now can conclude
by applying successively \ref{restgt} and \ref{gtsect} to $X_w$, for $w$ general
in $W'$. Indeed: from \ref{restgt}
we learn that the restriction $f_{w}:X_{w}\ra Y_w$ of $f$ to
$X_w$ is of general type. Further, for $z$ generic in
${c'}^{-1}(w):=Z_w'$, we know that
$f(X_{z})=f(X_w)=Y_w$, and so  by \ref{gtsect}, the Stein factorisation of
the restriction of $f$ to $X_z$ is of general type, as claimed $\square$

\begin{corollary} Let $f:X\mero Y$ be a special
fibration, and let $j:Z\ra X$ be such that $f\circ j:Z\mero Y$ is onto.
Assume that $Z$ is special. Then $X$ is special.
\end{corollary}

{\bf Proof:} Let $h:X\mero W$ be a fibration of general type (if any).
By \ref{domin} with $V=X$, there is a factorisation $k:Y\ra W$ such that
$h=k\circ f$. Thus: $j\circ h:Z\ra W$ is onto. We can thus apply \ref{gtsect},
which says that
the Stein factorisation $h':Z\ra W'$ of $h\circ j:Z\ra W$ is of general
type. But $Z$ being special by assumption, this
is a contradiction and $X$ is special.

\begin{remark} It is not true true in general that $X$ is
special if it admits a special fibration $f:X\ra Y$ with $Y$ special
(see \ref{amotexamp}, for example). But in some cases (if the fibres are for example,
rationally connected), this is true (see \ref{RCmultfree}).
\end{remark}

\begin{example}\label{d,k,spec}
For any $d>0$ and $k\in\lbrace -\infty,0,...,d-1\rbrace$, there exists projective manifolds of dimension $d$ and 
Kodaira dimension $k$. In particular, it is not true that a special manifold has
nonpositive Kodaira dimension.
\end{example}
To get examples with $k\geq 0$, just take indeed a general member of the linear
system:

 $|\mathcal{O}$$_P(d-k+2,m)|$
on $P:=\Bbb P^{d-k+1}\times \Bbb P^k$ for large m, and
such that this member has a section over
the base $\Bbb P^k$. 

If $k=-\infty$, just take $\Bbb P^d$.\\

We now come to a very important consequence of \ref{domin'}

\begin{corollary} Let $X\in \cal C$. There is {\bf at most} one fibration defined on $X$ which is both 
special and of general type. If it exists, this fibration special and is 
dominated by any other special fibration defined on $X$. 
It is also, moreover, of general type and dominates any other general type fibration defined on $X$.
\end{corollary}

The proof is immediate, from \ref{domin'}. In other words, such a fibration is the ``lowest special" and the 
``highest of general type" on $X$. The existence of such a fibration on any $X$ is the main result of tis paper. 

These two descriptions provide us with two means of construction: by consideration of 
chains of special subvarieties, one geometrically constructs the ``lowest special" fibration on any $X$. This is the way 
used in section \ref{core}. Dually, by making fibre products of general type fibrations, one constructs the 
``highest general type" fibration on $X$. This is the approach followed in  section \ref{geomcons}. In both cases, to 
show that the fibration so constructed has the missing property (special if general type, and conversely), we need the orbifold 
additivity result \ref{cnmtheor}.

Recall now:

\begin{definition}\label{defalmhol}
 Let $f:X\ra Y$ be a surjective
meromorphic map between normal compact irreducible analytic complex spaces.

We say that $f$ is {\bf almost holomorphic} 
if $f(J)\neq Y$, where
$J$ is the indeterminacy locus of $f$. 
\end{definition}

More precisely: if $X'\subset
X\times Y$
is the graph of $f$, and $f':X'\ra Y$ the (restriction of the) second
projection, then $f(J):=f'(J')$, with $J'$ being the set of all $x'\in
X'$ such that
$p^{-1}(x)$ does not reduce to $x'$, or equivalently: is positive 
dimensional. (Here $p:X'\ra X$ is the first
projection (which is a proper modification), and $x:=p(x')$).

\begin{theorem}\label{gtalmhol} Let $f:X\ra Y$ be a meromorphic fibration 
of general type, with $X\in \mathcal C$ smooth.
Then $f$ is almost holomorphic.
\end{theorem}

\begin{remark} The smoothness assumption is essential, as
shown by the example of the cone $X$
   over a projective manifold of general type $Y$.
The conclusion  of the preceding theorem should however hold if the singularities of $X$ 
are terminal, or even canonical.

\end{remark}

{\bf Proof:} Resolve the indeterminacies of $f$ by a
sequence of smooth blow-ups $u:X'\ra X$ (by H.Hironaka's results),
with $f':=f\circ u:X'\ra Y$ holomorphic. If $f$ is not almost holomorphic,
some irreducible component $V$ of the exceptional divisor of $u$
is mapped surjectively onto $Y$ by $f'$, in such a way that the
fibres $V_z$ of the restriction $u'$ of $u$ to $V$
are mapped to positive-dimensional subvarieties $f(V_z)$ of $Y$. This
contradicts \ref{domin+}, because $V$ has two maps $u':V\ra Z:=u(D)\subset
X$,
and $f":V\ra Y$, the restriction of $f'$ to $V$. Now, by smoothness of
$X$, $u'$ is special (its generic fibre is a rational variety, so
apply \ref{spec1}. Moreover, by \ref{gtsect}, $f"$ is of general type 
(possibly after Stein factorisation).

  From \ref{domin+}, we thus have a factorisation $\phi:Z\ra Y$ with $f"=\phi\circ u'$.
   But this precisely contradicts $dim(f(V_z))>0$, and we get the claim $\square$

\subsection{\bf General Type Fibrations and Bogomolov Sheaves}

\begin{definition}\label{bogsheav} Let $X\in \mathcal C$. A rank one coherent
subsheaf $F$ of $\Omega ^p_X, p>0$ is said to be
a ($p$-dimensional) {\bf Bogomolov sheaf} on $X$ if $\kappa(X,F)=p$.
\end{definition}

The properties of these Kodaira dimensions have been discussed in section \ref{diffib} to which we thus refer. 

By the results of that section, any (equivalence class of a) general type fibration $f$ defined on $X$ 
canonically determines a Bogomolov sheaf $F_f$ on $X$. We shall now see the converse direction.

By the results of [Bo], any $p$-dimensional Bogomolov sheaf determines a meromorphic 
fibration $f_F:X\ra Y_F$ with $dim(Y_F)=p$,
and such that
$F=f_F^*(K_{Y_F})$ at the generic point of $Y_F$. 

The proof given there applies 
only to $X$ projective
(because of the argument of cutting by transversal hyperplane sections),
but (as is well-known) can be easily modified to apply to any $X$ compact K\" ahler (or in 
$\mathcal C$), as follows:

\begin{theorem} \cite{Bo} \label{Bog} Let $F\subset \Omega_X^p, p>0$,
be a Bogomolov sheaf on $X\in \mathcal C$. Let $f_F:X\ra Y_F$ be the
fibration defined by the linear system $\mid L^{\otimes m} \mid$ for $m>0$
sufficiently large and divisible. We can assume that
$f$ is holomorphic. Then $F=f_F^*(K_{Y_F})$ at the generic point of $Y_F$.
\end{theorem}

{\bf Proof:} We can, using the covering trick argument of \cite{Bo},
reduce to the case when $m=1$, which we now treat.

We can thus select $(p+1)$ sections $s_i, i=0,1,...,p$ of $F$ which 
are analytically independent
(ie: the linear system they define is $f_F$ up to stein factorisation, and so
has $p$-dimensional image). Because $F$ has rank one, there exists meromorphic
functions $y_i, i=1,...,p$ such that $s_i=y_i.s_0$.

By Hodge theory ($X$ being K\" ahler, or even just in $\mathcal C$), the
holomorphic $p$-forms $s_i, i=0,...,p$ on $X$ are closed. From which
we get: d$s_0$=d$y_i\wedge s_0=0, i=1,...,p$.

  The last equality shows by simple algebraic arguments the existence of a
meromorphic function $g$ on $X$ such that $s_0=g.(dy_1\wedge ...\wedge dy_p)$.
The first equality shows that $g=f^*(h)$, for some meromorphic 
function $h$ on $Y$,
and so the claim, since the argument applies to $i>0$ as well $\square$

We can thus sum up the preceding observations as  follows:

\begin{theorem}
\label{gtBog}
Notations being as above, for any $X\in 
\mathcal C$, there are inverse bijective correspondances between Bogomolov sheaves $F$ on $X$ and 
(eqivalence classes of) general type fibrations $f$ defined on $X$. These correspondances are defined as fillows:
\begin{enumerate}

 \item[1.]If $f$ is of general type, then $F_f$ is a Bogomolov sheaf on $X$ .

 \item[2.]If $F$ is a Bogomolov sheaf on $X$, then $f_F$ is a 
fibration of general type.
\end{enumerate}
\end{theorem}

A direct application (and motivation) is:

\begin{theorem}
\label{specnobog} 
The manifold $X\in\mathcal C$ is special if and only if there is 
no Bogomolov sheaf on $X$. 
\end{theorem}

{\bf Proof:} The Bogomolov subsheaves on $X$ correspond bijectively to fibrations 
of general type with domain $X$ $\square$

\begin{corollary}\label{Ktriv} The manifold $X$ is special in the following two cases:

\begin{enumerate}
\item[1.] $X$ is rationally connected (see section \ref{rc} for this notion).

\item[2.] $X$ is a compact K\" ahler manifold with trivial (or torsion) canonical bundle, or even, more generally, if 
$c_1(X)=0$.
\end{enumerate}
\end{corollary}

{\bf Proof:} In both cases, it is shown in [C95] that $\kappa^+(X)\leq 0$, which means (in particular) 
that a coherent rank one subsheaf of $\Omega_X^p,p>0$ has Kodaira dimension negative or zero. Thus $X$ has no 
Bogomolov subsheaf. It is thus special, by \ref{specnobog}. 
The result in [C95] depends on Calabi-Yau's Theorem. But in the projective case, one can get algebro-geometric 
proofs using Miyaoka's generic semi-positivity Theorem $\square$

We shall see later that the weaker condition $\kappa(X)=0$ is actually sufficient for $X$ to be special. 

Notice that the property shown in [C95] in the above two cases is considerably stronger 
than the absence of Bogomolov sheaves. This is not surprising, in view of the fact that these manifolds are the 
building blocks of the class of special manifolds, but do not exhaust this class, by far.

\subsection{\bf  General Type Reduction}

 \subsubsection{\bf Ordering of Fibrations} \label{ordfib}

Recall from \ref{fib} that a  meromorphic fibration 
$f:X\mero Y$ canonically defines (see [C85]) a meromorphic
map $\phi_f:Y\ra {\cal {C}}(X)$,
with $\cal C$$(X)$ the Chow-Scheme of $X$, by sending the generic 
$y\in Y$ to the point of
$\cal C$$(X)$ parametrising the reduced fibre $X_y$ of $f$ over $y$. 
Let $\Phi_f\subset \cal C$$(X)$
be the image of $Y$ by $\phi_f$: it is a compact irreducible analytic 
subset of $\cal C$$(X)$ bimeromorphic
to $Y$ such that its incidence graph is bimeromorphic to $Y$.

The above correspondance induces a bijective map between equivalence 
classes of fibrations and
compact irreducible analytic subsets of $\cal C$$(X)$ with incidence 
graph bimeromorphic to $X$.

Recall also we said that $f$ is {\bf almost holomorphic} if the image by $f$
(naturally defined using the graph of $f$)
of the indeterminacy locus $I(f)$ of $f$ is not all of $Y$ (see \ref{defalmhol} for details). It is 
easy to show that $\Phi_f$ is
an irreducible component of $\cal C$$(X)$ if $f$ is almost 
holomorphic (see [C85]).

We now introduce an order on the set of (equivalence classes of) fibrations with domain $X$.

We say that $f$ {\bf dominates} the fibration $g:X\ra Z$ if there 
exists a meromorhic fibration
$\phi: Y\mero Z$ such that $g=\phi\circ f$. Equivalently: each fibre of 
$f$ is contained in some fibre of $g$.
  This binary relation defines an ordering on the set $\cal F$$(X)$ of 
all equivalence classes of fibrations
(seen as a subset of $\cal C$$(\cal C$$(X))$).

There is now an easy 

\begin{lemma}\label{uppfib}

If $\Lambda\subset \cal F$$(X)$ is any 
subset, it has in the ordered set
$\cal F$$(X)$ a least
upper bound, denoted $\Lambda^+$.

\end{lemma} 

 {\bf Proof:} $\Lambda^+$ is so constructed: let 
$\Lambda_0:=\lbrace \lambda _1,...,\lambda _N \rbrace\subset \Lambda$ 
be
finite such that the product map $f:=f_{\lambda _1}\times...\times 
f_{\lambda _N}$ has an image of maximal dimension. Then
take for $\Lambda ^+$ the (fibration part of the) Stein factorisation of $f$.
$\square$

Finally, the preceding construction shows that if any element of 
$\Lambda$ is almost holomorphic, so is the least upper bound 
$\Lambda^+$ of the family.

\begin{example}\label{fibprodgt} If $\Lambda$ consists of fibrations onto 
varieties of general type, then $\Lambda^+$ is also a fibration onto a variety of general type. 
\end{example}

This is easily reduced to the case when $\Lambda$ has two elements, and then reduces to 
showing that if $Z\subset Y\times Y'$ is a subvariety of a product of two varieties of general type, 
then $Z$ itself is of general type if it is mapped surjectively to $Y,Y'$ by the 
first and second projections. This results easily from the additivity theorem for fibrations 
with base of general type (generalised orbifold versions will be proved in section \ref{orbadd} below) $\square$

\begin{definition}\label{defgtred}

 For any  $X\in \mathcal C$, let $gt_X:X\mero GT(X)$ be
the least upper bound in $\mathcal F$$(X)$ of the family $\Lambda _X$ of 
all (equivalence classes of)
fibrations of general type $f_j:X\mero Y_j$. 

(if $X$ is special, we
just take for $gt_X$ the constant fibration). 

We call $gt_X$ the {\bf general type reduction of $X$}.

\end{definition}

From \ref{ordfib} above, we deduce:

\begin{proposition} \label{gtredalmhol}
 {\it For any $X\in \mathcal C$, smooth, the 
map $gt_X:X\ra GT(X)$ is almost holomorphic.}
\end{proposition}

\begin{question}\label{gtredgt?} Let $X\in \mathcal C$. Is $gt_X$ either constant, or a fibration of general type? 
\end{question}

This question has a positive answer if so does the next one:

\begin{question}\label{gtfibprod} Let $u:X\ra U$ and $v:X\ra V$ be fibrations
of general type. Is then $f:X\ra W$ of general type if so are $u$
and $v$?
(Here $f$ is the Stein factorisation of the product map
$(u\times v):X\ra W'$, with $W':=(u\times v)(X)\subset U\times V$).

\end{question}

It is obvious that if \ref{gtfibprod} has a positive answer, so does \ref{gtredgt?}, by the 
construction of $gt_X$. We shall see in \ref{core'} that these two questions have an affirmative answer. 
But this answer rests on the very different techniques of section \ref{orbadd}.

\begin{theorem}\label{gtrel} Let $X\in \mathcal C$, and
$f:X\mero Y$  be any fibration.  Then
$f$ admits a {\bf relative gt-reduction}.

This means that there exists a unique factorisation $f=h\circ g$ of $f$ by
fibrations $h: Z\mero Y$ and $g: X\mero Z$ such that for $y$ general in
$Y$, the restriction
$g_{y}:X_{y}\mero Z_{y}$ of $g$ to $X_y$ is the gt-reduction of $X_y$.
\end{theorem}

{\bf Proof:} This construction is actually in essence already in [C8O] to
which we refer for more details. We can and shall assume that $X$ is smooth 
and $f$ holomorphic, due to the bimeromorphic invariance of the notions involved.

We shall actually show a more precise version:

\begin{lemma}\label{gtrel'} Let $f:X\ra Y$ be a fibration, with $X\in \mathcal C$.

  After a generically base change $v:Y'\ra Y$ and proper modifications that we notationally 
ignore, there exists finitely many factorisations $f=h_i\circ g_i$, 
$i=1,2,...,N$,
with $g_i:X\ra Z_i$, $h_i:Z_i\ra Y$, such that:

\item[1.] The restriction $g_{i,y}:X_y\ra (Z_i)_y$ of each $g_i$ to the general fibre 
$X_y$ of $f$ is of general type.

\item[2.] if $g:X\ra Z$ is the 
(fibration part of) the Stein
factorisation of the product map $g_1\times ...\times g_N:X\ra 
Z_1\times _{Y}Z_2\times ....\times Z_{N-1}\times_{Y} Z_N$, then:

the restriction $g_{\mid X_y}:X_y\ra Z_y$ of $g$ to the general fibre 
$X_y$ of $f$ coincides with the $gt$-reduction $gt_{X_y}:X_y\ra 
G(X_y)$
of that fibre.
\end{lemma}

{\bf Proof:} Let, for $y\in Y^*$, 
$gt_{y}:X_{y}\ra Z_{y}$ be the
$gt$-reduction of $X_{y}$, where $Y^*$ is the Zariski open subset of
$Z$ over which $f$ is smooth.
For such a $y$, let $Z'_{y}$ be the family of fibres of $gt_{y}$,
defined as the image of the meromorphic map from $Z_{y}$ to $\mathcal
C$$(X_{y})$ sending a generic
point of $Z_y$ to the point in $\mathcal C$$(X_{y})$ parametrising its
reduced $gt_y$-fibre in $X_y$.

Because $gt_y$ is an almost holomorphic map by \ref{gtalmhol}, $Z'_y$ is an 
irreducible
component of $\mathcal C$$(X_{y})$.

Consider now the Zariski closed subset $\mathcal C$$(X/Y)$ of $\mathcal
C$$(X)$ consisting of cycles contained in some fibre of $f$. It is
naturally equipped with the
holomorphic map $f_y$$:\mathcal C$$(X/Y)\ra Y$ sending such a cycle to
the fibre containing it. (Strictly speaking, one may need to weakly 
normalise first, to make $f_*$ holomorphic, but this does not change the argument). 

Because $X\in \mathcal C$, the irreducible components of $\mathcal C$$(X/Y)$
are compact. 

Assume the fibre of $f$ is not special, for $y$ in a subset of $Y$ 
which is of second category, in Baire's terminology. 
 (Being of second category means: not contained in a countable union of 
closed subsets with empty interior. 
As we shall see later, the right topology here in our context 
is the Zariski topology, not the metric topology).

Because $X\in \mathcal C$, the irreducible components of $\mathcal C$$(X/Y)$
are compact. By the countability at infinity of $\mathcal C$$(X/Y)$, 
there is an irreducible
component $\Gamma'$
of $\mathcal C$$(X/Y)$
mapped surjectively onto $Y$ by $f_*$, and such that the $f_*$-fibre
$\Gamma'_y$ of $\Gamma'$ over $y$ has a component equal 
to $Z'_y\subset \mathcal C$$(X_y) $, the family of fibres of a fibration 
of general type $g_{i,y}:X_y\mero Z'_y$. This map is almost 
holomorphic, by \ref{gtalmhol}.

The Stein factoristation of $f_*$ restricted to $\Gamma'$ 
gives a finite base change for $Y$. This base changes we shall 
notationally
ignore, here, because they are irrelevant to our problem. So we deal 
as if the generic fibres of $f_*$$_{\mid \Gamma'}$ were irreducible.

Thus, for some $y\in Y^*$, the fibre $\Gamma'_y$ of 
$f_*$$_{\mid \Gamma'}$ is the family of fibres of some almost 
holomorphic
meromorphic fibration $g_{i,y}=\gamma_y: X_y\ra Z_{y,\gamma}$ on $X_y$. By 
the (obvious) openness of almost holomorphicity, one deduces the 
existence
of such a $\gamma_y$ for the generic $y\in Y$. And so, using the 
graph of the family, we get a factorisation $f=\delta\circ \gamma$, 
with
fibrations $\gamma:X\ra Z_{\gamma}$, and $\delta:Z_{\gamma}\ra Y$.

By our assumption, $\gamma_y$ is a fibration of general type for $y$ in 
$S\subset Y$ of second Baire category in $Y$. From \ref{semicont}, we 
conclude that
$\gamma_y$ is still of general type for $y$ general in $Y$.

The construction of the $g_i$'s is now obvious, by observing that if 
the map $g$ resulting from a finite family of $g_i$'s, $i+1,2,...,N$, 
does not
induce $gt_{X_y}$ on the general $X_y$, there exists, by the same argument as above, 
a component 
$\Gamma'$, inducing a general type fibration on the general 
$X_y$, and
such that its (Stein factorised) fibre product over $Y$ with the preceding ones will increase
 the dimension of the resulting $Z$. Contradiction.

This shows the lemma, and so \ref{gtrel} $\square$

\begin{definition}\label{SABCat'}

A subset $A\subset V$ of a complex analytic space is said to be of 
{\bf second Zariski category in $V$} if it is not contained in a countable union 
of Zariski closed subsets with empty interior of $V$.
(Notice that the definition makes sense in the algebraico-geometric context as well).
\end{definition}

From the proof \ref{gtrel}, we immediately get:

\begin{corollary}\label{specZR} Let $f:X\mero Y$ be a fibration, with $X\in \mathcal C$.

Assume that $dim(GT(X_y))=d$, for $y\in A$, where $A$ is of second Zariski category in $Y$.
 Then, this equality holds for the general point $y$ of $Y$.

In particular, if $X_y$ is special  for $y$ in a subset of second Zariski category in $Y$, the 
general fibre of $f$ is special.

\end{corollary}

{\bf Proof:} Let $f=h\circ g$ be the $gt$-reduction of $f$. By 
assumption, $dim(GT(X_y))=d$ for $y\in A$.
But also $dim(g(X_y))=dim(GT(X_y))$ for $y$ general in $Y$, and $dim(g(X_y))=d$ for 
$y$ generic in $Y$. Thus $dim(GT(X_y))=d$ for $y\in Y$ general $\square$\\

\subsection{\bf The Algebraic Reduction}

As an application of the preceding arguments, we show Theorem \ref{fibalgred}:

\begin{theorem} \label{fibAlgred} Let $a_{X}:X\mero Alg(X)$ be the
algebraic reduction of $X\in \mathcal C$. The generic fibre of $a_X$ is
special.

\end{theorem}

{\bf Proof:} Assume not.  By lemma \ref{gtrel'} above, after a suitable 
finite base change over $Alg(X)$ (which we
notationally ignore because it preserves the algebraic reduction and 
dimension),
  there exists a non-trivial factorisation $a_X=h\circ g$ with $g$ a
  fibration inducing a fibration of general type over the general 
fibre of $a_X$. Write it: $a_X=h\circ g$,
with $g:X\ra Z$ and $h:Z\ra Alg(X)$. Then: $dim(Z)>dim(Alg(X))$.

By construction, the line bundle $(K_{Y}+\Delta(g))$ over $Z$ is
thus $h$-big. Thus $Z$ is Moishezon, as one sees considering the line
bundle
$L:=(h^{*}(k.H)+(K_{Z}+\Delta(g)))$ on $Z$, which is big for
$k$ a large and positive integer, and $H$ an ample line bundle on 
$Alg(X)$, which we oviously can
choose to be projective. (See for example,
the Proof of [U],Theorem (12.1)). But this contradicts the definition
of $a_X$, and proves the claim $\square$

 \subsection{\bf The Category of Orbifolds}\label{orbcat}

We very briefly discuss without proofs the extension of part of 
our considerations to orbifolds, restricting here to 
prepared orbifolds $(Y/\Delta)$ with $Y$ smooth and the support of 
$\Delta$ an s.n.c divisor of $Y$
(but ultimately, one needs to consider $klt$ orbifolds).

One of the main point is to define bimeromorphic equivalence. The 
right notion seems to be
derived from {\bf terminal modifications}. 

\begin{definition}\label{tmodif} The 
bimeromorphic holomorphic map
$v:Y'\ra Y$ is said to induce a bimeromorphic map: $v:(Y'/\Delta')\ra 
(Y/\Delta)$ if it is  {\bf terminal}
with respect to the orbifold structures, that is if: 
$K_{Y'}+\Delta'=v^* (K_Y + \Delta) + \sum _{j\in J} a_j. E_j$, where
(as usual) the $a_j$ are all positive, and $J$ is the collection of 
$v$-exceptional divisors on $Y'$.(One might also define similarly the 
notion of {\bf canonical modification}, of course).
\end{definition}

Notice that the orbifold Kodaira dimension is invariant under 
bimermorphic equivalence of orbifold, which is the one generated by
terminal modifications. 

One can define for any fibration $g:(Y/\Delta)\ra Z$ its orbifold 
base, as in \ref{compfib}. One can extend this notion to the case of meromorphic $g$, 
by first resolving the indeterminacies of $g$ by a terminal modification.

The Kodaira dimension of this fibration is then the minimum of the 
Kodaira dimensions of the orbifolds bases of fibrations
equivalent to $g$, these being defined on orbifolds $(Y'/\Delta')$ 
bimeromorphically equivalent to $(Y/\Delta)$ .

\begin{definition} The fibration $g:(Y/\Delta) \ra Z$ is of {\bf general type} if 
its Kodaira dimension is $dim(Z)>0$.

The orbifold $(Y/\Delta)$  is {\bf special} if it has no 
fibration of general type.
\end{definition}

Fundamental tools for the study of orbifolds are the locally free sheaves
$\Omega^p_{Y}(log(\Delta))$ of logarithmic forms along $\Delta$ (classically known 
when the multiplicities are
infinite, or said differently, when $\Delta$ is reduced). We shall 
not give the definition here, but simply say that sections of this sheaf can be symbolically 
written locally in the standard normal crossing coordinates for 
fixed $q$ as linear combinations of expressions of the form: 

$h.(dy_{j_1}/(y_{j_1}^{(1-1/m_{j_1})})\wedge...\wedge (dy_{j_s}/y_{j_s}^{(1-1/m_{j_s})})\wedge dy_{j_{s+1}}\wedge...\wedge dy_{j_q}$, with $h$ 
holomorphic, and:

 $1\leq j_1<...<j_s\leq r<j_{s+1}<...<j_q$ if $y_1....y_r$ is a local equation of $\mid \Delta\mid$, 
the multiplicities being 
given by the $m_i$'s.

 A section of $\Omega_Y^q(log(\Delta))$ is thus an $m$-th root of a well-defined holomorphic tensor, when $m$ is an 
integer divisible by each of the $m_i$'s. 

More precisely, a section $s$ of this sheaf is defined as a pair $(F,s)$, where $F$ is a rank one coherent subsheaf of 
$\Omega_Y^q(log\mid \Delta\mid)$, and $s$ is a holomorphic section of $F^{\otimes m}$, for some $m$ 
divisible by all the $m_i$'s, and such that $s\in (\Omega_Y^q(log\mid \Delta\mid))^{\otimes m}(-m.\Delta^*)$, where 
$\Delta^*$ is the $\mathbb Q$-divisor on $Y$ defined by: $\Delta^*:=\mid \Delta \mid-\Delta=\sum_{i\in I} (1/m_j).\Delta_j$.

By lifting to a $\Delta$-nice covering (see section \ref{orbadd}), these 
sections become standard $p$-forms. From which one deduces the important property of $d$-closedness of such $log((\Delta))$-forms.

As we did above, one can then also define directly the Kodaira dimension of a fibration by introducing 
the saturation of the differential sheaf defined by $g$ in $\Omega_Y(log(\Delta))$. Because we may only consider high 
and divisible multiples to define the Kodaira dimensions, one does not need to define precisely $\Omega _Y(log(\Delta))$ to 
define this orbifold Kodaira dimension, and directly look at rank one subsheaves $F$ of $\Omega_Y(log\mid \Delta\mid)$, and define as usual the 
Kodaira dimension of $(F^{\otimes m}(-m\Delta^*))$ .

Their relevanco to our topic is that special orbifolds $(Y/\Delta)$ 
are characterised by the absence
of Bogomolov sheaves on $(Y/\Delta)$, defined as in \ref{bogsheav} above 
when $\Delta$ is empty,
just replacing there $\Omega^p_Y$ by $\Omega^p_Y(log\Delta)$. 

The correspondance between Bogomolov sheaves on $(Y/\Delta)$ and general type fibrations  
on this orbifold extends to this orbifold context. The orbifold additivity theorem then applies in this context.

The construction of the core for an orbifold can then be made by the second approach we followed (as the highest 
special fibration). The geometric approach seems more delicate than for varieties, and certainly needs some extra arguments, 
because one needs to take into account the order of contact of the subvarieties with the orbifold divisor.

\section{\bf The Core}\label{core}

\subsection{\bf The Core: Functoriality Properties}

We use the notations of section \ref{geomquot}.

\begin{definition}\label{defcore} Let $X\in \mathcal C$ be normal. 
Let $A:=A(X)\subset \mathcal C$$(X)$ be the family of special subvarieties of $X$. 
It is Z-regular. Let then $T(A)$ be the family of its components 
(see proposition \ref{strat}), and let $c_X:X\mero C(X)$ be the $T(A)$ quotient of $X$. 
This almost holomorphic fibration will be called {\bf the core of $X$}. (see section \ref{geomquot})
\end{definition}

In general, not much can be said about the fibres of $c_X$ (except that two of 
their points are joined by a chain of special subvarieties of $X$, 
at least if question \ref{specspec} has a positive answer).

\begin{example} Let $X$ be the cone over a variety of general type $V$. 
Then $c_X$ is the constant map. But $X$ is by no means special, since 
it has a $\mathbb P^1$-fibration over $V$. Note that this fibration 
is not almost holomorphic. 
\end{example} 

This example shows the role of singularities. Indeed:

\begin{theorem} \label{c_X} Let $X\in \mathcal C$ be {\bf smooth}.
  Let $c_X:X\ra C(X)$ be the core of $X$. Then:
\begin{enumerate}
\item[1.] The general fibre of $c_X$ is special.

\item[2.] If $F$ is a general fibre of $c_X$, and if $Z\subset X$ is a special subvariety of $X$
meeting $F$, then $Z$ is contained in $F$ (such a fibre $F$ will be said $c_X$-general).

\item[3.] $c_X$ is almost holomorphic.

\end{enumerate}
\end{theorem}

\begin{remark} The above result should hold true when $X$ is singular, 
provided it has at most canonical singularities.
Notice that $c_X$ is not a bimeromorphic invariant, in general (as the preceding example shows). 
But it is easily seen from \ref{c_X} to be so if $X$ is smooth: 
then $c_Y=c_X\circ m$, if $m:Y\ra X$ is bimeromorphic.
\end{remark}

  When $X$ is smooth, we denote by $ess(X)$ the dimension of $C(X)$, and call it the
{\bf essential dimension} of $X$. Thus $ess(X)=0$ iff $X$ is special,
and $ess(X)=d=dim(X)$ iff $X$ is of general type, by Theorem \ref{ess=n}
below. In the first case,
the core is the constant map; in the second case, it is the identity
map.

\begin{theorem}\label{corealg} Let $X\in \mathcal C$ be normal. Let $c_X:X\mero C(X)$ be its core, 
and let $a_X:X\mero A(X)$ be its algebraic reduction. There exists a factorisation $b_X:A(X)\mero C(X)$ 
of $c_X$ such that $c_X=b_X\circ a_X$. 

In particular: $C(X)$ is always Moishezon. 
\end{theorem}

{\bf Proof:} By Theorem \ref{fibalgred}, the fibres of $a_X$ (which is a bimeromorphic invariant of $X$) 
are special, hence contained in the fibres of $c_X$ $\square$

We shall prove Theorem \ref{c_X} below. Before we start with some easy observations.

\begin{proposition} \label{functc_X} Let $c_X:X\ra C(X)$ be the core of 
$X$, smooth. Let $h:Z\mero X$ be any
meromorphic map. Assume $h(Z)$ meets some $c_X$-general fibre of $c_X$. There 
exists then a natural
meromorphic map $c_h:C(Z)\mero C(X)$ such that $c_h\circ c_Z=c_X\circ h$.
\end{proposition}

{\proof:} By the assumption, if $z\in Z$ is general, its image 
$x:=h(z)$ in $X$ belongs to a $c_X$-general fibre $F_z$ of $c_X$.
The fibre of $c_Z$ through $z$ is special. Thus
so is its image $V_z$ by $h$. Since $V_z$ meets $F_z$, it is 
contained in $F_z$, by property [2.] in Theorem \ref{c_X}. Hence the existence of $c_h$ $\square$

\begin{corollary} \label{functc_X'} Let $h:Z\mero X$ be as in proposition \ref{functc_X} above. 
Then $c_h$ as above exists in the following cases:
\begin{enumerate}
\item[1.] $c_X\circ h:Z\ra C(X)$) is surjective. (In particular if $h$ itself is surjective).

\item[2.] $Z\subset X$ is the general member of a family $(Z_t)_{t\in T}$ 
of subvarieties of $X$
such that the varieties $c_X(Z_t)$ cover $C(X)$.

\item[3.] $Z=X_y$ is a general fibre of $(\psi\circ c_X)$, where 
$\psi:C(X)\ra Y$ is any fibration. In this case, $c_Z$ is
simply the restriction of $c_X$ to $Z$. 
\end{enumerate}
\end{corollary}

Let us give some easy examples in which $c_X$ can be described:

\begin{proposition}\label{gtspec} Let $X\in \mathcal C$ be a manifold, and assume that 
$f:X\mero Y$ is a special fibration of general type. Then: $f=c_X$

In particular: there is at most one fibration both special and of general type on $X$.
\end{proposition}

{\bf Proof:} Because $f$ is special, there is a factorisation $g:Y\mero C(X)$ such that $g\circ f=c_X$. 
Indeed: the general fibre $F$ of $f$ is special and meets the general fibre $C$ of $c_X$. Thus $F\subset C$. 

But $f$ is of general type, and so by Theorem \ref{domin}, there exists a factorisation $h:C(X)\mero Y$ such 
that $f=h\circ c_X$. Thus $f=c_X$, as claimed.$\square$

\begin{remark}\label{remsplitconj} We shall later conjecture (and prove in some cases (see section \ref{decconj})) 
that $c_X$ is a fibration of general type. So that 
$c_X$ should always be characterised as the unique fibration of domain $X$ both special and of general type.
\end{remark}

\begin{corollary}\label{k=n} Let $X$ be a manifold of general type. 
Then: $c_X=id_X$ (the identity map of $X$), and so: $ess(X)=dim(X)$.
\end{corollary}

{\bf Proof:} Indeed: $id_X$ is then special and of general type. Apply then \ref{gtspec} $\square$

\begin{remark}\label{remk=n} We shall see later in Theorem \ref{ess=n} that the converse also holds true: 
if $ess(X)=dim(X)>0$, then $X$ is of general type.
\end{remark}

\begin{corollary}\label{k=n-1} Let $X\in \mathcal C$ be a manifold, with $n:=dim(X)$. 
Then $ess(X)=(n-1)$ in the following two cases:
\begin{enumerate}
\item[a.] $\kappa(X)=(n-1)$, and $J_X$ is a fibration of general type. Then: $c_X=J_X$.
\item[b.] $R(X)$ (the rational quotient of $X$) is of dimension $(n-1)$ and of general type.
\end{enumerate}
\end{corollary}

{\bf Proof:} Indeed, in case (a.) (resp. (b.)), the generic fibre of $J_X$ (resp. $r_X$) 
is an elliptic (resp. rational) curve. 
The fibration $J_X$ (resp. $r_X$) is thus special. The other conditions imply that 
it is also of general type. It is thus the core of $X$. 
In particular, $ess(X)=(n-1)$ $\square$

\begin{remark}\label{remk=n-1} We shall see later in Theorem \ref{ess=n-1}, 
as a consequence of orbifold additivity theorems, 
that the converse also holds true: 
if $ess(X)=dim(X)-1>0$, then $X$ is of the type (a) or (b). For the case 
$ess(X)=(n-2)$ see below \ref{ess=n-2}.

Finally, in the same section, we shall also give a conjectural description of 
the core as a canonical and functorial 
composition of orbifold rational quotients and Iitaka fibrations.
\end{remark}

\begin{corollary}\label{fibcurv} Let $f:X\mero C$ be a special fibration, 
with $X\in \mathcal C$ and $C$ a curve.
Then: either $f$ is of general type, and $f=c_X$, or $f$ is not of general type, and $X$ is special.
\end{corollary}

{\bf Proof:} In the first case, the claim follows from Proposition \ref{gtspec}.
 In the second, from the fact that if 
$g:X\mero Z$ is a fibration of general type, there exists by \ref{domin} a 
factorisation $h:C\mero Z$ such that 
$g=h\circ f$. But $C$ is curve, and $f$ is not of general type. Thus $Z$ is 
a point, and $g$ is not of general type. 
Contradiction.$\square$

The core can be constructed in a relative setting, as well, by a simple application of 
Theorem \ref{relquot}:

\begin{theorem} \label{relc_X} Let $f:X\ra Y$, with $X\in \mathcal C$. There exists
a unique factorisation $f=g_f\circ c_f$ by two fibrations $c_f:X\ra C(f)$ and
$g_f:C(f)\ra Y$ such that, for $y\in Y$ general, the restriction 
$c_f:X_y\ra C(f)_y$
is the core of $X_y$.

We call the factorisation $f=g_f\circ c_f$ the {\bf core of f}.

\end{theorem}

We now prove Theorem \ref{c_X}:

\subsection{\bf Construction of the Core as the Lowest Special Fibration}

{\bf Proof of Theorem \ref{c_X}:} For this, we shall simply apply 
theorem \ref{quotstab} to the family $A(X)$ of 
special subvarieties of $X$, supposed to be smooth an to be in $\mathcal C$. 

We know that $A(X)$ is Z-regular (see section \ref{Zreg} for this notion). It is thus 
sufficient to show that $A(X)$ is also 
stable (see \ref{stab} for this notion). The property [stab2] is obtained by applying corollary \ref{specsurj}. 

The property [stab1] is the content of the next Theorem, which thus establishes 
at the same time Theorem \ref{c_X}:

\begin{theorem} \label{specstab} Let $T\subset T(A(X)) \subset \mathcal C$$(X)$ be a special
family as above, with $X$ smooth in $\mathcal C$. Assume each irreducible component
of $T$ is $X$-covering. Then:
\begin{enumerate}
\item[1.] $X$ is special if $X$ is $T$-connected (ie: $q_T$ is the constant map).
\item[2.] The general fibre of $q_T$ is special.
\end{enumerate}
\end{theorem}

{\bf  Proof:} The second assertion is a consequence of the first, because 
the if $X_y$ is a general fibre of $q_T$, then $X_y$ is smooth since 
$q_T$ is almost holomorphic, and the family $T_y$ consisting of $t\in T$ 
such that $V_t\subset X_y$ is a finite union of covering families  
of $X_y$ with general member special, and such that $X_y$ is $T_y$-connected. 

We thus only need to show the first statement: assume there existed a meromorphic fibration
$f:X\ra Y$ of general type. By \ref{gtalmhol}, $f$ is almost holomorphic,
since $X$ is supposed to be smooth.

   By Theorem \ref{domin} applied to each irreducible component $V_i$ of $V$, we
have a meromorphic factorisation $\phi:T\ra Y$ such that
$\phi\circ g=f\circ \psi:V\ra Y$.

Assume first that $f$ is holomorphic. Thus $f$ is constant on every
$V_t$, and so on every $T$-chain. Because $X$ is $T$-connected, $f$
takes the same value on two arbitrary points of $X$. Thus $f$ is
constant and $Y$ is a single point, in contradiction with the fact
that it is of general type.
Thus $X$ is special, as claimed.

If $f$ is only almost holomorphic, the same argument applies,
provided we choose an $f$-regular point $y\in Y$: for every $t\in T$,
if
$V_t$ meets $X_y$, then $V_t$ is contained in $X_y$ (because of the
factorisation property $\phi\circ g=f\circ \psi$, and the usual rigidity lemma.
More precisely: approximate $V_t$ in $\mathcal C$$(X)$ by a sequence 
$V_{t_n}$, such that
$V_{t_n}\subset X_{y_n}\subset f^{-1}(U)$, for some Stein-or affine in the algebraic category-
neighborhood $U$ of $y$ in $Y$, and $n$ large.
This is possible because the generic $V_{t'}$ is contained in a fibre of $f$.
Then $V_t\subset X_y$, by an easy argument, based on the fact that 
$f$ is holomorphic on $f^{-1}(U)$).

So
we get that the
generic member, hence every member of the family $T$ is contained in
some fibre
of $f$ (here the notion of fibre of $f$ is the usual Chow-scheme theoretic 
one, defined in section \ref{fib}).

Thus every $T$-chain meeting $X_y$ is contained in $X_y$.
Because $X$ is $T$-connected, $X$ is contained in $X_y$ and $f$ is
constant, and so not of general type, as assumed. Contradiction.
$\square$

We can now establish the following often useful characterisations of
the core:

\begin{theorem} \label{charcore} Let $f:X\mero Y$ be a special fibration,
with $X\in \mathcal C$, such that for any special fibration $g:X\mero Z$,
there exists a factorisation $\phi:Z\mero Y$ of $f:=\phi\circ g$. Then
$f=c_{X}$ is the core of $X$.

In particular, if $f:X\mero Y$ is a special fibration of general type,
it is the core of $X$.
\end{theorem}

The relevant diagram is:

\centerline{
\xymatrix{ 
X \ar[r]_-{f}\ar[d]_-{g} & Y\\  
Z  \ar@/_/[ur]_-{\phi}\\ 
}
}

{\bf Proof:} First, because $c_X$ is special, there exists a factorisation $\psi:C(X)\mero Y$. 
Let $F$ be a $c_X$-general fibre of $c_X$. By the existence of $\psi$, it is contained in some fibre $G$ of $f$. 
But $G$ is special, because $f$ is. By the defining property of $c_X$, we have the reverse inclusion $G\subset F$. 

The last assertion follows now from theorem \ref{domin}, because $f$ being of general type, 
for any special fibration $g$, theorem \ref{domin} shows that the factorisation $\phi$ exists.
$\square$

\subsection{\bf Rationally connected manifolds}\label{rc}

We now come to our first basic example of special manifolds: the rationally connected ones. 

We recall first their definition and some properties.

Recall
([C3],[Ko-Mi-Mo]) that an irreducible compact complex space $X$ 
is said to be {\bf rationally connected}
if any two generic points
of $X$ are contained in a rational chain of $X$ (ie: a connected
projective curve of $X$ , all irreducible components of which are
rational (possibly singular) curves).

Examples of rationally connected manifolds include unirational and
Fano Manifolds, twistor spaces. This property is bimeromorphically
stable
among manifolds, but not among varieties (the cone over a projective
manifold which is not rationally connected will again provide such an
example).

(Of course,the above definition can be given for algebraic varieties
defined over arbitrary fields).
We refer to [Ko-Mi-Mo] for some of the fundamental properties of this
class of manifolds.

Rational connectedness has a slightly different characterisation, by the 
following fundamental result ([G-H-S]):

\begin{theorem}\label{GHS} ([G-H-S]): Any fibration $f:X\ra C$ over 
a projective curve $C$ with $X$ smooth projective and generic fibre rationally connected 
has a holomorphic section.
\end{theorem}

\begin{definition}\label{ratgen} ([C95]) Let $X\in \mathcal C$ be 
irreducible. We say that
$X$ is {\bf rationally generated} if for any surjective meromorphic 
map $f:X\mero Y$, $Y$ is uniruled.
\end{definition}

Any rationally connected $X\in \mathcal C$ is thus rationally generated. 
But, conversely:

\begin{theorem} \label{RCRG} Let $X\in \mathcal C$ be rationally 
generated. Then $X$ is rationally connected.
\end{theorem}

{\bf Proof:} By induction on the dimension: $X$ is uniruled (take $f=id_X$).
Let $r_X:X\ra R(X)$ be the rational quotient of $X$. Then $R(X)$ is also
rationally generated (obvious).
By induction, it is rationally connected. So it is in particular 
Moishezon (by [C81], which says in particular that if any two points of $X\in \cal C$ can 
be joined by a chain of curves, then $X$ is Moishezon). Thus $X$ too is
Moishezon (by [C85], which says among others that $X\in \cal C$ is Moishezon 
if there is a fibration $u:X\mero Y$ with $Y$ and the generic fibre $F$ of $u$ Moishezon, if $F$ has 
$q(F)=0$) , because the fibres of $r_X$ are Moishezon, with vanishing irregularity. 
We can assume $X$ to be projective, by the 
bimeromorphic invariance
of the rational generatedness. But then the conclusion follows easily from [G-H-S], which 
allows to lift rational curves from $R(X)$ to $X$ $\square$

\begin{theorem} \label{RCspec} Let $X\in \mathcal C$ be a rationally connected {\bf manifold}. 
Then $X$ is special.
\end{theorem}

Notice that the smoothness of $X$ is essential, as shown again by the cone
over a projective manifold of general type. 

{\bf Proof:} Let $c_X:X\mero C(X)$ be the core. Assume it is not the constant map. 
Let $F$ be a $c_X$-general fibre of $c_X$. Because $X$ is rationally connected, some 
rational curve in $X$ meets $F$, but is not contained in $F$. Contradiction.$\square$\\

We now extend these results to the relative case.

\subsection{\bf  The Rational Quotient and The Core}

We now turn to the study of the {\bf rational quotient} of $X$ from
the point of view of special varieties. The rational quotient of
$X\in \mathcal C$ was
introduced in [C3] as an application of $T$-quotients. It was also
independently constructed
in [Ko-Mi-Mo] (under the name of
``maximal rationally connected fibration" (M.R.C for short))
by a different method based on their ``glueing lemma" for rational 
curves, in the algebraic context.)

\begin{theorem}\label{ratquot'} Let $X\in \mathcal 
C$ be normal. There exists a unique meromorphic fibration $r_X:X\mero R(X)$,
called the {\bf rational quotient} of $X$ such that:
\begin{enumerate}
 \item[1.] The general fibre of $r_X$ is rationally connected.

 \item[2.] The general fibre of $r_X$ contains any rational curve of $X$ that it meets.
\end{enumerate}
As usual, $r_X$ is almost holomorphic.
\end{theorem}

The proof is given in section \ref{geomquot}; see theorem \ref{ratquot}. 
Notice that, by Theorem \ref{relratquot}, the rational quotient also exists in relative version.

\begin{corollary}\label{RCquotspec} Let $X\in \mathcal C$ be smooth. The
rational quotient $r_X$ of $X$ is then a special fibration. There exists 
a factorisation:  $(cr)_X:R(X)\mero C(X)$ of $c_X=(cr)_X\circ r_X$.
\end{corollary}

{\bf Proof:} This is simply because $X$, and so the generic fibre of $r_X$ is special, by 
Theorem \ref{RCspec}. This shows the first assertion. Notice the second is obvious, and does not require $X$ to be smooth.
$\square$

We have the following easy property (already shown for $c_X$). The proof for $r_X$ is similar:

\begin{proposition} Let $f:X\ra Y$ be a surjective
meromorphic map, with $X\in \mathcal C$.

Then $f$ induces fonctorial maps (always denoted $f_{*}$):
$f_{*}:R(X)\ra R(Y)$ and $f_{*}:C(X)\ra C(Y)$.
\end{proposition}

In the above proposition, taking $f:=r_X$, we get in particular a natural map 
$(r_X)_*:C(X)\mero C(R(X))$.
For the rational quotient, we have a particular property not valid
for arbitrary special fibrations:

\begin{theorem} \label{CR=C} Let $X$ be smooth and Moishezon. Let
$r_X:X\ra R(X)$ be the rational quotient of $X$. Let $c_{R(X)}: R(X)\ra C(R(X))$ be the core of $R(X)$. 

Then: 
$(r_X)_*:C(X)\mero C(R(X))$ is bimeromorphic.
Equivalently: $(c_{R(X)}\circ r_{X}):X\ra C(R(X))$ is the core of $X$.
\end{theorem}

In other words: $C(R(X))=C(X)$.

\begin{remark} The hypothesis that $X$ is Moishezon can certainly be 
weakened to: $X\in \cal C$. For this, it is sufficient to make the same weakening 
in the hypothesis for $G$ in \ref{RCmultfree} below.
\end{remark}

{\bf Proof:} We have natural fibrations $\phi:R(X)\mero C(X)$ and
$\psi=(r_{X})_{*}:C(X)\mero C(R(X))$ defined above.
We thus only need to show that $\psi$ is isomorphic, or that the
general fibre $F$ of $(c_{R(X)}\circ r_{X}):X\mero C(R(X))$ is
special. 

Observe that we have, by restricting $r_X$ to $F$ a map
  $r_{X}:F\mero G$, where $G=(r_{X})(F)$
is the corresponding fibre of $c_{R(X)}$. Thus $F$ is fibered over
$G$, which is special, with fibres which are generically rationally
connected. The claim thus follows from the next:

\begin{proposition} Let $r:F\ra G$ be a fibration with
$F\in \mathcal C$ such that $G$ is Moishezon and special, and the generic fibre of $f$ is
rationally connected. Then $F$ is special.
\end{proposition}

{\bf Proof:} Let, if any, $g:F\ra H$ be an admissible holomorphic fibration of general type.
Since $f$ has special fibres, there exists by Theorem \ref{domin}  a
factorisation
$\phi:G\ra H$ of $g=\phi\circ f$. But now by Lemmas \ref{RCmultfree} and \ref{multfg} below, we
see that $\Delta(g)=\Delta(\phi)$. Thus $\phi$ is of general
type, too. But
this contradicts $G$ being special. Such a $g$ thus does not exist, and
$X$  is special.

\begin{lemma}\label{RCmultfree} Let $f:F\ra G$ be  a fibration with
generic fibres rationally connected and $G$ Moishezon.  Then $f$ is multiplicity free
(ie: $\Delta(f)$ is empty).
\end{lemma}

{\bf Proof:} We may assume that $G$ is projective. 
The claim then follows immediately from [G-H-S] (and is actually 
the most difficult part of the proof), by considering the restriction of $f$ over a very ample curve 
of $G$ meeting transversally any irreducible component of $\Delta(f)$ $\square$

\medskip

The proof of the following is immediate from 
the definition of multiplicities and the computation 
of the base orbifold divisor of a composed fibration, as in 
section \ref{compfib}:

\begin{lemma}\label{multfg} Let $f:F\ra G$ and $\phi:G\ra H$ be
fibrations. Assume that $f$ is multiplicity free. Then
$\Delta(\phi\circ f)=\Delta(\phi)$. 
\end{lemma}

\subsection{\bf Surfaces}

We can describe the core of a surface as follows, in terms of its rational quotient or Iitaka-Moishezon 
fibration. This will be extended to threefolds in the next section. A conjectural description in arbitrary 
dimension will be given in section \ref{c=rj^n}.

Recall that for any compact connected complex manifold $X$ with $\kappa(X)\geq 0$, we denote by 
$J_X:X\mero J(X)$ its Iitaka-Moshezon fibration. Let also $\kappa'(X):=\kappa(J(X),J_X)$. Obviously: 

$\kappa(X)\geq \kappa'(X)\geq -\infty$.

\begin{theorem}\label{coresurf} Let $X$ be a compact K\" ahler smooth surface. 
Then its core $c_X$ is described as follows:
\begin{enumerate}
\item[1.] $\kappa(X)=2$: then $c_X=id_X$, and $ess(X)=2$
\item[2.] $\kappa(X)=\kappa'(X)=1$: then $c_X=J_X$ and $ess(X)=1$.
\item[3.] $\kappa(X)=1>\kappa'(X)$: $X$ is special. 
\item[4.] $\kappa(X)=0$: $X$ is special.
\item[5.] $\kappa(X)=-\infty$ and $q(X)\geq 2$: then $c_X=r_X$, $ess(X)=1$.
\item[6.] $\kappa(X)=-\infty$ and $q(X)\leq 1$: $X$ is special.
\end{enumerate}
\end{theorem}

{\bf Proof:} If $\kappa(X)=2$, the claim is given by \ref{k=n}. 
If $\kappa(X)=1$, the fibration $J_X:X\mero C=J(X)$ is special, 
and the claim follows from \ref{fibcurv}. If $\kappa(X)=0$, $X$ 
is special from the facts just recalled above. If $\kappa(X)=-\infty$, 
from the classification of surfaces, $X$ is bimeromorphic to 
$\mathbb P^1\times C$, $C$ a curve with $g(C)=q(X)$, and $r_X$ is the projection 
to $C$ if $q>0$, and the constant map if $q=0$. The claims are then obvious $\square$

\begin{corollary}\label{esssurf} Let $X$ be a K\"ahler surface. Either $X$ is special (and $c_X$ is the constant map); 
or $\kappa(X)\geq 1$ and $c_X=J_X$ (the Iitaka fibration); or 
$\kappa(X)=-\infty$, and $c_X=r_X$ (the rational quotient of $X$).

One can compute $ess(X)$ as follows:
\begin{enumerate}
\item[1.] $ess(X)=2$ iff $\kappa(X)=2$;
\item[2.] $ess(X)=1$ iff $\kappa(X)\in \lbrace 1,-\infty\rbrace$ and 
$\pi_1(X)$ is not virtually abelian;
\item[3.] $ess(X)=0$ iff $\kappa(X) \leq 1$ and $\pi_1(X)$
is virtually abelian. 

(Recall that a group $G$ is said to be almost (or virtually) abelian if it has a finite 
index sugroup which is abelian)
\end{enumerate}
\end{corollary}

{\bf Proof:} All claims are deduced immediately from \ref{coresurf}, except 
for the ones concerning the fundamental group, when $\kappa(X)\leq 1$. 

If $\kappa(X)=-\infty$, then $\pi_1(X)\cong \pi_1(C)$, with the notations of the proof of \ref{coresurf}. 
The assertion is obvious. If $\kappa(X)=0$, we know that 
$\pi_1(X)$ is almost abelian from classification theory. If $\kappa(X)=1$, the 
assertion follows from Lemma \ref{ellsurf}, applied 
to $J_X$ $\square$

{\begin{corollary} \label{specsurf} A compact K\"ahler surface $X$ is special if and 
only if it has a finite \'etale cover which is bimeromorphic to one of the 
following surfaces:
\begin{enumerate}
\item[1.] $\Bbb P_{2}(\Bbb C)$
\item[2.] $\Bbb P_{1}(\Bbb C)\times E$, with $E$ elliptic.
\item[3.] K3, or Abelian.
\item[4.] Elliptic over a curve $C$ with $m$ multiple fibres, $C$ either rational (then $m\leq 2$) or 
elliptic (then $m=0$).
\end{enumerate}
\end{corollary}

{\bf Proof:} The surfaces listed above are special, by \ref{coresurf} above. Thus so are their undercovers. 
Conversely, if $X$ is special, it has a finite \' etale cover 
in the preceding list. This is clear if $\kappa(X)\leq 0$, by classification, and from 
the next Lemma \ref{ellsurf} if $\kappa(X)=1$ $\square$

{\begin{lemma} \label{ellsurf} Let $f:X\ra C$ be a relatively minimal 
elliptic fibration  on the compact K\"ahler surface $X$. 

\begin{enumerate}
\item[1.] Let $f^*(c):=\sum_{j\in J}m_j.D_j$ be any (scheme-theoretic) 
fibre $X_c$ of $f$. Then: 
its multiplicity $m(c,f):=inf \lbrace m_j\rbrace$ is also equal to: 
$m^+(c,f):=gcd \lbrace m_j\rbrace$.
\item[2.] There exists a finite \' etale cover $u:X'\ra X$ such that 
if $(v\circ f')=(f\circ u)$ is the Stein factorisation of $(f\circ u)$, with:
$f':X'\ra C'$ connected and $v:C'\ra C$ finite, then: $f'$ has no 
multiple fibre if $g(C')\geq 1$, and at most $2$ multiple fibres of coprime 
multiplicities if $C'$ is rational. 
\item[3.] Moreover: $g(C')=\kappa(C',f')=\kappa(C,f)$ in the preceding situation.
\item[4.] $X$ is special if and only if $\pi_1(X)$ is almost abelian.
\end{enumerate}
\end{lemma}

{\bf Proof:} (1) follows from Kodaira's classification of singular 
fibres of elliptic fibrations ([B-P-V; Chap.V.7]). This equality actually 
also follows from an elementary argument in the more general 
case of fibrations with generic fibre a complex torus. 

Assertion (2) follows from [C98'] and [N]: [N] shows that if a 
curve $C$ with points $a_1,...,a_m$, affected with multiplicities $n_1,...,n_m$ 
is given, there exists a cover $C'$ of $C$ ramified above the $a_i$'s 
only, each point above any $a_i$ having ramification exactly $n_i$. 
The only exception is when $C=\mathbb P^1,m=1, 2$, and  when $n_1\neq n_2$ if $m=2$.  
In [C98'], it is shown (it is a simple computation) that 
the base change over $C'$ leads to the sought after \' etale cover $u:X'\ra X$. 

Then property (3) follows from Theorem \ref{behKod} for the second equality, and 
from the fact that $m=0$ if $g(C')\geq 1$. 

We show (4): if $X$ is special then $\kappa(C,f)\leq 0$, so that $C'$ is rational or elliptic.
We apply [C98'], which shows that the natural sequence of maps: 

$$\pi_1(F')\ra \pi_1(X')\ra\pi_1(C')\ra1$$ 

is exact, with $F'$ a generic fibre of $f'$ (so that $F$ is an elliptic curve, 
and $\pi_1(F')\cong \mathbb Z^{\oplus 2}$). 
Thus $\pi_1(X')$ is almost abelian if $X'$ is  special, which is true if so is $X$ 
(because $f'$ is special and $\kappa(C',f')=\kappa(C,f)$).

Conversely, assume that $\pi_1(X)$ is almost abelian. Then so is $\pi_1(C')$, and $C'$ 
is either rational or elliptic. Thus $\kappa(C,f)\leq 0$, and $X$ is special by \ref{coresurf} $\square$

\subsection{\bf Higher Kodaira Dimensions}\label{higherkoddim}

We shall define higher Kodaira dimensions of any connected manifold $X\in \mathcal C$ as follows 
(this works for compact connected manifolds as well, actually):

The first Kodaira dimension of $X$ is the usual one: $\kappa(X)$.

If $\kappa(X)=-\infty$, the second Kodaira dimension of $X$ is not defined.

Otherwise: $\kappa(X)\geq 0$, and $J_X:X\ra J(X)$, the Iitaka-Moishezon fibration of $X$ is defined. 

Let then: $\kappa'(X):=\kappa(J(X),J_X)$.

We of course have: $\kappa(X):=dim(J(X))\geq \kappa'(X)\geq \kappa(J(X))\geq -\infty$.

If $\kappa'(X)=-\infty$, the next Kodaira dimension ($\kappa"(X)$) is not derfined. 

Otherwise: $\kappa'(X)\geq 0$.

Let then $J^{(0,1)}_X:J(X)\ra J'(X)$ be the Iitaka fibration defined on $J(X)$ by
the $\Bbb Q$-divisor:

 $(K_{J(X)}+\Delta(J_X))$ (for any admissible model of
the fibration $J_X$). 

Define: $J'_X::=J^{(0,1)}_X\circ J_X:X\ra J'(X)$.

Define next $\kappa"(X):=\kappa(J'(X),J'_X)\geq \kappa(J'(X))\geq -\infty$.

Of course, we have: $\kappa(X)\geq \kappa'(X)\geq \kappa"(X)\geq -\infty$, for the invariants so defined. 

Continuing inductively, we can define a decreasing sequence of invariants:

$\kappa(X)\geq \kappa'(X)\geq ...\geq \kappa^{(r)}(X)\geq -\infty$, and iterated orbifold Iitaka fibrations 

$J^{(r)}_X:X\mero J^{(r)}(X)$. If the sequence is defined til $J^{(r)}_X$, define $\kappa^{(r+1)}(X):=\kappa(J^{(r)}(X),J^{(r)}_X)$ , 
and if this is nonnegative, define $J^{(r,r+1)}_X$ as being the Iitaka fibration defined by the $\mathbb Q$-divisor 
$(K_{J^{(r)}((X))}+\Delta(J^{(r)}_X))$ (for any admissible model of
the fibration $J^{(r)}_X$).  

We thus have fibrations $J^{(r-1,r)}_X:J^{(r-1)}(X)\mero J^{(r)}(X)$ such that $J^{(r)}_X= J^{(r-1,r)}_X\circ J^{(r-1)}_X$.

Observe also that $\kappa(F_r,J^{(r)}_X)=0$, if $F_r$ is a general fibre of $J^{(r-1,r)}_X$, by the standard property of 
Iitaka fibrations.

This sequence stops at the first term (if any) equal to $-\infty$, and is stationary if any two 
terms $\kappa^{(r)}(X)=\kappa^{(r+1)}(X)$ are equal (and nonnegative, necessarily). 
This happens clearly if and only if the corresponding map $J^{(r)}_X$ is of general type.

The following is easily shown by induction on $r$:

\begin{proposition} The sequence of higher Kodaira dimensions is invariant under bimeromorphic maps and finite \' etale covers.
\end{proposition}

We shall later extend these notions, and even conjecture that these higher Kodaira dimensions are invariant under deformation (for $X$ K\" ahler).

As an illustration for the introduction of these invariants, we show the:

\begin{proposition}\label{ess=n-2} Let $X\in \mathcal C$ be smooth of dimension $n\geq 2$. 
Then $ess(X)=n-2$ in each of the following cases (a-e). Moreover, the core $c_X$ and its
 generic fibre $F$ (a special surface) are described as follows:

\begin{enumerate}
\item[a.] $\kappa(X)=(n-1)$, $\kappa'(X)=\kappa"(X)=n-2$. Then $c_X=J'_X$, and $\kappa(F)=1$, $\kappa'(F)=0$.
\item[b.] $\kappa(X)=(n-1)$, $\kappa'(X)=-\infty$, there exists 
$r:J(X)\mero Z$ with $dim(Z)=(n-1)$ such that $f\circ J_X$ is of general type. 
Then $\kappa(F)=1$, and $\kappa'(F)=-\infty$.
\item[c.] $\kappa(X)=\kappa'(X)=n-2$. Then $c_X=J_X$, and $\kappa(F)=0$.
\item[d.] $R(X)$ (the rational quotient of $X$) is of general type, of 
dimension $(n-2)$. Then $F$ is a rational surface.
\item[e.] $R(X)$ has dimension $(n-1)$, and $\kappa(R(X))=\kappa'(R(X))=(n-2)$. 
Then $F$ is birationally elliptic ruled.
\end{enumerate}
\end{proposition} 

{\bf Proof:} Case (a): Indeed, the fibre of $J'_X$ is special, because 
$F$ has an elliptic fibration $J:F\ra C$ with $\kappa(C,J)=0$, so the assertion 
follows from \ref{coresurf}. By assumption, $J'_X$ is of general type, 
because $\kappa'(X)=\kappa"(X)$. To show that $\kappa(F)=1$, 
use the easy addition theorem, applied to $X$ and $J'_X$: it says that 
$(n-1)=\kappa(X)\leq \kappa(F)+dim(J'(X))=\kappa(F)+(n-2)\leq dim(F)+(n-2)=(n-1)$.

We shall skip the proofs of the other cases, which are easier or similar $\square$

\begin{remark} We shall also nearly show the converse, as a consequence 
of additivity theorems, in section \ref{c=rj^n}. 
Actually, the converse holds under the general additivity conjecture.

The formulation of case (b) is unnatural. A natural formulation 
would require the conjectural 
notion of rational quotient for orbifolds. See \ref{c=rj^n}.
\end{remark}

\subsection{\bf Threefolds}

We shall describe the core of a compact K\" ahler threefold. For this we shall need 
Theorem \ref{k=0} shown later in section \ref{geomcons}, which says that $X$ is special if $\kappa(X)=0$, in all 
dimensions. 

\begin{theorem}\label{corenonspec3f}
 Let $X\in \mathcal C$ be a nonspecial threefold. The core $c_X$ of $X$ is a fibration of general type. 
Moreover, one can describe $c_X$, its generic fibre $F$, and $ess(X)$ as follows:
\begin{enumerate} 
\item[1.] $ess(X)=3$ iff $\kappa(X)=3$. Then $c_X=id_X$.
\item[2.] $ess(X)=2$ in the following two cases:
\begin{enumerate}
\item[a.] $\kappa(X)=\kappa'(X)=2$. Then $c_X=J_X$ is an elliptic fibration of general type.
\item[b.] $\kappa(X)=-\infty$, $R(X)$ is a surface of general type, 
and $c_X=r_X$ is a $\mathbb P^1$-fibration over $R(X)$.
\end{enumerate}
\item[3.] $ess(X)=1$ iff one of the following cases occurs:
\begin{enumerate}
\item[a.] $\kappa(X)=2,\kappa'(X)=\kappa"(X)=1$. Then $c_X=J'_X$ is a fibration of general 
type onto a curve, with $F$ a special surface with $\kappa(F)=1$ and $\kappa'(F)=0$.
\item[b.] $\kappa(X)=2$, $\kappa'(X)=-\infty$, $X$ nonspecial. Then $c_X$ is a fibration of 
general type onto a curve with $\kappa(F)=1$ and $\kappa'(X)=-\infty$.
\item[c.] $\kappa(X)=\kappa'(X)=1$. Then $c_X=J_X$ is a fibration of 
general type onto a curve with $\kappa(F)=0$.
\item[d.] $c_X=r_X$ is a fibration onto a curve of general type, with $F$ a rational surface.
\item[e.] $c_X=J_{R(X)}\circ r_X$ is a fibration of general type over a curve with $F$ a birationally 
ruled elliptic surface.
\end{enumerate}
\end{enumerate}
\end{theorem}

We now give a very rough list of the special threefolds:

\begin{theorem}\label{spec3f}
Any special threefold $X\in\mathcal C$ is one of the following:
\begin{enumerate}
\item[1.] $\kappa(X)=2$.
\begin{enumerate}
\item[a.]$\kappa'(X)=1>\kappa"(X)$: $J'_X$ is a non-general type fibration over a curve with fibre $F$ 
a special surface with $\kappa(F)=1>\kappa'(F)\in \{0,-\infty\}$.
\item[b.] $\kappa'(X)=0$: $J_X$ is an elliptic fibration with a klt orbifold base a normal 
surface with torsion canonical bundle $(K_S+\Delta)$. (The log-Enriques case).
\item[c.] $\kappa'(X)=-\infty$: $X$ has either a non-general type fibration over a curve 
with generic fibre a special surface $F$ with $\kappa(F)=1,\kappa'(F)=-\infty$, or 
an elliptic fibration with base orbifold a klt normal surface with Picard number one, 
and log-Del Pezzo (ie:-$(K_S+\Delta)$ is ample).
\end{enumerate}
\item[2.] $\kappa(X)=1$: then $J_X$ is a non-general type fibration 
over a curve with generic fibre $F$ a surface with $\kappa(F)=0$. 
\item[3.] $\kappa(X)=0$
\item[4.] $\kappa(X)=-\infty$: $X$ is either rationally connected, or a 
fibration over an elliptic curve with generic fibre a rational surface, or a $\mathbb P^1$-fibration 
over a special surface $S$ with $\kappa(S)\geq 0$.
(The case where $X$ is simple non-Kummer conjecturally does not exist, but strictly speaking additionally 
belongs to the last 
part (4) of the above list, because $\kappa(X)\leq 0$, then).
\end{enumerate}
\end{theorem}

\bf Proof:} We shall prove both results \ref{corenonspec3f} and \ref{spec3f} at the same time.

The case where $\kappa(X)=-\infty$ is clear from the section \ref{ratquot} above, 
because $X$ is then uniruled if $X$ is projective (by [Miy]), and by [C-P] otherwise if $X$ is not 
simple.

So we proceed case-by-case, assuming $\kappa(X)\geq 0$.

If $\kappa(X)\geq 0$, we are done by Theorem \ref{k=0}. If $X$ is non-projective, we could 
 also have applied [C-P], theorem (8.1), which says that if $X$ has a
nonzero holomorphic $2$-form, it is
covered by either a torus, or by the product of an elliptic curve and
a $K3$-surface. Thus $X$ is special in this case, too. (The existence
of a nonzero $2$-form when $X$
is non-projective is a famous result of K.Kodaira).

We now classify the cases occuring according to the pairs
$(\kappa(X),\kappa'(X))$, $\kappa(X)\geq 1$.

When the two terms are equal and positive,
   we conclude from \ref{gtspec} that $J_X$ is the core, and of general type.

If $\kappa(X)=1$, and $ \kappa'(X)\leq 0$, we conclude from \ref{fibcurv} that
$X$ is special. 

We are thus left with the cases where $\kappa(X)=2$, $\kappa'(X)\leq 1$.

We thus now assume that $\kappa(X)=2$.

Assume first that $\kappa'(X)=1$. Consider the map
$ J'_{X}:X\ra J'(X)$ (using the notations of
\ref{higherkoddim}): it has general special fibres, because if
$F$ is its general fibre, and $G$:=$J_{X}(F)$, 
then the restriction to $F$ of $J_{X}$ defines $J':F\ra G$, which has
generic fibres elliptic curves, while $G$ is a curve and, by the
definition of
$J^{0,1}_{X}$, we have: $\kappa(G,J')=0$. We conclude from \ref{fibcurv} 
that $F$ is a special
surface. The
easy addition theorem shows that $\kappa(F)=1$.

If $\kappa"(X)=1$, $J'_X$ is a fibration both special and of general type.
   So we conclude from \ref{gtspec} that it is the core of $X$. 

Otherwise (if $\kappa"'X)\leq 0$), $X$ is special. We are thus done with
this case ($ess(X)=0$, or $1$).

We are  now left with the more difficult case when $\kappa'(X)\leq 0$. 
Assume that
$X$ is not special. Let $f:X\ra Y$ be a fibration of general type.
The fibration $J_{X}$ being
special, we get from (1.12) a factorisation $\phi:J(X)\ra Y$ of
$f=\phi\circ J_{X}$. Thus $Y$ has to be a curve.

To conclude the proof we thus just need to show that this does not
happen if $\kappa'(X)=0$, because of the Minimal Model Program applied to a klt surface orbifold, 
as described below. The assertion we need follows from Proposition \ref{k=oorbsurf} below. 
But we need first some definitions for its statement.

\subsection{\bf Orbifold Surfaces}\label{orbsurf}

A {\bf surface orbifold} will be a klt pair $(S,\Delta)$, with $S$ a normal projective surface, together with 
an orbifold divisor $\Delta$ on it. 

If $g:S\ra C$ is a holomorphic fibration, we define the {\bf orbifold base} of $g:(S/\Delta)\ra C$ as the 
pair $(C,\Delta_g:=\Delta(g,\Delta))$, where $\Delta_g$ is the Weil $\mathbb Q$-divisor on $C$ defined as follows: 

$\Delta_f:=\sum_{c\in C}(1-1/m(c,g,\Delta)).c$, with $m(c,g,\Delta):=inf_{j\in J}\lbrace m_j.n_j\rbrace$, and 
$g^*(c):=\sum_{j\in J}m_j.D_j$. Here $n_j$ is simply $(1-d_j)^{-1}$, if $d_j=(1-1/n_j)$ is the multiplicity of 
$D_j$ in $\Delta$. 

\begin{remark}\label{D(gf)} As in \ref{compfib}, we see that if $h:X\ra S$ is a fibration with $\Delta=\Delta(h)$, 
then $\Delta_f=\Delta(g\circ h))$ 
(on a suitable model of $h$).
\end{remark}

We define then as usual the canonical bundle and Kodaira 
dimension of $(C/\Delta_f)$ and $(S/\Delta)$ (by definition, $(K_S+\Delta)$ 
is supposed to be $\mathbb Q$-Cartier).

We then say that {\bf $g$ is of general type} if 
$\kappa(C/\Delta_f)=1$, and that $(S/\Delta)$ is {\bf special} 
 if $\kappa(S/\Delta)<2$, and if there is no holomorphic fibration of general type $g:(S/\Delta)\ra C$ onto a curve.

\begin{proposition}\label{k=oorbsurf} Let $(S/\Delta)$ be a 
surface orbifold. If $\kappa(S/\Delta)=0$, then $(S/\Delta)$ is special.
\end{proposition}

{\bf Proof:} Assume there exists a general type fibration $g$ 
on $(S/\Delta)$. We shall show that $\kappa(S/\Delta)\neq 0$. 

We thus apply the MMP to our initial pair $(S/\Delta)$. This produces 
a sequence of elementary contractions of the form: 
$k:(S/\Delta)\ra (S'/\Delta')$, with $(S'/\Delta')$ still a klt pair, 
such that after at most $\rho(S)-1$ steps, one gets 
for the final pair (denoted also $(S'/\Delta'))$) one of the three basic cases:

(1) $K_{(S'/\Delta')}$ is nef.

(2) There is a fibration $g':S'\ra C'$ onto a curve such that 
$-K_{(S'/\Delta')}$ is $g'$-ample, and $\rho(S')=2$.

(3) -$K_{(S'/\Delta')}$ is ample, and $\rho(S')=1$. (``log-Del Pezzo" case).

Notice that at each step, the Kodaira dimension of the pair $(S/\Delta)$ is preserved, and 
that the curve being contracted is rational smooth (because $S$ itself is klt).

We refer to [K-M] and [F-M] for the existence and usual properties of these reduction steps. 

In cases (2,3), we have $\kappa(X'/\Delta')=\kappa(X/\Delta)=-\infty$. So these cases do not occur, 
because we assumed that $\kappa(X/\Delta)=0$.

Thus: $K_{(S'/\Delta')}$ is nef. 

{\bf Claim:} $D':=(K_{S'}+\Delta')\equiv 0$. Indeed, we have : $(D')^2=0$, otherwise we had: $\kappa(X/\Delta)=2$. 
From [F-MK],(11.3), we get the claim (their theorem asserts that $\kappa:=\kappa(X'/\Delta')$ is equal to the numerical 
Kodaira dimension of $\kappa(X'/\Delta')$). 

The rest of the proof rests on the following two lemmas \ref{Egt} and \ref{elcont}. 
The first one will be proved at the end of this section.

\begin{lemma}\label{Egt} Let $g:(S/\Delta)\ra C$ be a fibration of 
general type, with $(S/\Delta)$ klt. Let $E$ be a rational curve 
mapped surjectively onto $C$ by $g$. Then: $\Delta.E>2$.
\end{lemma}

\begin{corollary}\label{elcont} Let $k:(S/\Delta)\ra (S'/\Delta')$ is the 
contraction of an irreducible smooth rational curve $E$, 
with: $S,S'$ klt surfaces, $(K_S+\Delta).E\leq 0$, and $\Delta':=k_*(\Delta)$.

If $g:(S/\Delta)\ra C$ is a fibration of general type, then there exists $g':(S'/\Delta')\ra C$ such that $g=g'\circ k$. 
Moreover, $g'$ is still of general type, and $(K_S+\Delta)=k^*(K_{S'}+\Delta')$ if $(K_S+\Delta).E=0$.
\end{corollary}

{\bf Proof (of \ref{elcont}):} The second assertion is clear, if the first one is. 

This is because $m(c,g',\Delta)\leq m(c,g',\Delta')$, since in the definition of the left hand side 
of the inequality, the infimum is
 taken over a smaller subset. 

The first assertion is clear also if $g(C)\geq 1$, 
because the rational curve $E$ contracted by 
$k$ cannot be mapped surjectively to $C$ by $g$. 

We shall show that this also cannot happen when $C$ is rational, because $g$ is of general type 
(actually, as the proof shows, the condition $\kappa(C,g)\geq 0$ is sufficient, even).

We use the following numerical conditions: 

(a) $E^2\leq 0$ ($E$ exceptional)

(b) $(K+E).E=-2$ ($E$ rational, smooth)

(c) $(K+\Delta).E\leq 0$ 

So we assume by contradiction that $E$ is mapped onto $C$ by $g$. 

We then get: $\Delta.E\leq-K.E=2+E^2\leq 2$, which contradicts \ref{Egt} $\square$

We now complete the proof of \ref{k=oorbsurf}: 

We then apply the Minimal odel program to $S'$, but relative to $K_{S'}$. At each step again only smooth 
rational $K$-negative curves are contracted. Let $(S",\Delta")$ be the resulting pair. it has, by \ref{elcont}, 
all properties of $(S',\Delta')$. In particular, $D":=(K_{S"}+\Delta")\equiv 0$, and $g":(X"/\Delta")\ra C$ 
is of general type. Put $K:=K_{S"}$.

Assume first that $\kappa(X)\geq 0$. 

In addition to the above properties, $K_{S"}$ is nef. By our assumptions, 
the generic fibre of $g"$ is elliptic, and $\Delta"$ is 
``vertical" (ie: contained in fibres of $g"$).

Thus: $K^2=0$, $D".\Delta"=0=(D")^2$. Thus: $(D")^2=0$. And so $D"$ is a union 
of complete fibres of $g"$, by Zariski's Lemma.

There exists thus an orbifold divisor $\delta"$ on $C$ such that $\Delta"=(g")^*(\delta")$. 

Because $\kappa(X")\geq 0$, and $\kappa(C/\delta")=1$, we easily 
get that $\kappa(X"/\Delta")\geq 1$ (see Lemma \ref{addklem1} below, for example). A contradiction.

We now treat the remaining case in which $\kappa(X)=-\infty$. 

Now: $-D"$ is $g"$-ample, and every fibre of $g"$ has an irreducible reduction. (As above, $\Delta"=(g")^*(\delta")$, 
for an orbifold structure $\delta"$ of general type on $C$. 

We thus have: $\rho(S")=2$, and the arguments of [F-MK], theorem 11.2.3, show that: 

$K_{S"}+(\Delta")^{hor}\equiv \lambda F"$, $(\Delta")^{hor}$ is the horizontal part of $\Delta"$, defined 
as usual, and $F"$ is any fibre of $g"$. Here $\lambda \in \mathbb Q$ is such that $\lambda \geq (2q-2)$, with $q$ the 
genus of $C$. 

Thus $D"\equiv (\lambda+deg(\delta")).F"$ has Kodaira dimension at least $1$, because $deg(\delta")>(2-2q)$. 
Contradiction $\square$

We still have to show Lemma \ref{Egt}. 

{\bf Proof of \ref{Egt}:} 

From Hurwitz's formula, we get: $-2=-2d+\sum_{e\in E}(r_e-1)$, where 
$d$ is the degree of the restriction $h:E\ra C$ of $g$ to $E$, and for each 
$e\in E$, $r_e$ is the ramification order of 
$h$ at $e$.

Fix $c\in C$. Write: 

$g^*(c)=\sum_{j\in J} m_j.D_j$;

$m:=inf\{m_j.n_j,j\in J\}$, 

$n_j:=(1-d_j)^{-1}$, $d_j$:=multiplicity of 
$D_j$ in $\Delta$. 

Let $E_c:=E\cap S_c$, where $S_c:=g^{-1}(c)$ is the fibre of $g$ over $c$.

{\bf Claim:} $\sum_{e\in E_c}(r_e-1)\geq (1-1/m).d-\Delta_c.E$, where $\Delta_c$ is the union of components 
of $\Delta$ contained in $S_c$, with their corresponding multiplicities.

Then Lemma \ref{Egt} is an easy consequence of the claim. Indeed: 

 $-2=-2d+\sum_{e\in E}(r_e-1)$

$\geq -2d+\sum_{c\in C}\sum_{e\in E_c}(r_e-1)$

$\geq -2d+g^*(\delta).E-\Delta.E$, if $\delta:=\Delta(g,\Delta)$

$=(deg(\delta)-2).d-\Delta.E$

But $deg(\delta)-2>0$, because $g$ is of general type. Hence the conclusion of \ref{Egt}. 

To complete the proof, we establish the preceding {\bf claim}.

{\bf Proof of the claim:}

$(1-1/m).g^*(c).E=(1-1/m).\sum_{j\in J}.m_j.(D_j.E)= \sum_{j\in J}(m_j-m_j/m).(D_j.E)$

$\leq \sum_{j\in J}((m_j-1)+(1-1/n_j)).(D_j.E)$ (because $m_j.n_j\geq m, \forall j\in J$)

$=(\sum_{j\in J}(m_j-1).(D_j.E))+\Delta.E$.

We are thus reduced to show that: 

$\sum_{j\in J}(m_j-1).(D_j.E)\leq \sum _{e\in E_c}(r_e-1)$.

Because $\sum_{j\in J} m_j.(D_j.E)=\sum _{e\in E_c} r_e=d$, we just need to establish that: 

$\sum_{j\in j}D_j.E\geq \sum_{e\in E_c}1$, which itself follows from the inequality: 

$\sum_{j\in J} (D_j.E)_e\geq 1, \forall e\in E_c$, where the intersection number $(D_j.E)_e$ is the local 
intersection number near $e$. (Recall that $S$ is only assumed to be normal). 

We now show this last inequality.

Indeed: we have $\sum_{j\in J} (m_j.D_j)_e=r_e$, and the conclusion follows from the inequality:

$m_j\leq r_e, \forall j\in J$.

To show this inequality, we make a base change over $h:E\ra C$. Let $X'$ be the normalisation of 
the fibre product $X\times_{C}E$, and let $k:X'\ra X$, $g':X'\ra E$ such that $k\circ g=h\circ g$ 
be the natural maps induced from this base change. 

Let $E'\subset X'$ be the lift of $E$ to $X'$: it is a section of $g'$. Let $e'$ be the point of $E'$ 
lying above $e$.

The components of $(g')^*(c')$ which contain $e'$ are thus reduced. This easily implies the conclusion, by 
looking at a generic point of $D_j$, near which the projection $g$ is locally given by the equation $g(t,z)=z^{m_j}$. 
The fibre product is thus locally given by an equation of the form: $\zeta^{r_e}=z^{m_j}$. 

Dividing by $d:=gcd(r_e,m_j)$, we can assume that $d=1$, since we normalised the fibre product. 

Thus we have on $X'$ local coordinates $(t,s)$ with $\zeta=s^{m_j}$, and $z=s^{r_e}$. 

Locally, the projection $g'$ is given in these coordinates by $g'(t,s)=\zeta$. From which we deduce that 
$m_j=1$ Thus $m_j$ divides $r_e$. This in particular proves the claimed inequality $\square$ 

\begin{remark} What we actually proved in this section is the orbifold Additivity for Kodaira dimensions 
in dimension $2$. The general case is stated in the next section at \ref{cnmtheor}.
\end{remark}

\section{\bf Orbifold Additivity}\label{orbadd}

\subsection{\bf Orbifold Conjecture $ C_{n,m}$}

We use the notations and notions introduced in section \ref{compfib}.

So, if $g:(Y/H)\ra Z$ is a holomorphic fibration from the manifold $Y$ 
equipped with the orbifold divisor $H$, we defined the base orbifold 
$\Delta=\Delta(g,H)$ on $Z$.

The fundamental property of this definition is that, when 
$H=\Delta(f)$, for some fibration $f:X\ra Y$, then $\Delta=\Delta(g\circ f)$
 (for suitable models of $f,g$ which can be choosen so that 
$f,g,g\circ f$ are admissible).  

\begin{definition}\label{high}
We shall say that $g:(Y/H)\ra Z$ is {\bf high} (resp. {\bf very high})
if there exists 
a modification $u_0:Y\ra Y_0$ with $Y_0$ smooth 
such that (a) and (b) (resp. (a) and (b')) 
below are satisfied:

(a) Every $g$-exceptional divisor of $Y$ is $u_0$-exceptional.

(b) $\kappa(Y/H)=\kappa(Y_0/H_0)$, with $H_0:=(u_0)_*(H)$.

(b') $K_Y+H\geq (u_0)^*(K_{Y_0}+H_0)$.
\end{definition}

The following properties are immediate (using the arguments of \ref{existadmiss} and \ref{compofib}, for (3)):

\begin{proposition}\label{prophigh}
(1) If $g$ is very high, it is high.

(2) If $g$ is high, then $\kappa(Y,K_Y+H+B)=\kappa(Y,K_Y+H)$, 
for any effective $g$-exceptional $\mathbb Q$-divisor $B$ on $Y$.

(3) If $H=\Delta(f)$, then one can choose models of $f,g$ in such a way 
that $g$ is high, prepared, that $\Delta(g\circ f)= \Delta(g, H)$, 
and $f,g,g\circ f$ are admissible.
\end{proposition}

Now we can state the orbifold additivity conjecture $C^{orb}_{n,m}$:

\begin{conjecture} $(C^{orb}_{n,m})$ \label{$C^{orb}_{n,m}$} Let $g:(Y/H)\ra Z$ 
be a holomorphic fibration between manifolds, as above, with $Y\in \mathcal C$. 
Assume $g$ is prepared and high.

Then: $\kappa(Y/H)\geq \kappa((Y/H)_z)+\kappa(Z/\Delta(g,H))$, where 
$z\in Z$ is general and $(Y/H)_z:=(Y_z/H_z)$.
\end{conjecture}

Of fundamental importance for the considerations of the present 
paper is the following special case, shown by suitably adapting the classical 
methods of proof (T.Fujita, Y. Kawamata, E. Viehweg):

\begin{theorem} $(C^{orb}_{gt})$ \label{cnmtheor} Let $g:(Y/H)\ra Z$ 
be a holomorphic fibration between manifolds, as above, with $Y\in \mathcal C$, and $Z$ projective. 
Assume $g$ is prepared, high and of general type (ie: $\kappa(Z/\Delta(g,H))=dim(Z)$).

Then: $\kappa(Y/H)=\kappa((Y/H)_z)+dim(Z)$, 
where 
$z\in Z$ is general and $(Y/H)_z:=(Y_z/H_z)$.
\end{theorem}

Of course, the above $C^{orb}_{n,m}$ is a simple generalisation and refinement 
of the classical conjecture of S.Iitaka, dictated by the constructions made in the previous chapters. 

Let us list some of its corollaries or special cases of the above conjecture: 

\begin{proposition}\label{c'nm} Assume $C^{orb}_{n,m}$ holds.
Let $f:X\ra Y$ and $g:Y\ra Z$ be fibrations, with $X\in\mathcal C$. Then: 
$\kappa(Y,f)\geq \kappa(Y_z,f_z)+\kappa(Z,g\circ f)$, where 
$z\in Z$ is general and $f_z:X_z\ra Y_z$ is the restriction of $f$.
\end{proposition}

{\bf Proof:} By \ref{compfib} and \ref {semicont}, we can choose models 
of $g$ and $f$ in such a way that $f,g,g\circ f$ are admissible prepared, with $g$ high, 
$\Delta(g\circ  f)=\Delta(g, \Delta(f))$, and such moreover that $f_z$ is admissible. 

We conclude then from $C^{orb}_{n,m}$ and the following equalities: 

$\kappa(Y/\Delta(f))=\kappa(Y,f)$, $\kappa(Y_z/\Delta(f)_z)=\kappa(Y_z,f_z)$, and:
 
$\kappa(Z,g\circ f)=\kappa(Z/\Delta(g\circ f))=\kappa(Z/\Delta(g,\Delta(f)))$ $\square$ \\

In the special case where $X=Y$, we get: 

\begin{proposition}\label{c"nm} Assume $C^{orb}_{n,m}$ holds.
Let $g:Y\ra Z$ be a fibration, with $Y\in\mathcal C$. Then: 

$\kappa(Y)\geq \kappa(Y_z)+\kappa(Z,g)\geq \kappa(Y_z)+\kappa(Z)$.
\end{proposition}

The extreme inequality is of course the classical Iitaka conjecture.

Let us first list some immediate consequences of theorem \ref{cnmtheor}: 

\begin{corollary}\label{cnmcork>=o} Let $f:X\ra Y$ and $g:Y\ra Z$ be fibrations, 
with $X\in\mathcal C$. Assume $(g\circ f)$ is of general type, then: 
$\kappa(Y,f)=\kappa(Y_z,f_z)+dim(Z)$. 
\end{corollary}

\begin{corollary}\label{cnmspeccase} .
Let $g:Y\ra Z$ be a fibration, 
with $Y\in\mathcal C$. Assume $g$ is of general type, then: 
$\kappa(Y)=\kappa(Y_z)+dim(Z)$.
\end{corollary}

We shall now give the proof of Theorem \ref{cnmtheor}. It is classically done in two steps: first an 
easy reduction to weak-positivity statements for direct images by $g$ of twisted pluricanonical forms; 
second the proof of semipositivity. The first step is entirely similar 
to the known cases, so we shall be brief on it. The second step is simply obtained by introducing the 
orbifold divisors at appropriate places in the classical proofs of Y.Kawamata and E.Viehweg. See also 
the initial work [Fuj78].

\subsection{\bf Reduction to Weak Positivity} \label{reducwp}

We start by briefly recalling the notion of weak-positivity introduced in [Vi 83] (see also the survey [E]).

A torsionfree coherent sheaf $\mathcal F$ on $Z$, projective, 
is said to be {\bf weakly-positive} (written (w.p) for short) 
if for any ample line bundle $A$ on $Z$, and every $a>0$, integer, 
there exists an integer $b>0$ such that $S^{ab}(\mathcal F)\otimes A^b$ is generated 
over some (nonempty) open subset $U$ of $Z$ by its global sections (defined over $Z$).

Here $S^{ab}(\mathcal F)$ denotes the extension to $Z$ of the sheaf 
denoted by the same symbol, naturally defined over the open subset where $\mathcal F$ is locally free.

\begin{remark} \label{remwp} The following properties are shown in [Vi 82]:

(1) If $\mathcal F$ is locally free and nef, it is w.p.

(2) Let $v:Z'\ra Z$ be bimeromorphic, and 
$\mathcal F\subset \mathcal G$ a inclusion of torsionfree coherent 
sheaves of the same rank on $Z'$. If $\mathcal F$ is w.p, then so is $v_*(\mathcal G)$.

(3) If $v: Z'\ra Z$ is a ramified flat covering, with $Z,Z'$ smooth, 
and if $\mathcal F$ is torsionfree coherent on $Z$, 
then $\mathcal F$ is w.p if so is $v^*(\mathcal F)$.
\end{remark}

We now state without proofs two lemmas, shown but not separately stated in [E],
 and in various more or less implicit 
forms in [Vi 82] and [Kw 81]. Together with the weak-positivity result shown in the next section, 
they imply immediately Theorem \ref{cnmtheor}.

\begin{lemma}\label{addklem1} Let $g:Y\ra Z$ be a fibration with $Z$
projective. Let $E$ and $L$ be $\mathbb Q$-divisors on $Y$ and $Z$
respectively, such that:
$L$ is big  and $\kappa(Y,E)\geq 0$. Then:

$\kappa(Y,E+g^{*}(L))=dim(Z)+ \kappa(Y_{z}+E_{\mid
Y_{z}})$, for $z$ general in $Z$.
\end{lemma}

The crucial place where weak-positivity enters is to check that 
 $\kappa(Y,E)\geq 0$:

\begin{lemma}\label{addklem2} Let $g:Y\ra Z$ be a fibration, $E$ a 
line bundle on $Y$ and $L$ a 
$\mathbb Q$-divisor on $Z$.

Assume $L$ is big and $g_*(E)$ is weakly positive and nonzero.

Then: $\kappa(Y,E+g^*(L))\geq 0$.
\end{lemma}

In the next section, we shall show: 

\begin{theorem}\label{semipos} In the situation and hypothesis of Theorem \ref{cnmtheor}, 
for any sufficiently divisible integer $m>0$, there exists an effective $g$-exceptional 
divisor $B$ on $Y$ such that the sheaf:
 
$g_*(m(K_{(Y/(Z/\Delta(g,H)))}+H)+B)$ is weakly positive on $Z$. 

If $g$ is very high, we can choose $B=0$.
\end{theorem}

Let us now explain how to deduce Theorem \ref{cnmtheor} from the preceding 
lemmas and Theorem \ref{semipos}: 

apply first \ref{addklem2} to 
$E:=m.(K_{(Y/(Z/\Delta(g,H)))}+H)+B$, with $m>0$ 
an integer sufficiently divisible, so chosen that $g_*(E)$ is nonzero (and w.p), and:
 
$L:=(m/2). K_{(Z/\Delta(g,H))}$, which is big by hypothesis. 

We conclude that $\kappa(Y, E+g^*(L))\geq 0$.

Thus $E'=E+g^*(L)+H$ is such that $\kappa(Y,E')\geq 0$.

Next apply \ref{addklem1} to $E'$ and $L$ to conclude that $\kappa(Y, K_Y+H+B/m)$ satisfies the inequality 
stated in \ref{cnmtheor}. Use finally the fact that $g$ is high to conclude the proof of \ref{cnmtheor},
 because $\kappa(Y, K_Y+H+B/m)=\kappa(Y, K_Y+H)$ $\square$

\subsection{\bf Orbifold Weak-Positivity}\label{orbwp}

Our objective in the next two sections is to establish Theorem \ref{semipos} stated and used above .

\begin{notation}
\end{notation} We consider thus a prepared holomorphic fibration
$g:Y\ra Z$, with $Y$ and $Z$ smooth and $Z$ projective.

That $g$ is prepared means that its non-smooth locus is contained 
in a simple normal crossing divisor of $Z$, and that the inverse image by $g$ 
of this non-smooth 
locus is also a divisor of simple normal crossings on $Y$.
We let $\Delta':=\Delta(g)$, and $\Delta:=\Delta(g,H)$ be the 
orbifold divisors for $g$, and $(g,H)$ respectively, the second 
one is defined in \ref{compfib} to which we refer for details.

Thus $\Delta':=\Delta(g)=(\sum _{i\in I}(1-1/ m_i).\Delta_i)$ also is
 of simple normal
crossings. We define also $K_{Z/\Delta'}=K_Z+\Delta'$, and
$K_{Y/(Z/\Delta')}$
:=$K_Y\otimes g^*(K_{Z/\Delta'})^{-1}$. 

We apply similar notations for $\Delta$. Notice that $\Delta\geq \Delta'$.

\begin{remark} Notice that $\Delta(g, H)=\Delta(g, H^{vert})$, 
where  $H^{vert}$ is the {\bf $g$-vertical part of $H$},
equal to $\sum_{k\in K'} (1-1/n_k). H_k$, where $K'\subset K$ 
consists of the components of
$\mid H \mid$ which are not mapped onto $Z$ by $g$.

We shall also denote by $H^{hor}$ the {\bf $g$-horizontal part of 
$H$}, defined such that $H=H^{vert}+H^{hor}$.
\end{remark}

Our objective in this section is to establish the following
orbifold generalisation of famous results of Y.Kawamata anf 
E.Viehweg, initiated by T.Fujita 
in the case where $Z$ is a curve:

We shall actually essentially just reduce our case to the cases they treated.

\begin{theorem}\label{semipos'} Let $g:Y\ra Z$ be as above, prepared.
Let $H$ be an orbifold structure on $Y$, and $\Delta:=\Delta(g,H)$ as above. 
Assume $Z$ is projective and $Y$ is K\" ahler compact.
Let $m>0$ be an integer such that all $\mathbb Q$-divisors involved are integral.

Then: $g_{*}(m.(K_{(Y/(Z/\Delta))}+H)+B)$ is weakly positive, for $B$ a 
suitable effective $g$-exceptional divisor on $Y$.

\end{theorem}

We shall obtain this result as the consequence of three intermediate steps. 

The first step is a generalisation of the standard weak positivity results for direct image 
sheaves of pluricanonical forms in K\" ahler geometry: 

\begin{theorem}\label{dsemipos} Let $g:Y\ra Z$ be a fibration as in \ref{semipos'} above, with 
the same assumptions. Let $D=\sum_{j\in J}d_j.D_j$ be a divisor with integer positive coefficients on $Y$, the 
$D_j$ being (as always) pairwise distinct. Assume also that the support of 
$D$ is a divisor of simple normal crossings, and that $D$ is horizontal, which means that each $D_j$ is mapped 
surjectively onto $Z$ by $g$. 

Let $m$ be a positive integer such that $m\geq d_j, \forall j\in J$. 
Then: $g_*(mK_{Y/Z}+D)$ is weakly positive.
\end{theorem}

This result will be proved in the next section. Observe that the restriction that $D$ is horizontal is obviously 
unnecessary, because in general, one can remove the vertical part of $D$, if any, and then get an injection from the 
direct image sheaf relative to the horizontal part of $d$ into the same sheaf, relative to all of $D$. Notice that 
there is then no restriction needed on the multiplicities of the vertical components of $D$, besides nonnegativity.

The second and third steps are given by the following lemma \ref{basechange} and proposition \ref{orbbasechange}:

\begin{lemma}\label{basechange} Let $g:Y\ra Z$, $D$ and $m$ satify the same assumptions as the 
preceding \ref{dsemipos}. If $v:Z'\ra Z$ is a flat finite map with $Z'$ smooth, and if $g':Y'\ra Z'$ 
is deduced from $g$ by smoothing the base change $\hat{Y} :=Y\times_{Z}Z'$ of $Y$ by $v$, there is a natural injection of sheaves: 

$g'_*(mK_{Y'/Z'}+D')\ra v^*(g_*(mK_{Y/Z}+D))$, where $D':=u^*(D)$, if $u:Y'\ra Y$ is the natural map obtained by composing 
the desingularisation $d':Y'\ra \hat{Y}$ with the base change map $\hat{u}:\hat{Y}\ra Y$.

\end{lemma}

\begin{proposition}\label{orbbasechange} Let $g,D,m$ be as in \ref{basechange} above. Let, in addition, 
$H=H^{vert}$ be an orbifold vertical divisor on $Y$, which means that no component of $H$ is mapped onto $Z$ by $g$. 

There exists a finite flat map $v:Z'\ra Z$, with $Z'$ smooth, such that if $g':Y'\ra Z'$ is constructed as in \ref{basechange} 
above from $v$, the above injection of sheaves extends to an injection: 

$g'_*(mK_{Y'/Z'}+D')\ra v^*(g_*(mK_{Y/(Z/\Delta(g,H)}+H)+D+B)$, for some effective $g$-exceptional divisor $B$ on $Y$.
\end{proposition}

We shall give below the proofs of the preceding lemma and proposition. Let us first show that they imply, together 
with \ref{dsemipos}, the Theorem \ref{semipos'}: write $mH=D+mH^{vert}$, with $D:=mH^{hor}$. This is an integral divisor, 
if $m$ is sufficiently divisible. Moreover,  $0<d_j:=(1-1/m_{j}).m<m$, for each $j$. So that \ref{dsemipos} applies to 
$g'$ and $D"$, if $g'$ is deduced from $g$ by any base change $v:Z'\ra Z$ as in Lemma \ref{basechange} above, and if 
$D"$ is the strict transform by $u:Y'\ra Y$ of $D$, so that it is the horizontal part (for $g'$) of $D':=u^*(D)$. 
Actually, by the obvious injection of sheaves, 
the same result applies also if $D':=u^*(D)$, because then only effective vertical divsors are added. 

If we now apply the Proposition \ref{orbbasechange}, we see that the conclusion of \ref{semipos'} holds.

We shall now prove \ref{basechange} and \ref{orbbasechange}.

{\bf Proof of \ref{basechange}:} We use the notations of \ref{basechange}. It is proved in [Vi83], Lemma 3.3, pp. 335-336 to 
which we refer, that in our situation, we have: $d_*(mK_{Y'/Z'})$ naturally injects into $\hat{u}^*(mK_{Y/Z})$. 

This statement implies that $d_*(mK_{Y'/Z'}+D')$ injects into $\hat{u}^*(mK_{Y/Z}+D)$, since $D'=u^*(D)$. 

We then just need to apply $\hat{g}_*$ to both sides, noticing that $\hat{g}_*(\hat{u}^*)=v^*(g_*)$, by flatness of $v$. 

Here $\hat{g}:\hat{Y}\ra Z'$ is, of course, deduced from $g$ by the base change $v$ $\square$\\

{\bf Proof of \ref{orbbasechange}:}  We start with the construction of $v$: 

\begin{definition}\label{dnice} Let $g$ as above. A
finite covering $v:Z'\ra Z$ is said to be ${\bf\widetilde{\Delta'}}${\bf -nice}
(in similarity to [Kw]) if :

(1) It is flat, and $Z'$ is smooth.

(2) $v^*(\Delta_i)=\widetilde {m_i}.\Delta_i"$, for some reduced divisor
$\Delta_i"\subset Z'$, this for any $i\in I$.

Here $\widetilde {m_i}$ is any integer divisible by lcm($m_{ij}$'s), 
these being the same as above, used to define the multiplicity of $g$ along $\Delta'_i$.
Observe in particular that $v^{*}(\Delta')$ is Cartier on $Z'$.

(3) $v^{-1}(\mid \Delta' \mid)$ is a divisor of normal crossings on $Z'$.
\end{definition}

By [Kw] and [Vi82], such coverings exist.

{\bf Construction:} Let $v:Z'\ra Z$ be a $\widetilde {\Delta'}$-nice covering.
Let $u:Y'\ra Y$ be the composition $u=u"\circ n\circ d$, in which:
$u"=\hat{u}:Y":=\hat{Y}=\widehat {(Y\times_{Z} Z')}\ra Y$ is the base change by $v$, 
$n:Y^n\ra Y"$ is the normalisation, and: $d:Y'\ra Y^n$ is a desingularisation 
isomorphic above the smooth locus of $Y^n$. One thus also has, with the previous notations of \ref{basechange}: $d'=n\circ d$.

We also denote by $\hat{g}=g":Y"\ra Z'$ and $g':=g\circ d:Y'\ra Z'$ the maps 
deduced from $g$.

Then \ref{orbbasechange} is an immediate consequence of \ref{remwp} and the following:

\begin{proposition}\label{injsheav}: In the preceding situation, there exists on $Z'$ a
natural injection of sheaves:

   $u_{*}:(g')_{*}(m.(K_{Y'/Z'}+D'+(g')^*((v)^{^*}(m\Delta)))$
$\ra v^{*}(g_{*}(m.(K_{Y/Z}+H^{vert})+D+B))$, this for any integer 
$m>0$ sufficiently divisible,
and some effective $g$-exceptional divisor $B$ on $Y$.
\end{proposition}

{\bf Proof:} The main point is that we just need to check this
injection on the complement of a codimension two subvariety $A$ of $Z'$,
since by \ref{basechange}, the result holds with $(g')^*(v^{*}(m.\Delta))$ deleted from
the left-hand side, and $H^{vert}$ deleted from the right-hand side, on all of $Z'$.
   This provided $\mathcal O$$_Y(B)$ is defined by the poles of maximal order acquired by 
an arbitrary extension across $A$
of a section of $(g')_{*}(m.(K_{Y'/Z'}+(g')^*((v)^{^*}(m\Delta)))$ 
defined outside of $A$.

Observe that we are working here on the fibre product $Y\times_{Z}Z'$ 
which is Cohen-Mac Caulay by the normal crossing
hypothesis on $\Delta'$, so that the poles of the sections 
considered actually occur in codimension one. (This observation is due to S.Druel).

We shall thus check the above injection only above the generic point
$z$ of some $\Delta_i$. Let $U$ be a sufficiently small 
analytic open neighborhood of $z$ in $Z$. 

Let $Y_U:=g^{-1}(U)$,
and let $W\subset Y$ be a small analytic neighborhood of $y$, a generic
point of $D_{ij}$, notations being as in (1.1.2).

Thus $m_{ij}$ divides $\widetilde {m_{ij}}:=m_{ij}.q$, for some 
integer $q=q_{i,j}$, by the definition of
a $\widetilde {\Delta'}$-nice covering.

Factorise $v:=v'\circ v^*$ over $U$, with $v':Z'\ra Z^*$, and $v^*:Z^*\ra Z$,
in such a way that $v^*$ (resp. $v'$) ramifies at order exactly 
$m_{ij}$ (resp. $q$)
above $\Delta_i$ (resp. $\Delta_i^*:=(v^*)^{-1}(\Delta_i)$).

Construct $Y^*$ from $Y$ over $Z$ by taking base change by $v^*$, 
followed by normalisation and then smoothing.
We have also a natural fibration $g^*:Y^*\ra Z^*$. 

Possibly modifying 
$Y'$, we thus get a factorisation
$u=u^*\circ u'$. Moreover, from the usual commutation properties, we 
have: $v'\circ g'=g^*\circ u'$, and
$v^*\circ g^*=g\circ u^*$.

The crucial property in this construction is that $u^*$ is \'etale 
over $W$, if sufficiently small.
This is an easy standard local computation 
which we already used several times before (see [C2], for example). 
Thus $u^*$ is \'etale over the generic point
of $D_{ij}$.

 From the following lemma \ref{hvert} below, we deduce that, over the 
generic point of $D_{ij}$, we have
a natural injection of sheaves: $m.(K_{(Y^*/Z^*)})\subset 
(u^*)^*(m.(K_{Y/(Z/\Delta)}+H^{vert}))$, for
any integer $m>0$ sufficiently divisible.

From this injection, we can deduce the following, by tensorising with $\cal O$$_Y(D)$ 
and its lift $\cal O$$_{Y^*}(D^*)$ to $Y^*$ by $(u^*)$: $m.K_{(Y^*/Z^*)}+D^*\subset 
(u^*)^*(m.(K_{Y/(Z/\Delta)}+H^{vert})+D)$

The very same argument as in the proof of \ref{basechange} above shows the existence, for any $j$ (on which 
the preceding factorisations depend), of a natural
injection of sheaves: 

$g'_*(m.(K_{Y'/Z'})+D'))\subset 
(v')^*((g^*)_*(m.(K_{Y^*/Z^*})+D^*))$.

By composing the above injections, and restricting over $W$, we see 
that the sections of:

$(g^*)_*(m.(K_{Y^*/Z^*})+D')$ actually belong to 
$(v^*)^*(g_*(m.(K_{Y/(Z/\Delta)}+H^{vert})+D))$.

(We identify
local sections of $g'_*(m.(K_{Y'/Z'})+D')$ over $Z'$ and sections of 
$(m.(K_{Y'/Z'})+D')$ over corresponding open subsets of $Y'$).

The conclusion now follows from Hartog's extension theorem, applied 
over $Y_U$, to extend the sections thus obtained across the
intersection of two or more such $D_{i,j}$'s $\square$

We used the following:

\begin{lemma}\label{hvert}  With the above notations,  over the generic 
point of $D_{ij}$, we have
a natural injection of sheaves: $m.(K_{(Y^*/Z^*)})\subset 
(u^*)^*(m.(K_{Y/(Z/\Delta)}+H^{vert}))$, for
any integer $m>0$ sufficiently divisible.
\end{lemma}

{\bf Proof:} We shall argue using, instead of sheaf injections, 
rather inequalities between $\Bbb Q$-divisors;
the inequality $A\geq B$ meaning as usual that $(A-B)$ is effective.

From the equalities:

  $K_{Y^*}=(u^*)^*(K_Y)$ (because $u^*$ is \'etale), and:

  $K_{Z^*}=(v^*)^*(K_Z+ (1-1/m_{i,j}).\Delta_i)$
(ramification formula), we deduce that:

$K_{Y^*/Z^*}=(u^*)^*(K_{Y/Z}-(m_{i,j}-1).D_{i,j})$, because 
$g^*(\Delta_i)=m_{i,j}.D_{i,j}$,
and $g\circ u^*=v^*\circ g^*$.

On the other hand:

$K_{Y/(Z/\Delta)}+H^{vert}=K_{Y/Z}-g^*((1-1/m'_i).\Delta_i)+(1-1/n_{i,j}).D_{i,j}$ 

=$K_{Y/Z}+(-m_{i,j}.(1-1/m'_i)+(1-1/n_{i,j})).D_{i,j}$

=$K_{Y/Z}+(-m_{i,j}+1+(m_{i,j}/m'_i)-(1/n_{i,j})).D_{i,j}\geq 
K_{Y/Z}-(m_{i,j}-1).D_{i,j}$,
since: $m_{i,j}.n_{i,j}\geq m'_i$, by the very definition of $m'_i$.

This concludes the proof, by applying $(u^*)^*$ $\square$

\subsection{\bf Twisted Weak Positivity}\label{twwp}

The aim of this section is to prove Theorem \ref{dsemipos}. For this, we shall refer to the proof 
given in [Vi83] of the classical case where $D$ is empty, and simply indicate the changes needed. 
I would like to heartily thank E.Viehweg, who gave me a decisive 
hint for the proof of the following lemma. 

The single change needed lies in the corollary 5.2 of [Vi83]. We restate this corollary in the form we need.

\begin{lemma}\label{viehw} Let $g:Y\ra Z$ be as above. Let $A^*$ be ample on $Z$, let $A=g^*(A^*)$.
 
Assume that $S^N(g_*(mK_{Y/Z}+D+mA))$ is generated by its global sections on some nonempty open subset of $Z$. 

Then: $g_*(mK_{Y/Z}+D+(m-1)A)$ is weakly positive on $Z$.
\end{lemma}

(One can thus withdraw one $A$ in the above situation. This is sufficient to conclude to the desired weak positivity. 
See [Vi83] for the ingenious details.)

{\bf Proof:} We shall closely follow the proof of the same corollary 5.2 (where $D$ is empty) in [Vi83]. 

The proof in fact reduces to add $D$ everywhere in an appropriate way.  

Observe also that the proof given there 
uses only the fact that $g_*(K_{Y/Z})$ is weakly positive when $Y$ is projective. The case when $Y$ is K\" ahler 
can be obtained from the different proof of this fact sketched in [Vi86], based on [Ko86] 
(see the references in [Vi86] for more details), because it uses only the Hodge-theoretic 
result due to P.Deligne that holomorphic forms with logarithmic poles on a compact K\" ahler manifold are $d$-closed.

Take now $L:= K+D+A$, with $K:=K_{Y/Z}$. Define next;

$M:=$Image $(g^* (g_*(mK+D+mA))\ra (m(K+D+A)))$

Assume that the base locus of $(NmL-(m-1).ND)$ does not contain any component of $D$. We can of course always easily reduce 
to this case, by diminishing the relevant $d_j$'s, without modifying the conclusion.

We now reproduce the proof of [Vi83,5.2]. Blowing up $Y$ if needed, we can assume that $M$ is a line bundle, 
and that $NmL=NmM+E+(m-1)ND$, for $E$ an effective divisor on $Y$, not 
containing any component of $D$, such that $E+D$ has a support of normal crossings. 
By hypothesis, $Nm.M$ is generated by global sections over a nonempty Zariski open subset of 
$Z$.

 We first treat the case in which $m>d_j,\forall j\in J$. 

The corollary 5.1 of [Vi83] then applies without any change and shows that the subsheaf:
 
$g_*(K+L^{(m-1)})\subset g_*(K+(m-1)L)$ is weakly positive, in the notations $L^{(m-1)}$ of [Vi82,(2.2)]. 

Here $L^{(m-1)}:=(m-1).L-[((m-1)/mN).(E+(m-1)ND)]$, the integral part of a $\mathbb Q$-divisor being computed just 
by taking the integral part of the coefficients, componentwise. 

From the constructions made, we see that $g_*(K+L^{(m-1)})\subset g_*(K+(m-1)L-[((m-1)^2/mN).D])$, and 
that they coincide over the generic point of $Z$. 

The sheaf on the right is thus weakly positive, too. 

Computing, we get:  

$K+(m-1)L-[((m-1)^2/mN).D]=K+(m-1)L-\sum_{j\in J}[(m-1)^2m_j/m].D_j$

$=mK+\sum_{j\in J}((m-1)m_j-(m-2)m_j+[m_j/m]).D_j=mK+D$, as desired. 

This ends the proof of the special case (all $d_j<m$) considered. 

When $d_j=m$ for certain components $\Delta_j$ of $D$, so that $D=D^-+m\Delta$, each component $D_j$ of $D^-$ having $d_j<m$,
 one just needs to replace $K$ by $K+\Delta$ in the above proof, 
using the fact that $g_*(K+\Delta)$ is weakly positive on $Z$, by [Kw,Thm 32] (which can be proved also using 
the above result of Deligne in the K\" ahler case) $\square$ 

Once this Lemma \ref{viehw} is obtained, the rest of the proof of the weak positivity of 
$g_*(mK+D)$ is the same as the one given in [Vi83] when $D$ is empty.

\section{\bf Geometric Consequences of Additivity}\label{geomcons}

\subsection{\bf Varieties with $\kappa$=0}

The second fundamental example of special manifold is given by the:

\begin{theorem} \label{k=0} Assume $X\in \mathcal C$ has: $\kappa(X)=0$. Then $X$ is special.
\end{theorem}

{\bf Proof:} This is an immediate application of Theorem \ref{cnmtheor}:
let indeed $f:X\ra Y$ be a fibration of general type. Then the results 
apply and give: $\kappa(X)\geq \kappa(X_y)+dim(Z)$. But $dim(Z)>0$, by 
hypothesis, and $\kappa(X_y)\geq 0$. Contradiction $\square$

\begin{remark} This theorem together with Theorem \ref{specnobog} shows 
that $X$ has no Bogomolov sheaf if $\kappa(X)=0$. Actually, more is expected: 
$\kappa^+(X)=0$ if $\kappa(X)=0$ (see [C95] for 
details and definitions). This equality is shown for 
$X$ projective (or K\" ahler, using Yau's results) if $c_1(X)=0$.
\end{remark}

\subsection{\bf The Albanese Map}

The following result uses indirectly only the easiest part 
of the Additivity theorem \ref{cnmtheor}. It nevertheless seems to belong to this section.

\begin{proposition}\label{albsurj} Let $X\in \mathcal C$ be
special. Let $\alpha_X:X\ra Alb(X)$ be the Albanese map of $X$. Then:
$\alpha$ is surjective, connected, and has no multiple fibres in
codimension one (ie: $\Delta(\alpha)$ is empty).

If $X$ is only w-special (see section \ref{wspec} below for this notion), 
then $\alpha$ is surjective and connected.
\end{proposition}

{\bf Proof:} Assume $\alpha:=\alpha_X$ is not onto. 

Let $Z\subset Alb(X)$ be
its image. After ([U],(10.9)) , there exists a fibration
$g:Z\ra W$ with $W$ of general type (thus $dim(Y)>0$). Let $\psi:X\ra
W'$ and $\sigma :W'\ra W$
be the Stein factorisation of $\phi \circ \alpha$=$\sigma\circ \psi$. Then $W'$
is of general type, too, since $\sigma$ is finite. This contradicts
the fact that $X$  is special. Thus $\alpha$ is onto, and $Z=Alb(X)$.

Let now $\alpha=\beta\circ \alpha '$ be the Stein factorisation of
$\alpha$, with $\alpha ': X\ra A'$ connected and $\beta : A'\ra Alb(X)$
finite.
A slight variation of the arguments of [K-V] shows that if $\beta$ is
ramified, there exists a fibration $\phi:A'\ra W'$, with $W'$ of
general type.
Considering $\phi\circ \alpha '$, we get as above a contradiction to the
fact that $X$ is special.
Thus $\beta$ is unramified, hence isomorphic, by the universal
property of $\alpha$, which is thus connected and onto.

Assume now that the fibration $\alpha : X\ra Alb(X)$ has
$\Delta : =\Delta(\alpha)$ nonempty. Let $\Delta '$ be any component of
$\Delta$.

There exists a connected submersive quotient map $q:Alb(X)\ra
A:=Alb(X)/B$, for some subtorus $B$ of $Alb(X)$, such that
$\Delta '=q^{-1}(D')$, for some big $\Bbb Q$-divisor $D'\subset A$.

Consider the fibration $f:=q\circ \alpha:X\ra A$: the support of
$\Delta(f)\subset A$ obviously contains $D'$.
We conclude from \ref{k>=0} that $f$ is of general type (since $K_A$ is trivial and
$D'$ is big on $A$). Again this contradicts the assumption of
$X$ being special. And concludes the proof if $X$ is special.

If $X$
is only w-special, the first two steps of the proof apply without any
change
to give the last assertion $\square$

\begin{question}\label{albfib} Let $X\in {\cal C}$ be a special manifold.
Are then the generic fibres of $\alpha_X$ special ?
\end{question}

One can easily show that this question has a positive answer if so does the conjecture $C_{n,m}^{orb}$.

\subsection{\bf Varieties of General Type}

Although the results of this section are easy consequences of the main result \ref{decconj} below, we give a 
direct proof.

\begin{theorem}\label{ess=n} For any manifold $X \in \mathcal C$, we have: 
$ess(X)=dim(X)$ if and only if $X$ is of general type.
\end{theorem}

In particular, $X$ is generically covered by a nontrivial family of special
submanifolds (the fibres of $c_X$) if $X$ is not of
general type. The case when $\kappa(X)\geq 0$ is clear, by theorem 
\ref{k=0} applied to the fibres of the
Iitaka-Moishezon fibration of $X$.
So the result applies nontrivially only when $\kappa(X)=-\infty$, and 
its proof gives the following:

\begin{theorem} \label{k=-infty} Let $X\in \mathcal C$ be such that: 
$\kappa(X)=-\infty$ and $dim(X)>0$.
Then $c_X:X\ra C(X)$ has a general fibre $X_y$ which is special of 
positive dimension, with:
$\kappa(X_y)=-\infty$.
\end{theorem}

The proofs are easy applications of Theorem \ref{cnmtheor}:

{\bf Proof of \ref{ess=n}} Assume first that $X$ is of general type. 
So the identity map $id_X$ of $X$ is a special fibration of general type.
It is thus the core by \ref{gtspec}.

Conversely, assume that $ess(X)=dim(X)$, or what is the same, that: 
$c_X=id_X$, and that $n:=dim(X)>0$, so that 
$X$ is not special.

There exists then a fibration of general type $f:X\ra Y$, with $m:=dim(Y)>0$.
We proceed by induction on $n>0$. So the assertion holds for 
manifolds of dimension smaller than $n$.

If $F$ is a general fibre of $f$, then $F$ is of general type.

Indeed: otherwise, we could
construct a relative core for $f$, by theorem \ref{relquot}, and 
deduce a nontrivial fibration on $X$ with
general fibre special, contradicting our initial hypothesis.

Thus $F$ is of general type. But now $f$ is a fibration of general type 
with general fibre itself of general type.

By theorem \ref{cnmtheor}, we get that $X$ is of general type, as claimed $\square$

{\bf Proof of \ref{k=-infty}:}
We get first from the preceding argument that 
$dim(X_y)>0$, otherwise $\kappa(X)=n$. Assume that $\kappa(X_y)\geq 0$. From
Theorem \ref{cnmtheor} again, we learn that $\kappa(X)=-\infty\geq 
dim(Z)+\kappa(X_y)\geq dim(Z)\geq 0$. Contradiction $\square$

In a similar way, we can get a very simple description of the next step: $ess(X)=dim(X)-1$.

\begin{theorem}\label{ess=n-1} Let $X\in \mathcal C$ be a manifold of dimension $n>0$. 
Then $ess(X)=(n-1)$ if and only if one of the following two cases occurs:

(a) $\kappa(X)=(n-1)$ and $J_X$ is a fibration of general type.

(b) $R(X)$ (the rational quotient of $X$) is of dimension $(n-1)$ and of general type.
\end{theorem}

{\bf Proof:} The fact that cases (a), (b) imply that 
$ess(X)=(n-1)$ was already shown in \ref{k=n-1}. 
So assume $c_X:X\mero C(X)$ has $(n-1)$-dimensional image. Let $F$ be the generic 
fibre of $c_X$: it is a special curve, hence rational or elliptic.

If $F$ is rational then $C(X)=R(X)$, by the fact that $r_X$ dominates $c_X$ (see 
 \ref{RCquotspec}). Thus $R(X)$ is $(n-1)$-dimensional. Moreover, we get 
from \ref{corealg} that $R(X)$ is Moishezon, and so also $X$ is Moishezon, by [C85]. 

Because $C(R(X))=C(X)$, by 
Theorem \ref{CR=C} (which applies because $X$ is Moishezon), we conclude 
that $c_{R(X)}=id_{R(X)}$, and so by the preceding Theorem 
\ref{ess=n}, we infer that $R(X)$ is of general type, as claimed.

Next, if $F$ is elliptic, then $c_X=J_X$, because $J_X$ dominates $c_X$, and is obviously 
dominated by $c_X$ if $\kappa(X)\geq 0$. We show this by induction on $n\geq 2$ (the case of 
curves is obvious).

Thus $X$ is non-special, by assumption. Let $f:X\mero Y$ be a 
fibration of general type, and let $F$
be its general fibre. By \ref{domin}, 
we have a factorisation $f=c_X\circ g$, for a 
certain fibration $g:C(X)\mero Y$. 

Moreover, 
because of the functoriality properties of the core, 
the restriction of $c_X$ to the general fibre $X_y$ of $f$ is 
the core of $X_y$. The induction hypothesis thus applies to 
$X_y$, and we get that $\kappa(X_y)=(d-1)$ if $d:=dim(X_y)$, 
and also that $J_{X_y}$ (which is the restriction of $c_X$ to 
$X_y$) is a fibration of general type.

We are thus in position of applying Theorem \ref{cnmtheor}, since 
$f$ is of general type. We get first that $\kappa(X)=\kappa(X_y)+dim(Y)=(n-1)$. 

Then, because the general fibres $X_y$ of $f$ have nonnegative Kodaira dimension, 
we get also from Theorem \ref{cnmtheor} that $J_X$ is of general type, since so is its 
restriction to the general fibre of $f$ $\square$

The next step, classification of cases when $ess(X)=(n-2)$ 
is similar (but some cases would require the general 
case of $C^{orb}_{gt}$).

\subsection{\bf The Decomposition Theorem}\label{decth}

This is the following asssertion, which motivates most of the present paper:

\begin{theorem}\label{decconj} Let $X$ be non special. Then $c_X$ is a
fibration of general type.
\end{theorem}

  Roughly speaking, this means that $X$``decomposes" into its
``special part" (the fibres of $c_X$)
and its``core" or ``essential part", the orbifold
$(C(X)/\Delta(c_X))$, which is either a point, or 
of general type. Hence the name.

The Decomposition theorem can also be restated in the following form: 

Any $X\in \mathcal C$ has a fibration both special and of general type. 
This fibration is unique, and it is the core of $X$. (Apply \ref{gtspec} to 
see the last two assertions).

\begin{remark} We already showed some cases of the preceding theorem \ref{decconj}:
\end{remark}

(1) Up to dimension $3$ as a byproduct of 
the description of the core given in low dimensions.

(2) When $X$ is of general type, and when 
$ess(X)=dim(X)$: in this last case, showing 
that $X$ is of general type is actually 
the same as showing the Decomposition 
theorem in this case. The crucial 
step was actually applying $C^{orb}_{gt}$. 

This is not an accident: we shall show 
right below that the Decomposition Conjecture is 
in fact a consequence of $C^{orb}_{gt}$.

(3) Exactly the same remarks apply to 
the classification of the cases in which $ess(X)=(n-1)$ or 
$(n-2)$. 

(4) We shall show below two consequences of the 
Decomposition theorem: $C(X)$ is Moishezon, and 
 meromorphic multisections of $c_X$ are of general type.

\begin{theorem} \label{calg} Let $X\in \mathcal C $, and let $a_X:X\ra 
Alg(X)$ be its algebraic reduction.

There exists a factorisation $\phi:Alg(X)\ra C(X)$ of $c_X=\phi\circ 
a_X$. Alternatively: $C(X)$ is Moishezon.
\end{theorem}

{\bf Proof:} By Theorem \ref{fibalgred}, the general fibre of $a_X$ is 
special $\square$

Let us now establish the following weak version of the Decomposition theorem 
(\ref{gtsect} shows that it is a consequence of this result):

\begin{proposition} Let $j:Z\mero X$ be a meromorphic map 
such that $c_X\circ j:Z\mero C(X)$ is surjective and
generically finite, $c_X:X\ra C(X)$
being, as usual, the core of $X\in \mathcal C$. 

Then $Z$  is a variety of general type. 
\end{proposition}

{\bf Proof:} By the functoriality of the core,
 and the surjectivity of $c_X\circ j$, we 
get the existence of a factorisation $c_j:C(Z)\mero C(X)$.
Because $c_X\circ j= c_j\circ c_Z$ is generically finite, so is 
$c_Z$. The claim thus follows from \ref{ess=n} $\square$

{\bf Proof of Theorem \ref{decconj}:} Proceed by induction on $d:=ess(X)$. When $d=0$, the 
result trivially holds.

So assume it does when $ess(F)<d$. Assume that $ess(X)=d$, and that 
$\kappa(C(X))\geq 0$.

By assumption ($d>0$), $X$ is not special. So let $f:X\mero Y$ be of 
general type. We have, by Theorem \ref{domin}
a factorisation $\psi:C(X)\mero Y$ of $f=\psi\circ c_X$.

By the characteristic property of the core, the restriction of $c_X$ 
to the general fibre $X_y$ of $f$ is the core of $X_y$
(Otherwise $X$ would contain special subvarieties strictly larger 
than the generic fibre of $c_X$, a contradiction).

By the induction hypothesis, we conclude that 
the restriction of $c_X$ to the general fibre $X_y$ of $f$ is a fibration of 
general type.

We conclude from Theorem \ref{cnmtheor} that $c_X$ itself is 
a fibration of general type, as claimed $\square$

\subsection{Finite \' etale Covers}

Let us indicate that the Decomposition theorem implies the invariance of 
$ess(X)$ under finite \'etale covers:

\begin{theorem}\label{essetal} Let $u:X'\ra X$ be a finite \'etale 
cover. Let $c_u:C(X')\ra C(X)$ be the induced map.
  Then $c_u$ is generically finite. 
In particular: $ess(X)$ is invariant under finite \' etale covers; 
and a finite \' etale cover of a special manifold is 
special. 
\end{theorem}

Remark that the proof of this seemingly easy statement requires here fairly deep tools. It were 
interesting to know if there is an easy proof of it.\\

{\bf Proof:} We can assume that $u$ is Galois, of group $G$. Due to 
its uniqueness,
the map $c_{X'}$ is $G$-equivariant for an appropriate action of $G$ 
on $C(X')$.
Let $h:C(X')\ra C'(X)$ be the $G$-quotient. We have natural maps 
$c'_X:X\ra C'(X)$ (by $G$-invariance),
and $v:C'(X)\ra C(X)$ (since the fibres of $c'_X$, as images of those 
of $c_{X'}$ by $u$ are special.

Moreover: $c_X=v\circ c'_X$, by the general properties of $c_X$.

By the Decomposition theorem, $c_{X'}$ is a fibration of general type. And so 
$c'_X$ is also of general type, by \ref{behKod}.
 From \ref{domin}, we infer the existence of a factorisation $w:C(X)\ra 
C'(X)$ such that $c'_X=w\circ c_X$.
Thus $C'(X)=C(X)$, and $h=c_u$ is a finite map $\square$

\subsection{\bf Essential and Bogomolov Dimensions}\label{essbog}

\begin{definition} For $X\in \mathcal C$, let $B(X)$:=max $\{ 
p> 0\mbox{, such that there exists
a Bogomolov sheaf }  F\subset \Omega^p_X \}$ ; 
if there is no 
Bogomolov sheaf on $X$,
define $B(X):=0$.
\end{definition}

From \ref{specnobog}, we see that $B(X)=0$ if and only if $ess(X)=0$. In 
general, we have:

\begin{corollary} For any $X\in \mathcal C$, we have: $ess(X)= B(X)$.
\end{corollary}

{\bf Proof:} Let $F$ be any Bogomolov sheaf of dimension $p>0$ 
on $X$. The associated fibration is 
thus of general type, and so dominated by $c_X$. Thus $ess(X)\geq p$ 

Conversely, because $c_X$ is a fibration of general type, the Bogomolov sheaf associated to it 
has $p=ess(X)$. Hence the equality $\square$

Let us say that if $F,G$ are  Bogomolov sheaves on $X$ that 
$F$ {\bf dominates} $G$ if the fibration $f$ (of general type) 
associated to $F$ dominates 
the fibration $g$ (of general type) associated to $G$.

\begin{remark} The Decomposition theorem may then 
be restated in the following form:

 For any $X\in \mathcal C$, $ess(X)=B(X)$, there exists on $X$ 
a unique maximum Bogomolov sheaf, and it is defined by the core of $X$.
\end{remark}

\subsection{\bf Fibre Products}\label{fibpr}

We answer here the question \ref{gtfibprod}

\begin{theorem}\label{gtfibpr} Let $X\in \mathcal C$, and $f:X\mero Y$,$g:X\mero Z$ be two fibrations 
of general type. Let $h:X\mero W$ be the connected part of the Stein factorisation of the product map 
$f\times g:X\mero Y\times Z$. 

Then $h$ is a fibration of general type.
\end{theorem}

{\bf Proof:} This is an immediate generalisation of the argument in 
the classical case (when $Y$,$Z$ are themselves of general type). Indeed: let 
$g':W\mero Z$ and $f':W\mero Y$ the fibrations such that $f=f'\circ h$ and $g=g\circ h$. 

By the equality (deduced from \ref{cnmtheor}): $\kappa(W,h)=\kappa(W_z,g'_z)+\kappa(Z,g=g'\circ h)=\kappa(W_z,g'_z)+dim(Z)$,
 we just need 
to check that the restriction of $g':X_z\mero W_z$ is of general type, for $z$ general in $Z$.

Observe the statement is bimeromorphically invariant. We thus assume that $Z$ is projective (it is certainly 
Moishezon, since $g$ is of general type). Take $S\subset Z$ be an intersection of generic hyperpane sections 
such that $(g')^{-1}(S):=T\subset W$ is mapped surjectively and generically finitely onto $Y$ by $f'$. 

Let $V:=h^{-1}(T)\subset X$. The restriction$f_V:V\mero Y$ of $f$ to $V$ has thus the restriction $h_V:V\mero T$
 of $h$ to $V$  as Stein factorisation. It is thus a fibration of general type, by \ref{gtsect'}. 

If $g':T\mero S$ has positive-dimensional fibres, the conclusion follows from \ref{restgt}. If not, this means that 
$g$ dominates $f$, in which case the assertion is trivial $\square$

\subsection{Construction of the Core as the Highest General Type Fibration}\label{alternc_X}

We sketch here a second, shorter, construction of the core. It is more abstract, avoiding the consideration of 
chains of special varieties. It is something like a reconstruction ``a posteriori", once one already knows
from previous more down-to-earth work, the existence of the core. 

\begin{theorem}\label{core'} Let $X \in \cal C$. Then $X$ admits a fibration both of general type, and special. 
This fibration is unique up to equivalence. It is the core of $X$.
\end{theorem}

{\bf Proof:} Everything is clear from \ref{domin}, once the existence is known. We know from \ref{gtfibpr} above, together 
with \ref{defgtred} that $X$ admits a unique fibration of general type $gt_X:X\mero GT(X)$ 
wich dominates all general type fibrations with domain $X$.

It remains to show that the general fibres of $gt_X$ are special. Assume not. From the construction of the 
relative $gt$-reduction in \ref{gtrel}, we deduce the existence of a factorisation $gt_X=h\circ g$, with 
$g:X\mero G$ and $h:G\mero GT(X)$ such that the restriction $g_a:X_a\mero G_a$ to the general fibre $X_a$, $a\in GT(X)$, 
is the $gt$-reduction of $X_a$. Thus $g_a$ is a fibration of general type, since $X_a$ is assumed to be nonspecial. 

But then, we get from theorem \ref{decconj} that $g:X\mero G$ is a fibration of general type. This contradicts the facts 
that $dim(G)>dim(GT(X))$, by hypothesis, and that $gt_X$ dominates any fibration of general type defined on $X$ $\square$

\section{\bf Geometric Conjectures}

We shall state in this section some conjectures or 
questions naturally arising from our previous considerations.

\subsection{\bf Orbifold Rational Connectedness}\label{orbrc}

We consider in this section an orbifold $(Y/D)$ with $Y\in \cal C$ compact smooth and $D$ an orbifold divisor 
of normal crossings on $Y$. We define several variants of rational connectedness, directly adapted from the non-orbifold case. 

If $g: (Y/D)\ra Z$ is a fibration, we define its orbifold base $(Z/\Delta(g,D))$ as in section \ref{compfib}. 
If $g$ is only meromorphic, we replace $(Y/D)$ by a terminal modification $\mu:(Y'/D')\ra (Y/D)$ such that 
$g':=g\circ \mu$ is holomorphic, and then define 
the orbifold base of $g$ to be that of $g'$. We then define $\kappa(Z,g,D)$ as the 
minimum of the Kodaira dimensions of the orbifold bases so obtained.

The orbifold general fibre of $g:(Y/D)\ra Z$ is defined as the orbifold $(Y_z/D_z)$, with $Y_z$ the general fibre of $g$, and 
$D_z:=i^*(D)$, where $i:Y_z\ra Y$ is the inclusion map.

An {\bf orbifold rational curve} on $(Y/D)$ is an irreducible rational curve $C$ on $Y$ such that 
$\kappa(\hat{C},K_{\hat{C}}+n^*(D))=-\infty$, where $n:\hat{C}\ra C$ is the normalisation map.

\begin{definition}\label{deforbrc} We say that $(Y/D)$ is: 

(1) {\bf uniruled} if it is covered by a family of curves, whose generic member is an orbifold rational curve.

(2) {\bf rationally connected} if two generic points of $Y$ can be joined by a chain of orbifold rationall curves.

(3) {\bf rationally generated} if, for any fibration $g:(Y/D)\mero Z$, the orbifold base of some holomorphic 
model of $g$ is uniruled. 

(4) {\bf $\kappa$-rationally connected} if, for any fibration $g:(Y/D)\mero Z$, $\kappa(Z,g,D)=-\infty$.
\end{definition}

We note, as usual, for short: R.C (resp. R.G; resp. k-RG) for rationally 
connected (resp. rationally generated; resp. $\kappa$-rationally generated). Notice that R.G was introduced 
in [C95] in the non-orbifold context, and that the notion of $\kappa$-rational generatedness is directly inspired 
from the definition of $\kappa^+$ in the same paper.Actually, one can even define: 

\begin{definition} For $(Y/D)$ as above, define: $\kappa_+(Y/D):=$ max$\{\kappa(Z/\Delta(g,D))\}$, where $g$ runs through 
all fibrations $g:(Y/D)\mero Z$ as above.
\end{definition}

Then $\kappa$-rational generatedness is defined by the equality: $\kappa_+=-\infty$.

\begin{remark} One can also introduce in the orbifold case many other invariants, based on the vanishing 
of holomorphic `` covariant" tensors, which conjecturally characterise rational connectedness. We do not do this 
here, because one needs first to define these notions in the orbifold case. So we introduced only the simplest one: $\kappa_+$.
\end{remark}

One could, for example, define $\kappa_+(Y/D)$:=max $\{\kappa(Y,F)\}$, where $F$ ranges over all rank one coherent subsheaves of 
$\Omega_Y(log(D))$. (See section \ref{orbcat} for this notion). Another variant of orbifold rational connectedness would be defined by : 
$\kappa_+(Y/D)=-\infty$.

Standard arguments that we do not need to reproduce easily show: 

\begin{proposition}\label{comprc} For any orbifold $(Y/D)$ as above, one has: 

(1) $\kappa(Y/D)=-\infty$ if $(Y/D)$ is uniruled.

(2) $(Y/D)$ is R.G if it is R.C anf is k-RG if it is R.G.
\end{proposition}

The reverse implications depend on the answers to the following two questions.

\begin{question}\label{questrc}

(1) Assume $\kappa(Y/D)=-\infty$. Is then $(Y/D)$ uniruled?

(2) Let $g:(Y/D)\ra \mathbb P^1$ be a fibration with generic orbifold fibre rationally connected, and such that 
$\kappa(\mathbb P^1/\Delta(g,D))=-\infty$. Does $g$ have a section, that is: an orbifold rational curve mapped isomorphically 
to $\mathbb P^1$ by $g$?
\end{question}

Of course, the second question is directly inspired from [G-H-S]. Again, standard arguments show that: 

If question (1) has a positive answer, then: R.G and k-R.G are equivalent.

If question (2) has a positive answer, then: R.C and R.G are equivalent properties.

\begin{remark} Question (1) seems quite delicate (if true): it is open in the non-orbifold case.  
It is open even in the surface case $(dim(Y)=2)$, and seems to present great 
subtlety. It is similar to the case handled in [K-MK], of log-Del Pezzo surfaces with 
a reduced boundary $D$. But their methods might apply to treat the above surface case.
\end{remark}

\begin{remark} Another natural question in this context is whether the 
``glueing Lemma" of [Ko-Mi-Mo] extends to this context: 
for example, can one add to an orbifold 
rational curve $C'$ in $(Y/D)$ sufficiently many ``free" orbifold rational curves if $(Y/D)$ is uniruled, 
in such a way that the union deforms to an (irreducible) orbifold rational curve? 
\end{remark}

We shall now justify the introduction of the notion of $\kappa$-generatedness, by showing that it permits to define 
the notion of rational quotient (even in the orbifold case) without solving question (1) above, this definition however depends 
on the solution of conjecture $C_{n,m}^{orb}$.

\begin{conjecture}\label{k=-inftyconj}$((\kappa=-\infty)$-Conjecture) Assume that $\kappa(Y/D)=-\infty$.

There exists then a unique fibration $r_{Y/D}:Y\mero R(Y/D)$ such that:

{\bf (1)} the general orbifold fibre $F_r:=(Y_r/D_r), r\in R(Y/D)$ of $r_{Y/D}$ has $\kappa_+(Y_r/D_r)=-\infty$. 

{\bf (2)} $\kappa(R(Y/D),r_{Y/D},D)\geq 0$.
\end{conjecture}

\begin{proposition}\label{rbq} A fibration as in \ref{k=-inftyconj} is unique if it exists. 
The existence follows from 
Conjecture $C^{orb}_{n,m}$.
\end{proposition}

{\bf Proof:} Assume there exists another fibration $h:Y\mero Z$ with the same properties. Let $u:Y\mero W$ be the 
Stein factorised product of $r_{Y/D}$ and $h$. Assume $h$ is not equivalent to $r_{Y/D}$. The restriction $h_r:Y_r\mero W_r$ then 
maps the orbifold $F_r$ to a positive-dimensional orbifold of nonnegative Kodaira dimension (by the property (2), 
assumed to hold for $h$). This contradicts (1) above for $r_{Y/D}$, and shows unicity. 

For the existence: consider the Stein factorised product $p:Y\mero V$ of all 
fibrations $h:Y\mero Z$ such that $\kappa(Z/\Delta(h,D))\geq 0$. By $C^{orb}_{n,m}$, $p$ itself has this property. 
Assume the general orbifold fibre of $p$ were not $\kappa$-RG. We could find (by the result on 
the relative quotients) fibrations $h:Y\mero U$ and $k:U\mero V$ such that $p=k\circ h$, and such that the general orbifold 
fibre of $k$ had nonnegative Kodaira dimension. By $C^{orb}_{n,m}$ again, this would imply that $\kappa(U/\Delta(h,D))\geq 0$. 

But this contradicts the maximality of $p$ among such fibrations, since $dim(U)>dim(V)$ $\square$

\subsection{\bf Orbifold Iitaka-Moishezon Fibration}\label{jorb}

The Iitaka-Moishezon fibration of an orbifold $(Y/D)$ as in section \ref{orbrc} above can be defined without difficulty 
by considering the linear system $\mid m.(K_{Y/D}\mid$ for $m$ sufficiently big and divisible. 

We shall then denote by $J_{Y/D}:(Y/D)\mero J(Y/D)$ the resulting Iitaka-Moishezon fibration. 

By the usual properties of Iitaka fibrations, one has: $\kappa(Y_s/D_s)=0$, for its general orbifold fibre.

\subsection{\bf Orbifold Rational quotient and Iitaka Fibration of a Special fibration} \label{orbRJ}

We describe in this section an alternative construction of the core 
of an arbitrary $X$,
along the program described in the introduction. This 
construction actually coincides in dimensions
2 and 3 with the construction of the core given in \ref{coresurf} and \ref{corenonspec3f}.

The idea is quite simple: assume we have a {\bf special} fibration $f:X\mero Y$. 
We try to understand the obstructions to constructing a non-trivial $h:Y\mero Z$ 
such that $g:=h\circ f:X\mero Z$ is still special, this in terms of the invariant 
$\kappa(Y,f)$. 

This will certainly be impossible 
if $f$ is of general type, because then $f$ is the core of $X$, by \ref{gtspec}. 

The other cases can actually be reduced (as in the classical Iitaka-Moishezon 
classification Program) to the cases $\kappa=0$, or $-\infty$. And we conjecture 
below that such a non-trivial $g$ then exists. 

This leads to a conjectural decomposition of the core in a canonical and functorial 
sequence of orbifold rational and Iitaka fibrations.

First, we need to give an orbifold generalisation 
of Iitaka fibrations and
rational quotients for special fibrations.

We consider in the sequel a {\bf special} fibration $f:X\mero Y$ 
(thus with $X\in \mathcal C$). 

\begin{theorem} \label{k=oconj} Assume that $\kappa (Y,f)=0$. Then $X$ is special.
\end{theorem}

{\bf Proof:} This follows from \ref{cnmtheor} and its corollaries. Assume indeed $X$ were not 
special. By \ref{domin}, there would exists a fibration $g:Y\mero Z$ such that $g\circ f:X\mero Z$ were of general type. 

But then, we would have: 

(1) $\kappa(Y_z,f_z)\geq 0$, because $\kappa(Y,f)=0$, by the easy addition theorem.

(2) $0=\kappa(Y,f)=\kappa(Y_z,f_z)+dim(Z)\geq dim(Z)>0$, because of (1). Contradiction $\square$

 From this result, we deduce immediately the notion of {\bf Iitaka fibration} of the special fibration $f:X\mero Y$: 
it is just the fibration $J_f:Y \ra 
J(f)$ associated to the $\mathbb Q$-divisor $(K_Y+\Delta(f))$, first replacing $f$ by any of
its admissible models. Alternatively, it is also the Iitaka-Moishezon fibration of the orbifold base $(Y/\Delta(f))$ 
of any of its admissible prepared models.

The important property is that the general orbifold fibre of $J_f$ has Kodaira dimension 
zero (by the remarks made in section \ref{jorb} above.

From this follows immediately, using \ref{k=oconj}, that the composite fibration $J(f)\circ f:X\mero J(f)$ is special.

Next we need to investigate the case when $\kappa(Y,f)=-\infty$. Using the preceding conjecture \ref{k=-inftyconj}, we get:

\begin{conjecture}\label{k<oconj} Let $f:X\mero Y$ be special fibration, with $\kappa(Y,f)=-\infty$. 
There exists a unique fibration $r_f:(Y/\Delta(f))\mero R(f)$ such that:

(1) The general orbifold fibre of $r_f$ is $\kappa$-rationally generated.

(2) $\kappa(R(f)/\Delta(r_f,\Delta(f))= \kappa(R(f),r_f\circ f)\geq 0$.

We call $r_f$ the {\bf rational quotient} of $f$.
\end{conjecture}

Notice that, assuming the existence of $r_f$, we have: 
$dim(R(f))<dim(Y)$, and that $f^{-1}(Y_r):=X_r$ is then special, for $r\in R(F)$ general.

 By a slight extension of the definition, we
also define the rational quotient $r_f$ of $f$ to be the  identity 
map of $Y$ if $\kappa(Y,f)\geq 0$.

{\bf Caution:} The Iitaka base fibration (resp. the rational 
base quotient) differ from
the (relative) Iitaka fibration (resp. rational quotient) of $f$, 
which are defined for any fibration $f$,
  and provide a factorisation of $f$, inducing on the general fibre of $f$ the corresponding fibration of that fibre.

Of course, if the questions \ref{questrc} have a positive answer, the general orbifold fibre of $r_f$ is rationally conected.

\subsection{\bf The Core as an Iterated Orbifold 
Rational Quotients and Iitaka Fibrations}\label{c=rj^n}

We assume in this section the conjecture \ref{k=-inftyconj}.

\begin{definition}{\bf Special Reduction of a Special Fibration} \label{specred}
Let $f:X\ra Y$
be a {\bf special} fibration.
Assume that the above conjecture \ref{k=-inftyconj} hold.
Define $s_f:J_{r_f}\circ r_f:Y\mero J(r_f\circ f):=S(f)$, and
call it {\bf the first special reduction} of $f$.
\end{definition}

In vague terms: $s_f$ is obtained by first applying $r_f$ to $Y$, and then $J(r_f\circ f)$ to $(R(f)/\Delta(r_f\circ f))$.ñ
Notice that $s_f:X\ra S(f)$ has still its general fibre special.
So the construction can be iterated, and we get a sequence of special 
fibrations:
$s^0_f:=id_Y$, with $Y:=S^0(f)$, and then inductively: 
$s^{k+1}_f:=s_{s^k_f}:X\ra S^{k+1}(f)$.
We call $s^k_f$ {\bf the $k^{th}$-special reduction of $f$}.

Notice that $s_f=f$ exactly when $\kappa(Y,f)=dim(Y)\geq 0$.

In this case, by \ref{gtspec}, $f$ is the core of $X$. And it is the 
constant map iff $\kappa(Y,f)=dim(Y)=0$.

 From this we immediately deduce the following, by applying the
above iteration starting with $f:=id_X$, the identity map of $X$,
thus getting the sequence of $k^{th}$-special reductions of $X$:

\begin{theorem} Assume that the conjecture \ref{k=-inftyconj} holds.
Let $X\in \mathcal C$, with $n:=dim(X)$.
Define $s^k_X:= s^k_{id_X}:X\ra S^k(X):=S^k(id_X) ,\forall k\geq 0$. Then
$s^n_X=c_X$, the core of $X$, and $c_X$
is thus either constant, or a fibration of general type.
\end{theorem}

Of course, in general, the sequence of fibrations $s^k_X$ will be stationary
before its $n^{th}$-term is reached. 

\begin{definition} For $X\in \mathcal C$ and $k\geq 0$, define:
$s^k(X):=dim(S^k(X))$, and $rs^k(X):=dim(R(s^k_X))$, where
$r_{s^k_X}:S^k(X)\ra R(s^k_X)$ is the rational quotient of 
$s^k_X:X\ra S^k(X)$,
the $k^{th}$-special reduction of $X$.
\end{definition}

Under the said conjecture, these are new bimeromorphic invariants of $X$. We call $s^k(X)$ the
{\bf $k^{th}$-special dimension} of $X$, and $rs^k(X)$ the
{\bf reduced $k^{th}$-special dimension of $X$}.

Standard easy arguments give:

\begin{proposition}  Assume that the conjecture \ref{k=-inftyconj} holds. 
The invariants $s^k(X)$ and $rs^k(X)$
are invariant under finite \'etale covers $u:X'\ra X$ of $X$. The same holds
for the associated fibrations: $s^k_{X'}$ (resp. $rs^k_{X'}$) is the Stein
factorisation of $s^k_X\circ u$ (resp. $rs^k_X\circ u$ ).
\end{proposition}

One has the following ``strictness" property for the core decomposition, which we quote without proof:

\begin{theorem}\label{corestrictness} Assume that conjecture \ref{k=-inftyconj} holds. Let $X\in \cal C$, and let $k\geq j\geq 0$ 
be integers. Let $S$ (resp. $R$) be the general fibre of $s_X^k$ (resp. $rs_X^k$). The restrictions of $s_X^j$ and $rs_X^j$ to 
$S$ and $R$ coincide respectively with $s_S^j$, $rs_S^j$, $s_R^j$ and $rs_R^j$.
\end{theorem}

\begin{remark}

The sequence (s) of $(2n+2)$ integers 
$dim(X)=n=s^0\geq rs^0\geq s^1\geq rs^1\geq...s^k\geq rs^k\geq ...s^n\geq rs^n$ 
partitions the class of $n$-dimensional manifolds in $\mathcal C$ 
into a certain (finite) number $c(n)$ of classes, 
according to the number of steps and the dimensions 
of the steps needed to 
decompose the core in orbifold rational quotients and Iitaka fibrations.

One has for example: $c(0)=1$, $c(1)=3$, $c(2)=8$, $c(3)=21$, by direct listing.  
\end{remark}

One can easily compute the numbers $c(n)$ as follows. Say a sequence (s) is special if its last term is zero.

Notice the sequence is stationary precisely if two consecutive 
terms $rs^k=s^{k+1}$ coincide, or if one is zero.

We say the sequence (s) {\bf ends} at $J$ (resp. $r$) if $s^k=rs^k<rs^{k-1}$ (resp. $s^{k+1}=rs^k<s^k$).

let $s(n)$ be the number of special sequences of $n$-dimensional manifolds. 

Let $s'(n)$ (resp. $s"(n)$) be the number of such special sequences which end at $r$ (resp. $J$). 

Thus: $s(n)=s'(n)+s"(n)$. 

We quote without proof:

\begin{lemma}  One has, for any $n\geq 0$:

(1) $c(n):=\sum_{m=0,1,...,n}s(m)$,
 
(2) $s'(n)=1+\sum_{h=1,...,(n-1)}s"(h)$,

(3) $s"(n)=1+\sum_{h=1,2,...,(n-1)}s(h)$. 
\end{lemma}

Because $s(0)=1$,$s'(1)=s"(1)=1, s(1)=2, c(1)=3$, we get inductively: 
$s'(2)=2,s"(2)=3, s(2)=5, c(2)=8$, for example. 
We already met these classes, actually, when listing 
the possibilities for the core of a surface.

\begin{corollary} For all $n\geq 0$, one has:

(1) $s'(n+1)=s(n)$, $s"(n+1)=2.s(n)-s(n-1)$, 

(2) $s(n+1)=3.s(n)-s(n-1)$, $c(n+1)=c(n)+s(n+1)$.
\end{corollary}

And, finally:

\begin{corollary} If $\alpha:=(3+\sqrt 5 /2)$, and $\beta:=(3-\sqrt 5 /2)$, one has, for all $n\geq 0$:

(1) $s(n)=((1+\sqrt 5) /(2.\sqrt 5)).\alpha^n+((\sqrt 5-1)/(2.\sqrt 5)).\beta ^n$.

(2) $c(n)=(\alpha^{n+1}-\beta^{n+1})/\sqrt 5$.
\end{corollary}

\subsection{\bf  Deformation invariance}\label{definv}

We next state the following, for $X\in \cal C$ smooth:

\begin{conjecture}( Deformation Conjecture) $ess(X)$ is a deformation invariant of $X$.
In particular: the class of special manifolds is 
stable under deformation.
\end{conjecture}

One may ask more precise versions: does the core of $X$ vary 
holomorphically with $X$ for suitable models?

Even weaker versions are difficult: is $ess(X)$ upper or lower 
semi-continuous under smooth deformations?

Is it lower semi-continuous under degenerations?

The following is important and maybe accessible: is a degeneration of special manifolds 
still special (in the sense that all of its irreducible components are special)?

 A more refined version of the above deformation conjecture 
is the following, which extends the classical conjecture 
concerning the
deformation invariance of the Kodaira dimension for compact K\"ahler manifolds:

\begin{question}  Assume that the $(\kappa=-\infty)$-Conjecture hold. 
Let $X\in \mathcal C$. Are the invariants
  $s^k(X)$ and $rs^k(X)$ deformation invariants of $X$? More precisely, do
the associated fibrations vary holomorphically (for suitable models) 
when deforming $X$?
\end{question}

Notice that this is the case for $k=0$, and
in the projective case (and deformation among projective manifolds) 
for $s^1_X$,
by the invariance of plurigenera, due to Y.T.Siu ([Siu']).

A slightly less general version, but which deals with invariants known to exist, 
is the following, concerning the higher Kodaira dimensions, defined in \ref{higherkoddim}: 

\begin{question} Let $X$ be compact K\" ahler. Assume the invariants:
  
$\kappa(X)\geq\kappa'(X) \geq ...\geq \kappa^r(X)\geq -\infty$ are defined. 

Are then these invariants defined also and constant for small deformations of $X$?
\end{question}

Of course, one could still extend the preceding questions 
to the category of orbifolds, appropriately defined.

\subsection{\bf Weakly-Special Manifolds}\label{wspec}

\begin{definition} $X\in \mathcal C$ is said to be {\bf weakly-special}
(or w-special for short), if there
is no pair $(u,f')$ in which $u:X'\ra X$ is a
finite \'etale covering, and $f':X'\mero Y'$ is a surjective meromorphic
map onto a manifold $Y'$ of
general type (in the usual sense) and positive dimension.
\end{definition}

  From \ref{behKod}, we get:

\begin{proposition}  If $X$ is special, it is w-special. 
\end{proposition}

{\bf Proof:} It is the same as in \ref{gtfibprod}. Assume by contradiction that a triple 
$(u,X',f')$ as in the preceding definition exists. We can assume the cover $u$ to be Galois, 
of group $G$. If $f'$ is $G$-equivariant, then $f'$ descends to $f:X\mero Y:=(Y'/G)$, and $f$ is of 
general type by \ref{behKod}. Otherwise, just replace $f'$ by the least upper bound $f"$ of the (finite) 
family $(f'\circ g), g\in G$.  Now $f"$ is $G$-equivariant, and maps $X'$ to $Y"$, which is of general 
type, by \ref{fibprodgt} $\square$

\begin{example} All the examples listed in \ref{specex} are
thus w-special.
\end{example}

Notice that, contrary to the case of special manifolds, it is
obvious that any
finite \'etale cover of a w-special manifold is
   again w-special. In particular, any finite \'etale cover of a special manifold
is w-special.

Conversely, it follows from \ref{specsurf} that special surfaces and are weakly special.

The same assertion for threefolds is expected to have a negative answer. More generally:

\begin{question} Does there exist manifolds which are
w-special, but not special?
\end{question}

Although it seems reasonable to expect a positive answer,
the construction of such examples is far from being immediate.

The following has been shown in \ref{albsurj}:

\begin{proposition} Let $X$ be w-special. The Albanese
map of $X$ is then surjective and connected. 

In particular: if $X$ is special, and $X'$ is a finite \' etale cover of $X$, the
 Albanese map of $X'$ is surjective and connected.
\end{proposition}

Observe that, in contrast to the case $X$ special, 
this result does not assert that the Albanese map of $X$ is multiplicity free if 
$X$ is w-special.

Thus if $X$ is w-special, but its Albanese map has
multiple fibres in codimension one, $X$ cannot be special.

The core as a weak analogue, described below.

\subsection{\bf The Weak Version of the Core}\label{wcore}

It is possible to construct an analogue of the core for 
weakly-special manifolds and
fibrations strongly of general type. We briefly indicate how to do 
this. The proofs are much easier
than in the case of the actual core, and will be omitted, because of being 
straightforward applications
of standards results concerning the classical (not orbifold) Kodaira dimension. 
These notions 
will not be used, anyway, but it seems worth mentioning them, for sake of completeness.

\begin{definition} We say that the fibration $f:X\mero Y$, with 
$Y\in \mathcal C$, is of {\bf strong
general type} if there exists a finite \'etale cover $u:X'\ra X$ such 
that $Y'$ is of general type,
where $f\circ u= v\circ f'$ is the Stein factorisation of $f\circ u$. 
Here $f':X'\mero Y'$ is connected and
$v:Y'\ra Y$ is finite.
Thus $X$ is weakly-special iff it has no fibration $f:X\ra Y$ of 
strong general type.
\end{definition}

Notice that, by \ref{behKod}, a fibration of strong general type is of general type.

\begin{proposition}  Let $f:X\ra Y$ and $g:X\ra Z$ be two 
fibrations of strong general type.
Let $h:X\ra V$ be the Stein factorisation of the product map $f\times 
g:X\ra W\subset Y\times Z$.
Then $h$ is of strong general type.
\end{proposition}

\begin{corollary} For $X\in \mathcal C$, there exists a 
unique fibration $gt_X^+:X\ra GT^+(X)$ which is
of strong general type and dominates (in the obvious sense) any 
fibration of strong general type $f:X\ra Y$.
We call it the strong general type reduction of $X$. It is thus the 
constant map iff $X$ is weakly-special.
\end{corollary}

\begin{corollary} For $X\in \mathcal C$, there is a unique 
$w_X:C(X)\ra GT^+(X)$
such that $gt_X^+=w_X\circ c_X$, if $c_X:X\ra C(X)$ is the core of $X$.
\end{corollary}

{\bf Proof:} By \ref{domin}, if $f:X\ra Y$ is a fibration of general 
type, there is a unique
factorisation $\theta:C(X)\ra Y$ such that $f=\theta\circ c_X$. Apply 
this to $f=gt_X^+$, which is of general
type, since of strong general type.

\begin{remark}
\end{remark}

(1) The general type reduction $gt_X:X\ra GT(X)$ of $X$ (see \ref{defgtred}) 
dominates
$gt_X^+$ (and is dominated by $c_X$). 

(2) For any fibration $f:X\ra Y$, it is possible to construct a 
relative strong general type reduction
of $f$, as in the general type case.

(3) Applying the above reduction to $f:=gt_X^+$, we would get
that the general fibre of $gt_X^+$ is weakly-special,
provided one could show the following : if $g\circ f:X\ra Z$
is of strong general type, and if so is the restriction $f_z:X_z\ra 
Y_z$  of $f:X\ra Y$
to the general fibre $X_z$ of $g\circ f$,
then $f:X\ra Y$ is of strong general type.

This statement seems however very unplausible 
(it would be much better to have counterexamples!), 
which is why we were led
  to consider the more complicated notion of
special varieties and general type fibrations.

(4) As D.Abramovitch informed me, in [Ab] a reduction map analogous to the preceding 
$gt_X^+$ was constructed, but with fibrations onto general type varieties in place of strong general 
type fibrations (ie: ignoring \' etale covers),
inspired by ``Harris conjecture". It enjoys nice properties, but is 
very unstable with respect to
finite \'etale covers (which however preserve Kobayashi pseudo-metric 
and the arithmetic behaviour).

\subsection{\bf Classical $gcd$-Multiplicities}\label{gcdmult}

We can still introduce a third notion of special manifolds, which interpolates between the preceding 
two ones. It is based on the classical definition of multiple fibres, using 
$(gcd)$ instead of $(inf)$. We shall be brief on this.

\begin{definition} Let $f:X\ra Y$ be a fibration, and let $\Delta\subset Y$ 
be an irreducible
reduced divisor. Let  $f^*(\Delta)=:(\sum_{j\in J} m_j.D_j)+R $, 
where $f(D_j)=\Delta,\forall j$, while
$f(R)$ has codimension two or more in $Y$. Define multiplicities as 
follows: 

$m^-(\Delta,f):=gcd(m_j,j\in J)$,
and $m(\Delta,f):=min(m_j,j\in J)$. Obviously, $m^-$ divides $m$.
\end{definition}

Now using this definition of multiplicity, one can define, exactly 
as we did, with the same proofs, the notion of base orbifold and 
Kodaira dimension of a fibration. This leads to the ($gcd$-versions) 
of fibration of general type, and special manifold. The $gcd$-core 
can be constructed also with the same properties. 

One property only is 
(possibly) lost: $gcd$-fibrations of general type 
are no longer naturally in bijective 
correspondance with Bogomolov sheaves. 

On the other hand, the additivity result \ref{cnmtheor} and its proof 
(even becoming slightly simpler at a point) remain. So it may happen 
that the $gcd$-core is still a fibration of $gcd$-general type. 

At this point, one may wonder:

(1) Are the two theories actually different?

(2) If yes, which one is the ``right" one?

Presently, there is no counterexample to distinguish them. But it seems 
plausible that the two theories differ. We shall make some 
general easy comparisons between them right below.

Concerning the second question (b) above, assuming a positive answer to (a), it seems 
that the version given here (with $inf$, not $gcd$) is the right one. This is supported by 
the correspondance with Bogomolov sheaves, and the direct link between 
$inf$ fibrations of general type and geometric positivity properties of cotangent sheaves. 

This link between 
geometry and positivity is the main theme of Classification theory.

Other impressions support this choice of $inf$ versus $gcd$: the $inf$-notion 
is better suited for the study of Kobayashi pseudo-metric. Also, one 
may expect on the arithmetic side the fact that $gcd$-special manifolds are not always 
potentially dense. 

The only feature in favour of the $gcd$-theory is that it is better suited for the 
computation of the fundamental group. This property cannot however compensate the other drawbacks.

Let us end up with some obvious remarks concerning the comparisons between 
the three notions of specialness introduced:

A special manifold is $gcd$-special (because a fibration of general type is of
$gcd$-general type).

A $gcd$-special manifold is w-special (by the $gcd$-version of \ref{behKod}).

Although no counterexamples are known, it seems plausible that the three notions are distinct. 

The easiest counterexample would be a w-special, but not special threefold.

The following question is related to the comparison of special and 
$gcd$-special manifolds. A positive answer (not necessarily expected) would show 
the coincidence of the two notions):

\begin{question} Assume the general fibre $F$ of $f$ is special. 
Is it then true that $m(\Delta,f)$=
$m^-(\Delta,f)$, for any pair $(\Delta \subset Y)$?
\end{question}

\begin{remark}
\end{remark} 

It is sufficient to answer the above question 
when $Y$ is a curve and $m^-(\Delta,f)=1$:
one would like to know if one of the $m_j's=1$, then.
The question above is motivated by the fact that it has a positive 
answer when $F$ is rationally connected (by [G-H-S], the above assertion being actually 
the main difficulty in the proof),
  and when $F$ is an abelian variety (or a complex torus), by an easy argument
(the addition of local multisections). The next interesting case to consider is 
that of $F$ a $K3$-surface or a special surface with $\kappa=1$. The answer does not seem
to be known in these cases.

\subsection{\bf Irreducible Special Manifolds}\label{irrspec}

\begin{definition} Let $X\in \mathcal C$ be a special connected manifold.
 We say 
$X$ is {\bf irreducible special}
(resp. {\bf strongly irreducible special}) if it is special, and if 
every fibration $f:X\ra Y$
(resp. covering family $((V_t)_{t\in T})$ of $X$) with general fibre 
special (resp. special general member $V_t$) is trivial
(ie: has fibre or base dimension zero); (resp. if $V_t$ has dimension 
or codimension zero in $X$).
\end{definition}

{\bf Examples:} If $X$ is rationally connected of dimension 
$2$ or more, it is not strongly irreducible special.

But it is possible that some rationally connected threefolds 
(necessarily bimeromorphic to some terminal $\Bbb Q$-Fano threefolds
with Picard number one) may be irreducible special. The answer does 
not seem to be known.

The general quintic threefold is strongly irreducible (because not 
generically covered by images of elliptic curves or abelian surfaces).

\begin{remark} If $X$ is irreducible special, then 
$\kappa(X)=0$, or $-\infty$. (By \ref{k=n}, \ref{k=0}, and the Iitaka-Moishezon 
fibration).
\end{remark}

Recall the standard:

{\bf Uniruledness Conjecture:} If $\kappa(X)=-\infty$, then $X$ is uniruled.

Observe that it is obviously equivalent to:

{\bf Irreducibility Conjecture:} If $\kappa(X)=-\infty $, and if $X$ is 
special, of dimension $n\geq 2$, then $X$ is not strongly irreducible special.

\subsection{\bf Manifolds with nef anticanonical bundle}

   Question: Let $X$ be a K\"ahler manifold with $-K_X$ nef (see
[D-P-S] for this notion). Is then $X$ special?
   A positive answer to this
question would imply most of the conjectures concerning these manifolds.

Certain important consequences of a positive answer to this question 
have been established:
Q.Zhang has shown in [Zh] the surjectivity of the Albanese map in the projective case.

And M.Paun has shown (in the difficult K\" ahler case, see [Pa]) the virtual 
nilpotency of the fundamental group.
  Combined with Q.Zhang's result and [C4], it shows that if $X$ is 
projective and has nef anticanonical
bundle, then it has virtually abelian fundamental group.

\section{\bf The Fundamental group}

\subsection{\bf The Abelianity conjecture}

Our considerations here are guided by the following:

\begin{conjecture} \label{abconj}Let $X$ be special. Then  $\pi_1(X)$ is
virtually abelian.
\end{conjecture}

   Recall that a group is said to be {\bf virtually abelian} (one also
says: almost abelian) if it has a subgroup of finite index which is
abelian.

This conjecture is motivated by the fact that it should be true for klt orbifolds  
either rationally connected or with $\kappa=0$, and the fact that a special manifold 
should be (after section \ref{c=rj^n}) a tower of orbifold fibrations with orbifold fibres 
in the preceding two classes. See section \ref{motconjd} for the motivation of a similar conjecture, 
concerning the Kobayashi pseudo-metric.

\begin{examples}\label {exabel}
\end{examples}

{\bf (1) Curves:} A curve is special iff its fundamental group
is virtually abelian.

{\bf (2) Rationally connected manifolds:}
They have trivial fundamental group
([C3],[C5],[Ko-Mi-Mo], and are special.

More generally, if $r_X$ is the rational quotient of $X$ 
(see \ref{ratquot}), 
then it induces an isomorphism between
the fundamental groups of $X$ and $R(X)$ by [Ko]. So that $X$ and
$R(X)$ are simultaneously special and have simultaneously virtually
abelian fundamental groups. The above conjecture should  thus essentially
reduces to the case where $\kappa(X)\geq 0$ if 
the uniruledness conjecture is true.

{\bf (3) Manifolds with $a(X)=0$ or with $\kappa(X)=0$} are
special, and standard conjectures say they should have a virtually
abelian fundamental group. These conjectures are thus
special cases of conjecture I above.

Notice that Conjecture \ref{abconj} is established for manifolds
with $c_1(X)=0$ (which are indeed special,
and have almost abelian fundamental group, by the Calabi-Yau theorem).

This provides strong support for (and motivates) the case where $\kappa=0$.

Let us give another example of result supporting
Conjecture \ref{abconj} in the case $a(X)=0$:

\begin{theorem} ([C95]) Let $X\in \mathcal C$ be such that:
$a(X)=0$ and $\chi (X,\mathcal O$$_{X})\neq 0$. Then:
$\pi_{1}(X)$ is finite (in particular: conjecture \ref{abconj} holds for $X$).

More generally, if $F$ is the generic fibre of the algebraic
reduction map of any $X\in \mathcal C$, with $\chi (X,\mathcal O$$_{X})\neq
0$, the image in $\pi _{1}(X)$ of $\pi _{1}(F)$ is finite. (Recall
that $F$ is special, by \ref{fibalgred}).
\end{theorem}

{\bf (4) Surfaces and Threefolds:} By \ref{specsurf} and \ref{spec3f}, Conjecture \ref{abconj}
holds for K\"ahler surfaces and projective threefolds with
$\kappa\neq 2$ (and more precisely, for all special threefolds,
except maybe for the ones with $\kappa=2$, and having a model of 
the Iitaka-Moishezon fibration with orbifold base a log-Enriques or log-Fano 
normal surface. In this case, the conjecture is open because of an insufficient knowledge 
of these orbifolds).

This case shows clearly that the solution of Conjecture \ref{abconj} above should be extended 
to the orbifold situation.  Recall that the fundamental group of 
$(Y/\Delta)$ is defined as the quotient of the fundamental group of the complement of 
$\Delta$ in $Y$, modulo the normal subgroup generated by elements of the form 
$\lambdaÑ{j} ^{m_j}$, where $\lambda_{j}$ is a small loop winding once aroud $\Delta_j$, 
if: 

$\Delta:=\sum_{j\in j}(1-1/m_j).\Delta_j$.

\begin{question} Let $(Y,\Delta)$ be a klt-orbifold, with $Y$ normal.
 Assume that  $-(K_{Y}+\Delta)$ is ample or torsion. Is then the
fundamental group of
$(Y/\Delta)$ trivial in the first case, and almost abelian in the second case?

More generally is the
fundamental group of
$(Y/\Delta)$ almost abelian if $\kappa(Y/\Delta)=0$, or if $(Y/\Delta)$ is 
 $\kappa$-rationally generated (see section \ref{orbrc} for this notion)?
\end{question}

Notice that a positive answer for surfaces solves the remaining cases of
the abelianity conjecture in dimension three.

{\bf (5) Manifolds with $-K_X$ nef} are conjectured to have a
virtually abelian fundamental group. This conjecture is thus also a
consequence
of conjecture \ref{abconj} under a positive answer to the question asked above of 
whether such manifolds are special. The
main results concerning this conjecture are in [D-P-S], [D-P-S'], [Pa]
and [Zh].

{\bf (6) K\"ahler manifolds $X$ covered by $\mathbb {C}^d$} are
special (by the results of the next section \ref{ko}, 
and S.Iitaka conjectured them to be covered by complex tori. 
Easy arguments (given below) 
show that this happens precisely if $\pi_1(X)$ is virtually abelian. So
this conjecture is also a special case of the abelianity conjecture above.\\

We are not presently able to solve substantial examples of
the abelianity conjecture. But, at least, the solvable and linear quotients of the
fundamental group can be dealt
with.

\subsection{\bf Linear and Solvable Quotients}

Recall first (\ref{albsurj}):

\begin{theorem}\label{alb'} Let $X$ be any finite \' etale cover of a special manifold.

 Let $\alpha_X:X\ra
Alb(X)$ be the Albanese map of $X$.
Then it is surjective, connected. Moreover, $\alpha_X$  has no
multiple fibres in codimension one if $X$ is special.
\end{theorem}

\begin{remark} Due to \ref{k=0}, this extends and sharpens
slightly a result of [Kw].
\end{remark}

\begin{remark} Assume that conjecture \ref{abconj} holds (or alternatively,
that $X$ is special and has almost abelian fundamental group).
Easy standard arguments then show, using the surjectivity of the Albanese
map and its absence of multiple fibres in codimension one, that the universal
cover of $X$ is holomorphically convex
(in accordance with the conjecture of Shafarevitch),
and that its Remmert reduction makes it proper over an affine
complex space, which is the universal cover of its Albanese variety. 
\end{remark}

Now using [C4] (see also [A-N],[B],[Sim]), one gets:

\begin{corollary} Let $X\in \mathcal C$ be special, and let $c:X'\ra X$
to be an \'etale finite cover. Then:

(1) The derived group $D\pi_1(X)$ of $X$ is finitely generated.

(2) The Green-Lazarsfeld set $\Sigma^{1}(X')$ consists of
finitely many torsion characters of $\pi_1(X)$.
\end{corollary}

  Again using [C4], corollary 4.2' (based on Deligne's
semi-simplicity and invariant cycle theorems), and \ref{alb'}, we get:

\begin{theorem}\label{solv} If $X$ is w-special and if
$\mu:\pi_1(X)\ra G$ is a surjective morphism of groups, with $G$
solvable and torsionfree, then $G$ is virtually abelian. 
(We say that torsionfree solvable quotients of $\pi_1(X)$ are almost abelian).
\end{theorem}

{\bf Proof:} This is an easy consequence of \ref{alb'} (which says that
the Albanese map of $X$ is surjective and connected) and of corollary
4.2' of [C4] (which says that if the Albanese map of $X'$ is
surjective for any finite \'etale cover $X'$ of $X$, then the
conclusion of \ref{solv} is true for $X$) $\square$

Remark this applies in particular when $a(X)=0$ (it is shown by a
different method in [C4]), and when $\kappa(X)=0$).

Now, using a result of [Z], itself based on results of N.Mok [Mk] (see
also [Siu]), and \ref{solv}, we get:

\begin{theorem} \label{lin} Let $X$ be w-special
, and let $\rho:\pi_1(X)\ra Gl(N,\mathbb C)$ be a  linear
representation.
Then $G:=Image(\rho)$ is virtually abelian.
(We say that linear quotients of $\pi_1(X)$ are almost abelian).
\end{theorem}

  So, counter-examples to conjecture \ref{abconj}, if any, can't be detected
by linear representations of their fundamental groups.

{\bf Proof:} Let $G'$ be the Zariski closure of $G$. Replacing $X$ by
a suitable \'etale cover, we assume that $G'$ is connected.
We have an exact sequence of groups:

$$1\ra R\ra G'\ra S\ra 1 , $$

with $S$ semi-simple and $R$ solvable. Consider the representation
$\rho'=\sigma \circ \rho:\pi_{1}(X)\ra S$, where $\sigma:G\ra S$ is the
above quotient.
By [Z], thm.5,p.105, we get that $S$ is trivial (ie: reduced to the
unit element). Thus $G\subset R$ is solvable, and can be assumed
to be torsionfree, passing to a suitable finite \'etale cover, by
A.Selberg's theorem. The conclusion now follows from \ref{solv} $\square$

\begin{corollary} Conjecture \ref{abconj} holds when $\pi_{1}(X)$
is linear (ie: has a faithfull representation $\rho$.)
\end{corollary}

\subsection{\bf Manifolds Covered by $\mathbb C^n$}

An easy consequence of \ref{lin} is:

\begin{theorem}\label{liniitak} Let $X \in \mathcal C$ be covered by $\Bbb
C^d$, so $d:=dim(X)$. Assume that $\pi_1(X)$ is either solvable, or linear 
(ie: admits a faithful linear
representation in some $Gl(N,\Bbb C)$).
Then $X$ is covered by a complex torus. (ie: S.Iitaka's Conjecture
then holds for $X$).
\end{theorem}

{\bf Proof:} We know that $X$ is special by \ref{ko} below (one can also
get directly from [K-O] that $X$ is w-special, which is sufficient
for our purpose).
Thus $\pi_{1}(X)$ is virtually abelian by \ref{solv} or \ref{lin} above.

 Replacing
$X$ by a suitable unramified cover, the conclusion thus
follows from the next Lemma \ref{covtor} $\square$

\begin{lemma}\label{covtor} Let $X\in \mathcal C$ be covered by $\mathbb
C^d$, with $d:=dim(X)$. Assume that $\pi_{1}(X)$ is abelian. Then $X$
is a complex torus.
\end{lemma}

{\bf Proof:} Let $\alpha:X\ra A$ be the Albanese map of $X$. By
\ref{alb'}, it is surjective and connected. Let $F$ be a generic fibre of
$\alpha$.

It is sufficient to show that $dim(F)=0$, because $X$ is covered by
$\mathbb C^d$ (indeed, $\alpha$ is then birational. And so $X$ contains a
rational curve if $\alpha$ is not isomorphic. But this rational curve 
would lift after
normalisation to the universal cover of $X$. A contradiction).

Assume $dim(F)>0$; then $(image (\pi_{1}(F)\ra \pi_{1}(Y)))$ would be
infinite, again because $\mathbb C^d$ does not contain positive
dimensional
subvarieties. Because $\pi_{1}(X)$ is abelian, this implies by
dualising that the restriction map $H^{0}(X,\Omega^{1}_{X})\ra
H^{0}(F,\Omega^{1}_{F})$
is nonzero. A contradiction to the fact that $F$ is a fibre of $\alpha$ $\square$

We can slightly extend \ref{lin} to the relative case, which shows that interesting 
linear representations of $\pi_1(X)$, $X\in\mathcal C$, come from the core:

\begin{theorem}\label{lin"} Let $X\in \mathcal C$, let $c_{X}:X\ra C(X)$
be its core with $F$ its generic fibre,
   and let $\rho:\pi_{1}\ra Gl(N,\mathbb C)$ be a linear
representation. 

Let $G$ (resp.$H$) be the image of $\pi_{1}(X)$ by
$\rho$
(resp. the intersection with G of the Zariski closure of the image
$H$ of $\pi_{1}(F)$
by $\rho$). Then $H$ is normal in $G$, almost abelian, and there is a
natural quotient
linear representation $\rho': \pi_{1}(C'(X))\ra (G/H)\subset
Gl(N',\Bbb C)$, where
$C'(X)$ is a suitable finite \'etale cover of $C(X)$.
\end{theorem}

{\bf Proof:} This is an immediate consequence of
\ref{lin}, using the well-known fact ([C2], for example) that in the
situation of \ref{lin"}, if $K:=Ker(\pi_{1}(X)\ra \pi_{1}(C(X))$, then
$\pi_{1}(F)_{X}\subset K$ , and the quotient group is generated by
torsion elements. Thus, if
$K':=\rho(K)$, since $H$ is almost abelian by \ref{lin}, and the
quotient group $(G/H)$ is linear
by standard linear algebraic group theory, we have an embedding
$(G/H)\subset Gl(N',\Bbb C)$, for some $N'$, and
we see that  $(K'/H)$ is the trivial group by A.Selberg's theorem,
provided we replaced $X$ by some suitable finite \'etale cover.

 We thus get then an isomorphic morphism of groups:
 $\mu: (G/K')\ra (G/H)$, which gives the asserted representation $\rho'$ $\square$

\section{\bf An Orbifold generalisation of Kobayashi-Ochiai's Extension theorem}\label{orbko}

\subsection{\bf Statements}

We give here an orbifold version of a famous extension theorem,
 due to Kobayashi-Ochiai ([K-O]). 

The setting considered in all of this section \ref{orbko} is the following: 

\begin{notation}\label{notko}$V$ be a connected
complex manifold, $D$ a reduced divisor on $V$, and $U:=V-D$.
Let $\phi:U\mero X$ be a meromorphic map, with $X\in \mathcal C$.

We always assume  $\psi:=f\circ \phi:U\mero Y$ is nondegenerate (ie:
has maximal rank $p=dim(Y)$ somewhere)

Let $f:X\mero Y$ be a fibration, with $p:=dim(Y)$.
\end{notation}

Our main result here is:

\begin{theorem}\label{ko} Let $V,U,X,Y,\phi,f$ be as above. 

assume $f$ is a fibration of general type.Then:

{\bf(a)} $\psi$ extends meromorphically to $\psi':V\mero Y$.

{\bf(b)} For any $m>0$ sufficiently divisible and $s\in
H^0(Y,m(K_Y+\Delta(f)))$, $\psi^*(s)$ extends to:

$({\psi}){'}^*(s)\in  H^0(V,(\Omega_V^p)^{\otimes m}((m-1).D))$.
\end{theorem}

\begin{remark} The result of Kobayashi-Ochiai is the special case 
where $X=Y$, so that $X$ is of general type. 
\end{remark}

The proof will be given in the next sections. It is an 
orbifold modification of the original proof of [K-O].

We shall first give some applications of this result, which are 
criteria of transcendental nature for certain varieties to be special.

\begin{definition}\label{logspecdef} Let $U=V-D$ be the complement of a reduced divisor $D$ in
the connected complex manifold $V$. Say that the
pair $(V,D)$ is {\bf log-special} if: 
for $p>0$, there is no rank one coherent subsheaf
$L\subset (\Omega_V^p)$ and no $m>0$ such that
the complete linear system defined by $H^0(V,mL+(m-1)D)$
is of general type (ie: defines a meromorphic map of maximal rank
$p$).

\end{definition}

As a first application, we immediately get from \ref{ko}:

\begin{corollary}\label{logspec} Let $U=V-D$ be the complement of the divisor $D$ in 
the connected complex manifold $V$. Let $\phi:U\mero X$ be nondegenerate. 

Assume that the pair $(V,D)$ is log-special. Then: $X$ is special.
\end{corollary}

We now consider the case when $V=X$ and 
$D$ is empty. So the extension part of \ref{ko} is certainly 
trivial. Observe that the pair $(V,D):=(X,0)$ is log-special precisely if there is no Bogomolov 
sheaf on $X$. We however recover one direction of \ref{specnobog}:

\begin{corollary} Assume that $X$ has no Bogomolov sheaf. Then $X$ is special.
\end{corollary}

To get some concrete applications, we take more restrictive conditions than log-specialness.

\begin{definition}\label{logrc} Let $U=V-D$ be the complement of a reduced divisor $D$ in
the connected complex manifold $V$. Say that the
pair $(V,D)$ is {\bf log-RC} (for log-rationally connected) if:
 
$ H^0(V,(\Omega_V^p)^{\otimes m}((m-1).D))$ vanishes for any $m>0$.
\end{definition}

Of course, if $(V,D)$ is log-RC, it is log-special, if $V$ is projective rationally connected and 
$D$ is empty, the pair $(V,D)$ is log-RC. 

One easily checks the following (shown in a generalised form below):

\begin{example}\label{c^nlogrc}The pair
$(\mathbb P^n,\mathbb P^{n-1})$ is log-RC.
\end{example}

From \ref{ko}, we get:

\begin{corollary} \label{logrcspec} Assume that the pair $(V,D)$ is log-RC. Then: $X$ is special.
\end{corollary}

Specialising from Example \ref{c^nlogrc}, we get the:

\begin{definition} Say that $X$ is {\bf $\mathbb C^d$-dominable} (resp. 
{\bf covered by $\mathbb C^d$}) if there exists 
$\phi:\mathbb C^d\mero X$, nondegenerate (resp. if the universal cover 
of $X$ is isomorphic to $\mathbb C^d$).
\end{definition}

The terminology $\mathbb C^d$-dominable is from [B-L]

From \ref{logrcspec}, we then deduce:

\begin{corollary}\label{cddomspec} Let $X\in \mathcal C$ be {\bf $\mathbb
{C}^d$-dominable}. Then $X$ is special.

In particular: if $X$ is covered by $\mathbb C^d$, then $X$ is special.
\end{corollary}

  Remark that manifolds $X\in \mathcal C$ covered by $\mathbb {C}^d$
are conjectured by S.Iitaka to be covered
by a torus (ie: to have a finite
\'etale cover which is a torus). This is shown to be true for surfaces in
[I2], and for threefolds in [C-Z].
For $X$ projective of arbitrary dimension,
this follows (see [C-Z])
from the standard conjectures of the Minimal Model Program. See
\ref{liniitak} for another case where this conjecture is known to
hold true.

We now give a geometric criterion for the pair $(V,D)$ to be log-RC. 
This criterion may explain the name log-RC, and shows that being 
log-RC may be seen as a transcendental analogue of rational connectedness.

\begin{definition} \label{h'spec}
We say the pair $(V,D)$ is {\bf $H'$-special} if there exists a meromorphic
map $\phi:U\mero X$,where $U=W\times \mathbb {C}$,with $W$
quasi-projective, and moreover such
that $\phi$ is holomorphic and constant at the generic point of $W\times\{0\}$ 
(ie: $\phi(w,0)=a$, for some fixed $a\in X$ and any
$w\in W$).
\end{definition}

\begin{example} The pair $(\mathbb P^n,\mathbb P^{n-1})$ is $H'$-special. (Blow-up a point outside 
of $\mathbb P^{n-1}$ to see this).
\end{example}

\begin{proposition}\label{h'speclogrc} The pair $(V,D)$ is log-RC if it is $H'$-special.
\end{proposition}

The easy but lengthy proof is deferred to section \ref{h'specimpllogrc} below.

Although the following clearly results from 
the preceding observations, we state it separately:

\begin{corollary}\label{h'specspec} Let $X\in \mathcal C$ be {\bf
$H'$-special}. Then $X$ is special.
\end{corollary}

\begin{remark} Due to \ref{ko}, and in this
situation, it is tempting to ask if $X$ is special when $V$ itself is
special.
\end{remark}

This is not true, because transcendental holomorphic
maps do not preserve multiplicities along $D$. The following example
shows this.

\begin{example} Let $B,E$ be elliptic curves, and $a_i \in
B$ be distinct points, $i=1,...,4$. Let $B^*$ be the complement of 
these points in $B$.
   Let $V:=E\times B$ and $U:=E\times B^* \subset V$. Let now $m\geq
2$ be an integer. Let now $X$ be obtained from $V$ by
applying
to $V$ the four logarithmic
transformations of multiplicity $m$ respectively to its fibres
$E\times \lbrace a_i \rbrace$  over $B$.

\end{example}

Then $X$ is projective, and contains $U$ as a
Zariski open subset. However, $V$ is special, while $X$ is not.

So the specialness of $V$ does not imply that of $X$, even when $\phi$ is the
embedding of a Zariski open subset of both $V$ and $X$. (Except when 
$X$ is a curve, as
observed by M. Zaidenberg (this is an immediate consequence of [K-O])).

\subsection{\bf  Sketch of Proof of the Kobayashi-Ochiai's
Extension Theorem.}\label{pfko}

We shall very briefly first sketch the proof of the classical result of Kobayashi-Ochiai, and 
then indicate in section \ref{pforbko} the changes needed to obtain the orbifold version. 

We consider here only the extension part of the statement (ie: the fact that $\psi$ extends 
meromorphically to $V$). The second part (the control of the order of poles of $(\psi')^*(s)$) is 
shown in section \ref{vanord} by a lemma of independent interest.

Recall we want to show that if $\psi:U:=(V-D)\mero Y$ is nondegenerate, 
then $\psi$ extends meromorphically to $V$, provided $Y$ is of general type. 
The idea is thus that maps to general type varieties do not have essential singularities. 
This is analogous to the fact that a bounded function on the unit disc minus the origin extends across the origin. 
Now general type manifolds behave like manifolds with universal cover a bounded domain, in the sense 
that they admit pseudo-volume forms with negative Ricci form.

The proof (as explained in, [Kob],Chap.7) goes in several steps:

{\bf (1)} It is sufficient to deal with the special case  $V=\mathbb D^p$, 
$D=\{ 0\}\times \mathbb D^{p-1}$, so that:

 $U=\mathbb D^p-\{0\}\times \mathbb D^{p-1}:=\mathbb D^{(p)}$. 

This is simply because meromorphic maps to a projective variety (or even a K\" ahler manifold) 
extend meromorphically 
across codimension $2$ or more analytic subsets. 

The case where $dim(V)>p$ can be reduced to the equidimensional case (see [Kob], p.374-75).

Notations: For future reference, we denote by  $\mathbb D (r)$ the disk of radius $r$ centered at $0$ in 
$\mathbb C$, by $\mathbb D:=\mathbb D (1)$ the unit disk; 
$\mathbb D^p(r)$ and $\mathbb D^{(p)}(r)$ are constructed as $\mathbb D^{p}$ and 
$\mathbb D^{(p)}$above, replacing $\mathbb D$ by $\mathbb D (r)$.

{\bf (2)} Because $Y$ is of general type, it has a pseudo-volume form $w$ 
of Ricci curvature negative, bouded from above by a constant $-1/C<0$. 

This form is constructed out of a basis of a linear subsystem $L$ of $m.K_Y$ giving 
an embedding of $Y$ into some projective space (once its fixed components have been removed). 

Another property used in the proof is that if $w_j$ is the pseudo-volume form associated to 
a single element of the basis of $L$ above, the quotient $w_j/w$ is a smooth function on $Y$. 
This function is thus bounded. 

These notions and constructions are explained in more detail in section \ref{gentyp} below.

{\bf (3)} An elementary lemma (see \ref{extg} below), together withe argument exposed in (5) below, 
show that the special case (1) above is true, if one shows that $w$ being
 the preceding  pseudo-volume form on $Y$, then 
the integral of $\psi^*(w)$ converges on $\mathbb D^{(p)}(1/2)$.

{\bf (4)} By the celebrated Ahlfors-Schwarz lemma (see  \ref{a-s} below), we have the estimate: 
$w\leq C.\beta$, if $\beta$ is the homogeneous volume form on $\mathbb D^{(p)}$, 
with constant Ricci curvature $-1$, and $C$ is the positive constant of (2) above.

 One concludes from the elementary fact that $\beta$ has finite volume on 
$\mathbb D^{(p)}(1/2)$. 

{\bf (5)} Let us now conclude the proof in the equidimensional case: the forms $w_j:=\psi^*(s_j)$ of the linear system 
$\psi^*(m.K_Y)$  extend meromorphically to $\mathbb D^{p}$, by the steps (2,3,4,5). 
But the map $\psi$ is defined by the $N$-tuple 
of maps: $(w_0:w_1:...:w_N)$ with values in $\mathbb P^N$. The map $\psi$ itself thus 
extends meromorphically to $\mathbb D^{p}$, 
as claimed.

We shall refer to [Kob] for the details not given here. For the main points (2,3,4), 
we shall introduce the basic definitions and properties in the sections \ref{gentyp} 
and \ref{psvol} below. Finally, the proof of the orbifold version will be given 
in section \ref{pforbko}. 

The proof of point (b) in the statement of Theorem \ref{ko} 
will be given in the separate section \ref{vanord} below. (Although it might possibly be 
deduced from [Kob], which contains a similar result in implicit form, we prefer 
to give an independent, more general statement).

\subsection{\bf Pseudo-Volume and Ricci Forms}\label{psvol}

We need to recall some classical notions. See, for example, [Kob],
Chapters 1 and 7 for more details .

Let $M$ be a complex manifold. We say that $v$ is a {\bf
pseudo-volume form }(with holomorphic
degeneracies) on $M$
if it is a form of type $(p,p)$
on $M$, with $p:=dim(M)$, such that locally, in any coordinate system
$(x_{1},....,x_{p})=(x)$, one has:
$v$=$\mid a\mid^{2}.T.vol_{(x)}$, with:
$vol_{(x)}$:=$(i^{p^2}.d(x)\wedge\overline {d(x)})$, $a$ is a holomorphic nonzero function,
and $T$ a smooth  (ie: $\mathcal C$$^{\infty}$) everywhere positive function locally
defined on $M$.

For such a form $v$, its {\bf Ricci form} is the real $(1,1)$-form: $Ricci(v)$:=
$-dd^{c}(log T)=-i\partial\overline{\partial} (log T)$, and its {\bf Ricci function} 
is: $K_{v}$:=$(-Ricci(v))^{m}/v$.

Both are well-defined independently of local coordinates,
because changing charts multiplies $T$ by the square of the modulus of the Jacobian,
the $log$ of which has vanishing $dd^c$.

  From the definitions, one gets also immediately the functoriality of these notions: 
if $f:N\ra M$ is
a dominant holomorphic map, with $p=dim(N)$,
and if $v':=f^{*}(v)$, then the definitions
apply to $v'$, and one has both: $Ricci(v')=f^{*}(Ricci(v))$ and $K_{v'}=
f^{*}(K_{v})$.

\begin{remark}\label{remric} The Ricci form is thus the special case for $-K_M$ of the 
notion of curvature form of a singular hermitian metric on a line bundle, neglecting the singular part 
(cohomologically significant, however) of the curvature current.
\end{remark}

Assume that $v$ has holomorphic degeneracies. Then we say that $v$
has {\bf negatively bounded Ricci curvature } if $(-Ricci(v))$ is an
everywhere positive or zero $(1,1)$-form,
and if $-K_{v}\geq C$ for some constant $C>0$, everywhere on $M$.

Observe that $K_{v}$ is then nonpositive from the condition on
$Ricci(v)$, and that $K_{v}$ may take the value
$(-\infty)$ at points where $a$ vanishes. Note that if $M$ is
compact and $-Ricci(v)$ everywhere positive, then $v$ has negatively
bounded Ricci curvature. 

Observe too that $v'$ has negatively bounded
Ricci curvature if and only if so has $v$, for $v'=f^{*}(v)$ as above.

\begin{example}\label{psvolex} We give the examples used in the proof. As said above, these are bounded domains.
\end{example}
(1) The unit disc. Take (in the linear coordinates) $v:=idx\wedge \overline {dx}/(1-\mid x\mid^2)^2$.
 
Then $Ricci(v)=-v$, and $K_v=-1$.

(2) The unit disc with origin removed $\mathbb D^*$. Take (in the linear coordinates):
 
$v:=idx\wedge \overline {dx}/(\mid x\mid)^2.(log(\mid x\mid)^2)^2$. 
One again has: $Ricci(v)=-v$, and $K_v=-1$. Actually the forms in the preceding two examples correspond under the 
universal covering map from $\mathbb D$ to $\mathbb D^*$. This can be checked explicitely using the Poincar\' e upper 
half plane as an intermediate step.

An important property used in the proof is that $\int _{\mathbb D^*(r)} v<+\infty$ if $r<1$. 
(Of course, the domain over which the integral is taken is the disc of radius $r$ with origin deleted);

(3) $M=\mathbb D^{(p)}$. Take (in the linear coordinates):

$v:=i^{p^2}dx\wedge \overline {dx}/(\mid x_1\mid)^2.(log(\mid x_1\mid)^2)^2.\prod _{2\leq j\leq p} (1-\mid x_j\mid^2)^2$.

One again has: $Ricci(v)=-v$, and $K_v=-1$. 

Again, one has: $\int _{\mathbb D^{(p)}(r)} v<+\infty$.

The above volume forms have thus negatively bounded Ricci function. 

We give in the next section the example used in the proof in an essential way.

\subsection{\bf Pseudo-Volume forms on Varieties of General Type}\label{gentyp}

  Let $Y$ be a $p$-dimensional manifold of general type, that is: with $K_Y$ big. 

From Kodaira's lemma, for $m$ large enough, we can write: $mK_Y=H+A$, with $H$ effective and 
$A$ very ample.

Let $L$ be a free linear subsystem of the complete linear system $\mid A\mid$, such that the 
associated regular map $\Lambda: Y\ra \mathbb P^N$ ($N=dim L$) is an embedding.

Let $(s):=(s_0,s_1,...,s_N)$ be a (complex) basis of $H+L$, and $h$ a section of $\cal O$$_Y(H)$ vanishing 
exactly on $H$.

Define $w:=(\sum_{j\in \{0,...,N\}}i^{m.p^2} s_j\wedge \overline{s_j})^{1/m}$. 

Then $w$ is a pseudo-volume form with holomorphic degeneracies along $H$, since it is written in 
local coordinates as: $w=\mid a\mid^2. v_1$, with 
$w_1=(\sum_{j\in \{0,...,N\}}\mid t_j\mid ^2)^{1/m}.vol_{(x)}$, 
if $a$ is a local equation for $H$, and $s_j=t_j.(vol_{(x)})^{\otimes m}, \forall j$ (using the notation 
introduced in section \ref{pfko} above).

Notice that $L$ being free, $w_1$ is a smooth volume form on $Y$. We now compute its Ricci form.

By its very definition, we have: $Ricci(w)=-\Lambda^*(\Theta)$, where $\Theta$ is (up to a positive constant) 
the curvature form of the Fubini-Study metric on the tautological line bundle $\cal O$$_{\mathbb P^N}(1)$. 

(This is in accordance with remark \ref{remric}: the Ricci form od $w$ is just $-\Lambda^*(\Theta)$, 
since the singular part of the metric is given by $h$).

The Ricci form of $w$ is thus negative everywhere. As a consequence: $K_w:= (Ricci(w)^p/v)$ is negative everywhere 
(and $-\infty$ on $H$). By the compactness of $Y$, we conclude that this function is negatively bounded by some constant 
$-1/C$, for some $C>0$.

\begin{remark}\label{wbound} For $j\in \{0,...,N\}$, define: $w_j:=(i^{m.p^2} s_j\wedge \overline{s_j})^{1/m}$. This is 
still a pseudo-volume form with holomorphic degeneracies on $Y$.  The function $w_j/w$ is obviously smooth on $Y$, and 
thus bounded. There exists $B>0$ such that $w_j\leq B.w, \forall j\in \{0,...,N\}$.
\end{remark}

\subsection{\bf The Lemma of Ahlfors-Schwarz}\label{a-slem}

Despite its simplicity, it is a very powerful tool to obtain bounds on the growth of pseudo-volume forms. 
We give the simplest version, sufficient for our applications.

\begin{lemma}\label{a-s} Let $v$ be a pseudo-volume form with holomorphic degeneracies on $\mathbb D^{(p)}$. 

Assume that $K_v\leq -1/C<0$ everywhere on $\mathbb D^{(p)}$. Then $v\leq \beta$, where $\beta$ is the 
volume form on $\mathbb D^{(p)}$ with constant Ricci function $-1$ defined in example \ref{psvolex} (3) above.
\end{lemma}

{\bf Sketch of proof:} See [Kob], (2.4.14) for details. One can reduce (by introducing the volume forms relative 
to polyradii $r<1$) to the case when the quotient (smooth) function $k:=v/\beta$ has a maximum at an interior 
point $b$ of (the closure of $\mathbb D^{(p)}$). Obviously, $v$ does not vanish at $b$

The real $(1,1)$-form $dd^c(log k)=i\partial\overline{\partial}(log k)$ is then nonnegative at $b$. (By functoriality, 
it is sufficient to check this on a complex parametrised curve. But then:
 $i\partial\overline{\partial}(log k)=(\partial ^2 (log k)/\partial s^2+\partial ^2 (log k)/\partial t^2).ds\wedge dt$ 
in real coordinates 
$x=s+it$, and the claim follows).

 Rescaling $v$, we assume that $C=1$. 
Now: $i\partial\overline{\partial}(log k)=i\partial\overline{\partial}(log V)-i\partial\overline{\partial}(log B)=
Ricci(\beta)-Ricci(v)\leq 0$, if $V,B$ are the functions such that $v=V.vol_{(x)}$ and $\beta=B.vol_{(x)}$.

We get: $Ricci(\beta)\leq Ricci(v)$. Hence also: $-v\geq K_v.v=Ricci(\beta)^{\wedge p}\leq Ricci(v)^{\wedge p}=
K_{\beta}.\beta=-\beta$, and 
the conclusion follows $\square$

We now state an extension criterion which is used crucially in the proof.

\begin{proposition}\label{extg} Let $s\in H^0(\mathbb D^{(p)}, m.K_{\mathbb D^{(p)}})$ be such that:

$\int_{\mathbb D^{(p)}(r)}(i^{m.p^2}s\wedge\overline{s})^{1/m}<+\infty$, for some $r<1$. 

Then $s$ extends to a meromorphic section 
of $mK_{\mathbb D^{p}}$. Moreover, the poles of this extension are of order at most $(m-1)$ along $\{0\}\times \mathbb D^{p-1}$.
\end{proposition}

The proof is given in [Kob], (7.5.7-8). The proof reduces to the case when $p=1$ and the elementary fact that a holomorphic 
function on the punctured unit disc which is L$_{(2/m)}$ on some $\mathbb D(r), r<1$ extends meromorphically with a pole of 
order at most $(m-1)$.

\subsection{\bf Proof of the Orbifold Version}\label{pforbko}.\\

{\bf Setting:} We consider in this section the following
data: $f:X\mero Y$ is a fibration of general type,
$\phi:U\mero X$ is a meromorphic map which is dominant (ie: submersive
at some point), $X,V$ are connected complex manifolds, with $X$ 
compact, and $j:U\subset V$ is a Zariski dense
open subset, complement of some divisor $D$ of $V$.

Let $\psi:=f\circ \phi:U\mero Y$. Blowing-up $X$ and $Y$ if needed,  we
shall assume that $f$ is holomorphic, that $Y$ is projective, and that $f$ is high (see \ref{high}).

We let $p>0$ be the dimension of $Y$, and
$\Delta:=\Delta(f)$ be the multiplicity divisor for $f$, as defined in \ref{multdiv}

Our purpose here is to establish that $\psi$ extends meromorphically to $V$.

\begin{proposition}\label{orbliftmK} Let $m>0$ be an integer such that
$m.\Delta$ is Cartier. Then:

(1) $f$ defines an injection of sheaves $f^{*}:\mathcal
O$$_{Y}(m.(K_{Y}+\Delta))\ra (\Omega ^{p}_{X})^{\otimes m}$,
   after suitable modifications of $X$ and $Y$.

(2) In particular, let $s\in H^{0}(Y,m.(K_{Y}+\Delta))$. Then
$f^{*}(s)\in H^{0}(X,(\Omega ^{p}_{X})^{\otimes m})$.
(Recall $Y$ has dimension $p$).
\end{proposition}

{\bf Proof:} Just as in the proof of \ref{injsheav}, because we have an
injection of sheaves (also denoted $f^{*}$):
$f^{*}:\mathcal O$$_{Y}(m.(K_{Y}))\ra (\Omega ^{p}_{X})^{\otimes m}$, it
is sufficient to check the above injection outside a
codimension two subvariety of $X$, provided $f$ is ``high".

So we first assume that we are near (in the analytic topology) the
generic point of $\Delta_i$.

But the assertion follows here from an easy local computation: indeed
$m_i$ divides $m=(m'.m_i)$ by assumption,
   and we can choose local coordinates:

 $(x)=(x_{1},...,x_{d})$ near the generic point of the
divisor $D_i\subset X$, a component of $f^{-1}(\Delta_{i})$
mapped surjectively to $\Delta_i$ by $f$
in such a way that in suitable local coordinates
$(y):=(y_{1},...,y_{p})$, near the corresponding generic point of
$\Delta_i$, one has:

 $f(x)=(y)=(y_{1},...,y_{p})$, with
$y_{1}=(x_{1})^{m_{i}+m"}$, and
$y_{j}=x_{j}$ for $j>1$, with $m"$ some appropriate nonnegative integer.

(Local equations for $D_i$ and $\Delta_i$ being respectively
$x_1$=$0$ and $y_1$=$0$).

Thus we get a local generator of $f^{*}(\mathcal
O$$_{Y}(m.(K_{Y}+\Delta)))$ in the form:
$f^{*}(d(y)/y_{1}^{(1-1/m_{i})})^{\otimes m}$, and this expression is
equal to:

$f^{*}(d(y)^{\otimes m}/y_{1}^{(m-m')})$
$=x_1^{(m_i+m"-1) \dot m-(m_{i}+m").(m-m'))}.d(x')^{\otimes m}$
$=(x_1)^{m'.m"}.d(x')^{\otimes m},$

up to a nonzero constant factor,
where $d(y):=dy_{1}\wedge ...\wedge dy_{p}$ and $d(x'):=dx_{1}\wedge
...\wedge dx_{p}$.

Hence the claim, since the exponent of $x_1$ in the last expression is obviously nonnegative.

We thus see that $f^*$ is well-defined outside the finite union $B$
of all divisors of $X$
which are mapped to codimension two or more analytic subsets of $Y$.
By a suitable composition of blow-ups $u:X'\ra X$ and $v:Y'\ra Y$, we
get $f':X'\ra Y'$  holomorphic such that $f.u=v.f'$,
and moreover such that the strict transform $B'$ of $B$
in $X'$ has all of its irreducible components mapped onto a divisor
of $Y'$ by $f'$ (depending on the component, of course).

Then $g:=v^{-1}\circ f=f'\circ u^{-1}:X\mero Y'$ is holomorphic outside
a codimension
two or more analytic subset $A$ of $X$ (the indeterminacy locus of $u^{-1}$).

But then $g^*(\mathcal O$$_{Y'}(m.(K_{Y'}+\Delta(f')))$ injects into
$(\Omega_X^p)^{\otimes m}$ over $(X-A)$.
   By Hartog's theorem (for example), this injection extends through
$A$. Now composing with
$u^*:(\Omega_{X}^{p})^{\otimes m}\ra(\Omega_{X'}^{p})^{\otimes m}$,
we see that $(f')^*$ injects $(\mathcal O$$_{Y'}(m.(K_{Y'}+\Delta(f')))$
into $(\Omega_{X'}^{p})^{\otimes m}$, as asserted.

The second assertion is an immediate consequence of the first $\square$

The first part of this argument (in simpler form) applies to give:

\begin{proposition}\label{liftdnice}
Let $u::Y'\ra Y$ be a 
${\Delta}$-nice
covering, in the sense of \ref{dnice}. But we choose $u$ such that it ramifies at order exactly $m_j$ above 
each component $\Delta_j$ of $\Delta(f)$. Then:

(1) $u$ defines an injection of sheaves: $u^{*}:\mathcal O$$_{Y}(m.(K_{Y}+\Delta))\ra \mathcal O$$_{Y'}(m.K_{Y'})$.

(2) In particular, let $s\in H^{0}(Y,m.(K_{Y}+\Delta))$. Then
$u^{*}(s)\in H^{0}(Y',m.K_{Y'})$.
\end{proposition}

We now come to the crucial point of our modification of the proof of
the extension theorem of Kobayashi-Ochiai. In a first step, we proceed as in section \ref{gentyp}, by the construction of 
a pseudo-volume form on $Y$ with meromorphic-not holomorphic- degeneracies.

\begin{notation}\label{orbgentyp} 
\end{notation}

Let $\phi:\mathbb D^{(p)}\mero X$ be such that $\psi:=f\circ \phi:\mathbb D^{(p)}\mero Y$ is 
nondegenerate. We need to extend $\psi$ meromorphically to $\mathbb D^{p}$. 

Since $(K_{Y}+\Delta)$ is big on $Y$, we can write, by Kodaira's 
lemma: $m(K+\Delta)$=$H+A$,
with $A$ very ample and $H$ effective on $Y$, $m$ being choosen sufficiently
   divisible (by the l.c.m of the $m_i$'s) and sufficiently big.

Let $h\in H^{0}(Y,\mathcal O$$_{Y}(H))$ be an element defining $H$.

Let $\delta\in H^0(Y,\cal O$$_Y(m.\Delta))$ be a section vanishing exactly on $m\Delta$, with the right multiplicities.

For $\sigma\in \mid A\mid$, the complete linear system defined by $H$,
we denote by:

 $s:=h.\sigma.\in H^{0}(Y,m.(K_{Y}+\Delta))$.

Let now $\sigma_{0},...,\sigma_{N}$ a basis of $H^{0}(Y,A)$; let also
$s_{j}:=h.\sigma_{j}$, for $j\in J:=\lbrace 0,...,N\rbrace$.

Let further $t_j:=s_j/\delta, \forall j$. These are meromorphic sections of $m.K_Y$.

Define $v:= (\sum_{j=0}^{j=N} (i^{p^2.m} t_j\wedge\bar{t_j}))^{1/m}$. Thus  $v$ is a global pseudo-volume form with meromorphic degeneracies 
on $Y$. 

Notice however that, by Proposition \ref{orbliftmK} above, $f^*(t_j)$ is, for any $j$, a holomorphic 
section of $(\Omega_X^p)^{\otimes m}$. This is simply because in local coordinates on $Y$, we can write: 

$t_j=(T_j/d). dy$, for some holomorphic function $T_j$, where $d$ (resp. $dy$) is a local generator of 
the ideal defining $m\Delta$ (resp. $K_Y$). (I thank M.Paun for this important observation). The assertion 
thus follows from \ref{orbliftmK} and the fact that $dy/d$ is a local section of $m.(K_Y+\Delta)$.

In particular, $\psi^*(t_j)$ is a holomorphic section of $mK_M$, with $M=\mathbb D^{(p)}$, and so 
$w:=\psi^*(v)$ is naturally a pseudo-volume form with holomorphic degeneracies on $M$.

For the same reason, $f^*(s_j)$ is, for any $j$, a holomorphic 
section of $(\Omega_X^p)^{\otimes m}$. 

Exactly in the same way we defined $w:=\psi^*(v)$, we can define a pseudo-volume form with holomorphic degeneracies $\bar{w}$ on $M$ by:
$\bar{w}:= (\sum_{j=0}^{j=N} (i^{p^2.m} \psi^*(s_j)\wedge\overline{\psi^*(s_j}))^{1/m}$. 

Notice that here, however, there is no pseudo-volume form $\bar{v}$ on $Y$ such that $\bar{w}=\psi^*(\bar{v})$.

The crucial point is now:

\begin{lemma}\label{ricbound} $Ricci(\bar{w}=Ricci(w)=\psi^*(-\Lambda^*(\Theta))$, 
if $\Lambda:Y\ra \mathbb P^N$ is the regular map defined by the complete linear system 
$\mid A\mid$ on $Y$, and $\Theta$ the curvature form of the Fubini-Study metric on $\mathbb P^N$.
 In particular, $Ricci(w)$ is thus negative everywhere, and there exists 
$C>0$ such that $K_{\bar{w}}\leq -1/C$ everywhere on $\mathbb D^{(p)}$.
\end{lemma}

{\bf Proof:} Write again $t_j=(T_j/d) .dy$, as above. We deduce
that: $v=\mid h/d\mid ^{2/m}.V.vol_{(y)}$, if 
$h$ is a local equation for $H$, and $vol_{(y)}$ a local volume form on $Y$.
 Here $V=(\sum_{j\in }\mid T_j\mid ^2)^{1/m}$.

Now, from \ref{orbliftmK} again, we see that $\psi^*(vol_{(y)}/\mid d\mid ^{2/m})=\mid G\mid ^2.vol_{(M)}$, 
for a certain holomorphic function $G$ defined on 
$\psi^{-1}(U)$, for $U\subset Y$, the open suset where the said trivialisations and charts are defined. 

Recall that $-i\partial \overline{\partial}V=-\Lambda^*(\Theta)$, by \ref{gentyp} above.

A computation of the Ricci curvature then gives the first assertion for $w$. 

The assertion for $\bar{w}$ is deduced from the fact that $\bar{w}=\psi^*(\mid d\mid^{2/m}). w$, so that 
these two pseudo-volume forms have the same Ricci-form. 

We then get: 

$K_w=\psi^*(-\Lambda ^*(\Theta))^{\wedge p}/\psi^*(v)=
\psi^*(-\Lambda ^*(\Theta))^{\wedge p}/\psi^*(\mid h.G^m\mid^{2/m}.V).vol_{M}$

And so, $K_{\bar{w}}=\psi^*(\mid d\mid^{-2/m}.K_w\leq (-1/C).(1/\psi^*(\mid h\mid^{2/m}.V))$, 

because $V$ is differentiable, and $h$ locally bounded from above everywhere on $Y$.

To express things differently, one can compute the Ricci function of $\bar{w}$ symbolically, as if it were of the form 
$\psi^*(\bar{v})$, with $\bar{v}:= (\sum_{j=0}^{j=N} (i^{p^2.m} s_j\wedge\overline{s_j}))^{1/m}$.
$\square$

By the argument given in section \ref{pfko}, to show \ref{ko}, it is sufficient to show that each of the $\psi^*(s_j)$'s extend 
meromorphically to $\mathbb D^{p}$, and this is true if the integral of $w$ over $\mathbb D^{(p)}(r)$ is convergent, 
for some $0<r<1$. 

To show this convergence, we simply need, by the Ahlfors-Schwarz Lemma, to show that $w$ has everywhere negatively bouded 
Ricci function (by some constant $-1/C<0$). But this exactly what the preceding Lemma \ref{ricbound} claims.
This concludes the proof of the first assertion of \ref{ko}.

\begin{remark} Observe that the preceding proof shows with minor adaptations (using \ref{liftdnice}), the following orbifold 
version of Theorem \ref{ko}: If $(Y/\Delta)$ is an orbifold of general type ($Y$ smooth and $\Delta$ suported on a normal 
crossings divisor), then any nondegenerate meromorphic map $\phi:U\mero (Y/\Delta)$ extends meromorphically to $V$. 
Here, to say that $\phi$ is a meromorphic map to the orbifold $(Y/\Delta)$ means that any point $u\in U$ at which $\phi$ is holomorphic, 
with $y:=\phi(u)$ a smooth point of the support of $\Delta$, lying on $\Delta_i$, then $\phi$ lifts 
locally around $u$ to a holomorphic map to the ramified cover of $Y$ near $y$, ramifying at order exactly $m_i$ above $\Delta_i$.
\end{remark}

\subsection{\bf Extension of Pluri-Canonical Meromorphic Forms}\label{vanord}.\\

We now want to control the order of the poles along $D$ of the
meromorphic extension to
$V$ of $\psi^{*}(s)\in H^{0}(U,(\Omega_{U}^{p})^{\otimes m})$, for
$s\in H^{0}(Y,m.(K_{Y}+\Delta))$. More precisely, notations being as
above:

\begin{proposition}\label{ordvan} For any
   $m>0$
sufficiently divisible and
$s$ any element of $ H^0(Y,m(K_Y+\Delta(f)))$,
$\psi^{*}(s)$ extends to
$(\psi')^{*}(s)\in  H^{0}(V,(\Omega_V^p)^{\otimes m}((m-1).D))$.
\end{proposition}

{\bf Proof:} Recall that the support of $\Delta$ is of normal
crossings. We can assume that $\psi'$ is holomorphic, because the
assertion holds for $V$
if it does for any of its blow-ups. 

Recall
that $s$ can be written in local coordinates as follows, s being a section of
$m(K_Y+\Delta(f))$:

$s=(d(y)/(y_{1}^{(1-1/m_{1})}....y_{p}^{(1-1/m_{p})}))^{\otimes m}$,
   if $\Delta$ has local equation the denominator of this expression.
(Some of the
$m_{i}$ may be one, and so the corresponding terms don't contribute).

Let now $(v)=(v_{1},...,v_{n})$ be local coordinates near a generic
point of some component $D'$ of $D$, such that $D'$ has local
equation $w:=v_{1}=0$.
We can assume that near this point the map $\psi'$ is given by
$\psi'(v)= (v_{1}^{r_{1}}.Y_{1},...,v_{1}^{r_{p}}.Y_{p})$, for
integers $r_{i}$ and nowhere
vanishing functions $Y_{i}$ of $(v)$.

Computing, we get (letting $r$ be: $r_{1}+...+r_{p}$): ${(\psi)}{'}^{*}(s)=
(w^{r-1}.\Theta/w^{r-s}.Y')^{\otimes m}$, with $w:=v_1$, for some
holomorphic $p$-form $\Theta$ on $V$,
and $Y'$ a nowhere vanishing function on $V$ (the product of the
$Y_{i}$'s), and $s:=r_{1}/m_{1}+...+r_{p}/m_{p}$. This shows the
proposition, because $s.m$ is
a positive integer, since each $m_{i}$ divides $m$ $\square$

\subsection{\bf Proof that $H'$-Specialness Implies Log-Rational Connectedness}
\label{h'specimpllogrc}.\\

We prove here proposition \ref{h'speclogrc}. 

We use the notations and setting of \ref{notko}, \ref{h'spec} and \ref{logrc}.

{\bf Proof:} We just need to show that any element $s'$ of
$H^{0}(V,(\Omega_V^p)^{\otimes m}((m-1).D))$ vanishes, provided it is
on $U$ of the form
$\phi^{*}(s")$ for some
$s"\in H^{0}(X,(\Omega_{X}^{p})^{\otimes m})$. For this, we shall
just show that its restriction to any $F$:=$\lbrace w\rbrace\times
\mathbb P^{1}(\mathbb C)\subset V$
has to vanish for any $w\in W$. The argument is similar to (and
motivated by) the one showing that sections of $(\Omega_X^p)^{\otimes
m}$ vanish
if $X$ is rationally connected (see [C5]).

Let $\Omega^q$ denote the restriction of $\Omega^q_V$ to $F$, for
$q\geq 0$, and let $N$ be the (trivial) dual of the normal bundle of
$F$ in $V$.

  From the exact sequence: $1\ra N\ra \Omega^{1}\ra\Omega^{1}_{F}\ra1$,
we get, for $q>0$:

$$ (@) 1\ra \wedge ^{q} N\ra \Omega^{q}\ra\Omega^{1}_{F}\otimes
\wedge ^{q-1}N\ra1$$ Notice that all these sequences are split.

We now fix $q=p>0$, and denote respectively  by $B:=\wedge ^{p} N$
the kernel, and $A:=\Omega^{1}_{F}\otimes \wedge ^{p-1}N$
the quotient of the corresponding exact sequence $(@)$.

  From this same split exact sequence, we get now for
$(\Omega^p)^{\otimes m}$ a decreasing (split) filtration by
subbundles:

$$\lbrace 0\rbrace:=W_{m+1}\subset W_1\subset ....\subset
W_0:=(\Omega^p)^{\otimes m},$$

with successive quotients: $W_0/W_1=A^{\otimes
m}$;.....;$W_k/W_{k+1}=A^{\otimes (m-k)}\otimes
B^k$;.....;$W_m=B^{\otimes m}$.

Notice that, for $k=0,...,m$, the bundle $A^{\otimes (m-k)}\otimes
B^k$ is isomorphic to $\mathcal O$$_{F}(-2(m-k))$ tensorised by
a trivial vector bundle $N_k$ on $F$.

We now take local coordinates $((w),t)$ on $V$ near $(w,0)$, where
$(w):=(w_1,.....,w_{n-1})$ are local coordinates on $W$, and $t$
is a a local coordinate near
$0\in \mathbb P^1$.

The map $\phi$ takes the form
$\phi((w),t)=(t.\phi_1,...,t.\phi_{d})$, for holomorphic functions
$\phi_i$ of $((w),t)$, with $d:=dim(X)$,
   in local coordinates near  $a\in X$, since
$\phi(W\times\lbrace 0\rbrace)=\lbrace a\rbrace$, by our
$H'$-speciality assumption.

Locally, if $\sigma"$ is a $p$-form on $X$ near $a$, then
$\sigma':=\phi^*(\sigma")$ is written:
$\sigma'=dt.\alpha+t.\beta$, for $\alpha$  and
$\beta$  forms in the $dw_{i}$'s. If now $\sigma"$ is a $p$-form,
then $\sigma'=t^{p-1}.dt\wedge\alpha+t^{p}.\beta$, with the same
conditions.

We now consider the case when $s'=\phi^{*}(s")$ for some $s"\in
H^{0}(X,(\Omega_{X}^{p})^{\otimes m})$: from the preceding remarks,
we can write
$s"=s_0"+...s_k"+...+s_m"$ , where each of the $s_k"$ is the tensor
product of $(m-k)$ terms of the form: $t^{p-1}.dt\wedge \alpha$, and of
$k$ terms of
the form: $t^{p}.\beta$.

So that, for each $k$, $s_k"$ is actually
the piece of degree $k$ of $s"$ for the above graduation, and so turns
out to be a section of the bundle $A^{m-k}\otimes B^{k}$ defined
above, and vanishing at $0$ to order
$(m-k).(p-1)+k.p$=$m.(p-1)+k$.

Now, we assumed that $s"$ has a pole of order at most $(m-1)$ at the
point at infinity of $F$. Thus $s_{k}"$ defines a meromorphic
section of
   $\mathcal O$$_{F}(-2(m-k))\otimes N_{k}$ , with $N_k$ a trivial bundle,
this section vanishing at $0$ to order $(m(p-1)+k)$ at least, and having
a pole of order $(m-1)$ at most at infinity.

This section thus has to vanish, because otherwise the number of poles should
be the opposite of the degree. But here: $(m-1)-m.(p-1)-k<2(m-k)$
for any integers $m,p>0$ and $0\leq k\leq m$.

This concludes the proof of \ref{h'speclogrc} $\square$

\section{\bf The Core and the Kobayashi Metric}

\subsection{\bf The Kobayashi Pseudo-Metric}\label{komet}

 We refer to [Kob] for a systematic
treatment of the following notions.

Let $X$ be a complex manifold. Let $d_X$ be the {\bf Kobayashi
pseudometric} on $X$ (which is the largest
   pseudometric $d$ on $X$ such that $h^*(d)\leq p_{\mathbb D}$ for any
holomorphic map $h:\mathbb D\ra X$, where $\Bbb D$ is the unit disc in
$\mathbb {C}$ and $p_{\Bbb D}$ is the Poincar\'e metric on $\Bbb D$).

For example, it is easy to check (using the Ahlfors-Schwarz lemma)
that $d_{\mathbb C^n}\equiv 0$, and $p_{\Bbb D}=d_{\Bbb D}$.

A fundamental property of $d$ is that holomorphic maps are distance
decreasing: if $h:X\ra Y$ is any holomorphic map, then
$h^*(d_Y)\leq d_X$. Hence if $h:\mathbb C\ra X$ is holomorphic, then
$d_X$ vanishes on the (metric) closure of the image of $h$.
Rational and elliptic curves have thus vanishing $d$.

In particular, if $h:\Bbb C\ra X$ is holomorphic, then $d_X$ vanishes on the
(metric) closure $H$ of the image of $h$. Thus $d_X$ will vanish 
identically if, for example,
$H=X$, or if any two points of $X$ can be connected by such a chain of $H_i$'s.

On the opposite side, $X$  is said to be {\bf hyperbolic} if $d_X$ is a
metric on $X$, defining its usual ``coherent" (or ``metric") topology.
The unit disc and the curves of genus 2 or more which it covers are
hyperbolic.

A fundamental result of R.Brody states that a compact complex
manifold $X$ is hyperbolic if and only if every holomorphic map
$h:\mathbb C\ra X$ is constant.

The study of projective hyperbolic manifolds is guided by the two
nearly converse conjectures of S. Kobayashi and S. Lang:
Kobayashi's conjecture states that a hyperbolic projective manifold
has ample canonical bundle (and is so of general type). Lang's 
conjecture states that
a manifold $X$ is of general type if and only if it is 
``pseudo"-hyperbolic (in Lang's terminology).
  Which means that there is a proper algebraic subset
$S\subset X$
such that any
holomorphic map $h:\mathbb C\ra X$ has image contained in $S$. See
[La 1,2] for more information,
also on the relations with arithmetic geometry.

Notice that being
pseudo-hyperbolic in the above sense is a priori weaker than being
hyperbolic modulo $S$, in Kobayashi's terminology (which means that
$d_X$ is a metric outside $S$). Although it seems plausible that both
properties are
equivalent.

\subsection{\bf Hyperbolically Special Manifolds}

Our aim here and in section \ref{splitd} below is to give a conjectural 
description of $d_X$ in the
general case.

\begin{definition}\label{hspecdef} $X\in \mathcal C$ is said to be {\bf $H$-special}
($H$ for Hyperbolically) if $d_X\equiv 0$.
\end{definition}

 Our investigations are guided by :

\begin{conjecture}\label{conj3h}{\bf (Conjecture III$_H$):} Let $X\in \mathcal C$. Then $X$ is special if
and only if it is $H$-special.
\end{conjecture}

We shall give a motivation for conjecture \ref{conj3h} in section \ref{motconjd} below.

In Section \ref{splitd}, we shall complete conjecture III$_H$ by conjecture IV$_H$, which
(hopefully) describes $d_X$ for any $X\in \mathcal C$.

\begin{examples}\label{hspecex} We list special manifolds known to be $H$-special.
\end{examples}

{\bf (1) Curves:} A curve is $H$-special iff special, 
that is: either rational or elliptic.

{\bf (2) Rationally connected manifolds} are $H$-special (and
special). This is obvious. More generally:

\begin{proposition}\label{rcorbd} Let $f:X\ra Y$ be a rationally connected 
fibration (ie: so is its generic fibre), with $Y$ projective.
Then:  $d_X=f^*(d_Y)$. In particular, $X$ is H-special iff so is $R(X)$.
\end{proposition}

{\bf Proof:} Clearly, there exists a pseudometric $d_Y^+$ on $Y$ such 
that $d_X=f^*(d_Y^+)$.

We need to show that $d_Y^+=d_Y$. The inequality $d_Y^+\geq d_Y$ is obvious
(by the decreasing property of $d_X$ under holomorphic maps).
By the algebraic approximation theorem of [D-L-S], it is sufficient
to show this when $Y$ is a curve. But then, by [G-H-S], $f$ has a section.
This provides us with the opposite inequality $\square$

\begin{remark} This is another instance in which it is shown that 
rationally connected manifolds are in many aspects ``negligible", also in fibrations: 
indeed we already saw that if $f:X\mero Y$ is a rationally connected fibration with $X$ projective, then 
$X$ and $Y$ have the same core, the same fundamental group, and the same 
Kobayashi (pseudo-)metric. Only the arithmetic aspects are still unknown. 
Although the function field case is, by [G-H-S], much better understood.

Notice that, by the birational invariance of $d_X$, the above result applies 
if $Y$ is Moishezon. In fact, as in \ref{RCmultfree}, it should still be true when $X\in \cal C$.
\end{remark}

{\bf (3) Manifolds with trivial canonical class} They provide
important examples of manifolds with $\kappa=0$. Then conjecture III$_H$
is obviously true for complex tori, true for K3-surfaces (see below).

Using Bogomolov decomposition theorem, it is sufficient to
test the irreducible holomorphic symplectic
manifolds and the Calabi-Yau manifolds.

The known symplectic manifolds form a very short list (The two series
of Beauville-Fujiki, and the two examples
of O'Grady). For all these examples, a direct verification seems
quite feasible, by reduction to the case of abelian surfaces and
K3-surfaces.

For others symplectic manifolds (if any), the method of [C90],
obviously underexploited there, gives a very positive indication.

The first significant example of Calabi-Yau manifolds, that of
quintics in $\Bbb P_{4}(\Bbb C)$, might be accessible also, by a
strategy explained to me by C.Voisin.

{\bf (4) Manifolds with $\kappa=0$:} This is a crucial case, motivated 
by the special case where $K_X$ is trivial; the
conjecture then implies (by \ref{k=0}) the vanishing of $d_X$.
This question is standard, and was already raised in [Kob 76].

{\bf (5) Finite \'etale covers} of an $H$-special $X$ are
$H$-special, too (by the invariance of Kobayashi pseudometric under
\'etale covers).

{\bf (6) $H"$-special manifolds:}

Let us now consider the opposite implication (that H-specialness 
should imply specialness).
Its difficulty lies in the fact that the vanishing of $d_X$ has no 
known geometric translation
in dimension $3$ or more. See \ref{varhspec} for a simple discussion of this 
topic. So we start with a stronger hypothesis that gives such an implication.

\begin{definition} We say that $X\in \mathcal C$ is
$H"$-special if there exists a holomorphic map
with dense image (in the metric topology) $\phi:U\ra X$,
where $U=W\times \mathbb {C}$, with $W$ quasi-projective , and
moreover such that $\phi(w,0)=a$, for some fixed $a\in X$ and any
$w\in W$.
\end{definition}

Then $X$ is obviously $H'$- and $H$-special if it is $H"$-special. 
By \ref{h'spec}, $X$ is also special.

  In particular, of course, manifolds $X$ which are
$\mathbb {C}^d$-dominable or covered by $\mathbb {C}^d $ are
$H$-special.They are special, too by \ref{cddomspec}

\subsection{\bf Similarities}\label{sim}

{\bf Three stability properties:}\label{stabprop} We list here, to stress
their analogies, three stability properties which are common
to the classes of special and $H$-special manifolds:

(a) Let $f:X\ra Y$ be a fibration, and $g:V\ra X$ be
holomorphic with $f\circ g:V\ra Y$ surjective. Then $X$ is special
(resp.$H$-special)
if so are $V$ and the generic fibres of $f$.

(b) Assume $X$ is $S$-connected (resp.$H$-connected). In this case,
$X$ is special (resp.$H$-special). Here being $S$-connected
(resp.$H$-connected)
   means that two generic points of $X$ can be joined by a chain of
subvarieties of $X$ which are special (resp. $H$-special).

(c) The Albanese map is surjective if $X$ is special (resp. $H$-special, 
provided its Albanese torus is simple). See \ref{albsurj} (resp.\ref{hspecalbsurj} below) 
for the special (resp. $H$-special) case).

Notice that Conjecture III$_H$ implies 
that specialisations of special varieties are special, 
because this specialisation stability holds for $H$-special varieties, by the 
metric continuity of the Kobayashi metric.

\subsection {\bf The Origin of Conjecture III$_H$}\label{motconjd} Conjecture III$_H$
is motivated by the following:

\begin{question}\label{dspecfib} Let $f:X\ra Y$ be an admissible fibration.
Assume that $d_F\equiv 0$ if $F$ is a general fibre of $f$.
Then $d_X=f^*(d_Y^+)$, for a certain pseudometric $d_Y^+$ on $Y$.

(a) Is then $d_Y^+=d_{(Y/\Delta(f)})$ the Kobayashi pseudometric of
the orbifold $(Y/\Delta(f))$ (defined in \ref{orbd} below)?

(b) Is $d_{(Y/(\Delta(f))}\equiv 0$ if either $\kappa(Y/\Delta(f))=0$, or if
there exists a family of rational curves $(C_t)_{t\in T}$ with 
$X_t:=f^{-1}(C_t)$ special
for $t\in T$ general, and such that any two points of $Y$
can be connected by a finite chain of such $C_t$'s?\\
\end{question}

Notice that the above theorem \ref{rcorbd} shows that the preceding question
has a positive answer for rationally connected
fibrations $f:X\ra Y$, if $Y$ is projective, because then $\Delta(f)$ is empty.

Let us show briefly how to deduce Conjecture III$_H$ from a
positive answer to the above question, assuming moreover Conjectures 
$\kappa=0,-\infty$ expressing the core as a composition of orbifold Iitaka 
fibrations and rational quotients.
We keep the notations of \ref{c=rj^n}.

We argue by induction on the largest number $k\geq 0$ such
that $s_X^k:X\ra S^k(X)$ is not the constant map:
if $( s_X^k)_R$ is the constant map, then: there
exists a family of rational curves $C_t$ as in the question (b) above.
By the positive answer to that question, the orbifold
  pseudometric on the base vanishes, while the fibres
have vanishing psudometric by the induction hypothesis.
Hence the conclusion, by (a). If $g:=r\circ s_X^k$ is not constant, then
$s_X^{k+1}$ is constant. The preceding argument,
applied to $JB_g$, instead of $( s_X^k)_R$
gives the conclusion.

   The following result, due to [B-L] and [L] can be rephrased as giving
a second example (abelian fibrations over curves) in which the above 
question \ref{dspecfib} has a positive answer:

\begin{theorem} Let $f:X\ra C$ be a fibration with
$C$ a curve, $X$ projective, and the generic fibres of $f$ being
abelian varieties. Then \ref{dspecfib} above has a positive answer. 

More precisely:
$d_X=c_X^*$($d_C^*$), where $d_C^*$ vanishes
if $f$ is not of general type, and otherwise is a metric on $C$ such that
$g^*(d_C^*)=d_{C'}$, where $g:C'\ra C$ is a finite Galois cover
with $C'$ of genus 2 at least, the map $g$ having ramification order
exactly $m_i$ above each of the points $c_i$ of $\mid\Delta(f)\mid\subset C$.
\end{theorem}

\subsection{\bf  Surfaces}

\begin{proposition} \label{hspecsurf}A compact K\"ahler surface not of general type is
$H$-special if and only if it is special, and if and only if its fundamental group
is virtually abelian.
\end{proposition}

 Thus conjecture III$_H$ holds for them, except for
potential surfaces of general type and vanishing $d_X$, the
existence of which is excluded by Lang's conjectures. 

For surfaces, conjecture III$_H$ is thus equivalent to Lang's conjecture.

{\bf Proof:} This follows easily from the results of [B-L] (which can be rephrased 
as saying that non-$K3$ special surfaces are $\mathbb C^2$-dominable) 
and the fact that
every $K3$-surface is $H$-special, by the metric continuity of
Kobayashi pseudometric used in conjonction with the usual argument of approximating an
arbitrary $K3$-surface by Kummer ones. (More generally, being
$H$-special is a property which is closed under deformation in the metric
topology). Thus special surfaces are $H$-special.

The reverse direction is clear, because a non-special surface not of general type 
has (by \ref{coresurf}) a finite \' etale cover which maps surjectively onto a curve of 
genus two or more, and cannot thus be $H$-special $\square$

\begin{remark} Let $X\in \mathcal C$ be $H$-special with
$\kappa(X)=0,1$. Then $X$ is special. Indeed, the case $\kappa=0$ (resp. $\kappa=1$) 
follows from \ref{k=0} (resp. \ref{fibcurv}).

\end{remark}.

\subsection{\bf Threefolds}

Our results are very partial. The complete solution of Conjecture III$_H$
in this case depends in particular on the solution of several other
difficult conjectures:
Lang's conjecture in dimension 3, its orbifold generalisation (as
stated in section \ref{orbd} below) for orbifold surfaces, and finally the case where
$\kappa=0$.

\begin{proposition}\label{hspec3f} Let $X$ be a projective threefold.

(1) Assume $\kappa(X)\neq 3$, $ess(X)\neq 2$, and
$X$ is $H$-special. Then $X$ is special.

(2) Assume $\kappa(X)\neq 0,2$, and $X$ is special. Then $X$ is
$H$-special, unless possibly in the following case: $\kappa(X)=1$ and 
$J_X$ has generic fibres wich
   are $K3$-surfaces.
\end{proposition}

Let us remark that the non-algebraic case presents some specific
difficulties which can however be solved at the cost of some lengthy
arguments. The case
$\kappa(X)=2$ left open in the above statement seems also accessible
by known methods. We hope to come back on this in a later work.

The proof that the special $3$-folds listed above are $H$-special 
does not present any difficulty, 
using the description of special $3$-folds given in 
\ref{spec3f}. The converse direction rests on the fact that if 
$f:X\mero Y$ is surjective with $X$ $H$-special, then $Y$ is $H$-special too.

\subsection{\bf The Fundamental group}

Conjectures \ref{abconj} and III$_H$ imply that if
$X$ is $H$-special, then its fundamental group should be virtually
abelian.

It is thus natural to conjecture independantly that being $H$-special
should force the fundamental group to be virtually abelian.

 Notice
that
this would imply for example that being $H'$-special, or
   $\mathbb {C}^d$-dominable, or being covered by $\mathbb {C}^d$
implies the virtual abelianity of the fundamental group. Recall  that
this last
implication, by far the easiest one, is equivalent to S. Iitaka's
conjecture, and is widely open. 

A typical example of a question
suggested by the preceding context is:

\begin{question} Let $f:\Bbb C\ra X$, $X\in \mathcal C$, be
a holomorphic map with Zariski dense image. Is then $\pi_1(X)$ almost abelian?
\end{question}

\subsection{\bf The Albanese map}\label{hspecalb}

  However, in complete analogy with the special case (see 
\ref{albsurj}), we have, using
strengthenings of a classical result of Bloch-Ochiai (surprisingly,
it does not seem
to have ever been stated that way):

\begin{proposition}\label{hspecalbsurj} Let $X\in \mathcal C$ be
$H$-special. Then its Albanese map $\alpha_X:X\ra Alb(X)$ is
surjective  if
$Alb(X)$ is a simple torus
\end{proposition}

Of course, the hypothesis that $Alb(X)$ is simple is artificial. It
could be removed, using another result of Y.Kawamata in [Kw"],
provided one
could positively answer question \ref{hypmod} below, which
is a natural strengthening of Brody's theorem.
(Alternatively, it is possible that a suitable sharpening of the
proof of Bloch-Ochiai's theorem
permits to remove this artificial hypothesis, without answering question \ref{hypmod} in general):

\begin{question}\label{hypmod} Let $X$ be a compact complex
connected manifold. Let $Z\subset X$ be a proper compact analytic
subset. Assume
that any holomorphic map $h:\mathbb C \ra X$ has image contained in
$Z$. Is then the Kobayashi pseudo-metric $d_{X}$ a metric on its
complement $X/Z$?
\end{question}

Notice also that
certain recent results of Noguchi-Winkelmann  (or at least the
methods there) might imply that in \ref{hspecalbsurj}
the Albanese map considered is surjective and has no multiple fibres
in codimension one,
as expected from Conjecture III$_H$ and Theorem \ref{albsurj}

\begin{corollary} If $X\in \mathcal C$ which is $H$-special (and with
simple Albanese torus), then:

(a) any torsionfree solvable quotient of $\pi_1(X)$ is almost abelian.

(b) any linear quotient of $\pi_1(X)$ is almost abelian.

\end{corollary}

{\bf Proof:} The proof is exactly the same as that of \ref{lin} (where 
 $X$ was special), just observing that any finite \'etale cover of $X$ satisfy the same 
hypothesis. Also we need to apply Zuo's result ([Z]) that if $\pi_1(X)$ had 
a non-trivial semi-simple linear representation, 
then $X$ were generically Brody hyperbolic (in particular: not $H$-special) $\square$

\subsection{\bf Geometric variants of $H$-speciality}

  Brody's theorem says that the obstruction to the hyperbolicity
of a compact $X$ is due to the presence of ``entire" curves in $X$,
that is images of nonconstant holomorphic maps $h:\mathbb C\ra X$. It
is thus natural to try to express the fact that $X$ is $H$-special by
the geometry of ``entire" curves lying on $X$. In this sense, being
$H$-special appears as a transcendental analogue of rational connectedness.

Unfortunately, the analogue of a deformation theory for entire curves
does not exist (and would certainly present great difficulties).
So we introduced above the notions of
$H'$-speciality, of $H"$-speciality, and of $\mathbb
{C}^d$-dominability to try to give some examples in which the 
vanishing of $d_X$
is deduced from the geometry of ``entire" curves.

\begin{remark}\label{varhspec} Let $X\in \mathcal C$. One has the
following implications between the different notions of speciality:
$\mathbb {C}^d$-dominability implies $H"$-speciality, which implies
both $H$-speciality and $H'$-speciality, which
implies  speciality (by \ref{ko} ).
\end{remark}

Of course, there are lots of other geometric variants of
$H$-speciality which one can introduce, replacing for example
the metric topology by the Zariski
topology, and/or considering jets (see [Gr]), and which can be
considered as approximations of $H$-speciality. One can ask about the 
relationships
between the various notions so defined. In particular, do the notions
relative to metric and Zariski topologies always coincide?

We shall just ask one question
of whether some strong reverse implications hold between the various notions 
introduced in \ref{varhspec}.

\begin{question} Assume $X$ is $H$-special. Is  $X$ then also
$H"$-special?
\end{question}

  This question is motivated by the fact that a positive answer
(which does not seem unrealistic, although far beyond reach if true)
would imply, using \ref{ko}, that
$H$-special manifolds are special, which is one of the two
implications of conjecture III$_H$.

 It is shown in [B-L] that for surfaces not
of general type, $\mathbb C^d$-dominability, $H'$-speciality and
$H$-speciality are equivalent.

This may however be a low dimensional phenomenon, similar to
the fact that
rational connectedness, unirationality and rationality coincide for 
complex surfaces.

This suggests the following obvious observations concerning
rationally connected manifolds:

\begin{remark}\label{rcunir} Let $X$ be a rationally connected
manifold. Then $X$ is $H"$-special.
If $X$ is unirational, then $X$ is
$\mathbb {C}^d$-dominable.
\end{remark}

  This suggests of course also the following (again far beyond
reach, since no rationally connected manifold is known which is not
unirational, although most of them are suspected not to be):

\begin{question}\label{cddominunir} Let $X$ be a rationally connected manifold.
If $X$ is $\mathbb {C}^d$-dominable, is $X$ unirational?
\end{question}

The
answer is not expected to be necessarily positive. It seems quite possible
that the class of rationally connected manifolds which are
$\mathbb {C}^d$-dominable is strictly intermediate between the other
two ones (unirational and rationally connected), as the following 
example may suggest:

\begin{example} The threefold $X$ constructed in [U,(11.7.1),p.137]
is obviously $\mathbb C^3$-dominated, and (less
obviously) rationally connected. But it is asked by K.Ueno whether it is
unirational or not.
\end{example}

Recall that $X$ is nothing but a smooth model of
the quotient of $E^3$ , with
$E$ the Gauss elliptic curve, by the cyclic group of order 4
generated by complex multiplication by $i:=\sqrt {-1}$ on each
factor.

In the other direction (discrepancy between rational connectedness
and $\mathbb {C}^d$-dominability), even the $\mathbb
{C}^3$-dominability
of the general quartic in $\mathbb {P}^4$ is far from obvious, if true.

\subsection{\bf The Core and the Splitting of the Kobayashi Metric}\label{splitd}

Let $X\in \mathcal C$, and let
$c_X:X\ra C(X)=Y$ be its core. By conjecture III$_H$, its general fibres are
$H$-special, since special. Thus $d_X$ should vanish on all of them, 
by Conjecture III$_H$ and
the metric continuity of $d_X$. There should thus exist
a unique pseudometric
$d_{C(X)}:=d_Y^*$ on $Y=C(X)$ such that $d_X=c_X^*(d_Y^*)$.

We give a first version of:

\begin{conjecture}{\bf Conjecture IV$_H$:}

The pseudometric $d_{C(X)}^*$ is a metric on some
Zariski nonempty open subset $C(X)^*$ of $C(X)$.
\end{conjecture}

\begin{remark}
\end{remark}

(1) Assume $X$ is of general type. Then $c_X=id_X$, and we
just recover the usual conjecture of Lang in its strong form.

(2) Recall that, after Conjecture II, $(C(X)/\Delta (c_X))$
is an orbifold of general type. Thus Conjecture IV$_H$ appears as
an extension of Lang's conjecture to the orbifold case.

(3) The preceding remark can be made more precise, if one
introduces as in \ref{orbd} below an orbifold  Kobayashi
pseudometric on any orbifold $Y/\Delta$. One may then ask whether
$d_{C(X)}^*=d_{C(X)/\Delta (c_X)}$. Then Conjecture IV$_H$ appears as
the orbifold extension of Lang's conjecture. See conjecture IV'$_H$
below. We hope to come back on this latter.

We shall now propose an orbifold definition of the Kobayashi pseudometric:

\begin{definition}\label{orbd} Let $X$ be a complex manifold and
$\Delta=\sum (1-1/m_{i}).\Delta_{i}$ be a $\Bbb Q$-divisor on $X$,
   with $m_i$ positive integers and
$\Delta_i$ irreducible distinct divisors.
  Let $\Bbb D$ be the unit
disc in $\Bbb C$, and $h:\Bbb D\ra X$ a holomorphic map.

We say that $h$ is an {\bf orbifold map} to $(X/\Delta)$ if, for
each $z\in \Bbb D$ such that $h(z)$ belongs to the smooth locus
of  $\mid\Delta\mid$, say to
$\Delta_i$, for some $i$, then for any nonnegative integer $j$ which
is not divisible by $m_i$, the
  $j^{th}$ derivative of $h$ at
$z$, seen as a tangent vector to $X$ at $h(z)$, is tangent to
$\Delta_i$ at $h(z)$.

We then define $d_{(X/\Delta)}$ as the largest pseudometric on $X$
such that $h^*$$(d_{(X/\Delta)})$$\leq d_{\Bbb D}$, for any orbifold
map from $\Bbb D$ to $X$.

We call $d_{(X/\Delta)}$ the {\bf orbifold Kobayashi pseudometric of
$(X/\Delta)$}.
\end{definition}

(I thank J.P. Demailly for supporting my impression that such a
definition is meaningful, and correcting my first attempt).
   Remark that by the same arguments as in the classical case, this
pseudometric enjoys the same
properties as when $\Delta$ is empty (in which case we just recover
  the usual Kobayashi pseudometric on $X$). 

In particular, we have the
following obvious refinement of the distance decreasing property of
holomorphic maps:

\begin{proposition}  Let $f:X\ra Y$ be a fibration.
Assume that $\Delta(f)=\Delta^-(f)$. 
Then: 

$d_X\geq f^*(d_{(Y/\Delta(f)})$.
\end{proposition}

{\bf Proof:} For any holomorphic map $h:\Bbb D\ra X$, $f\circ h:\Bbb
D\ra Y$ is an orbifold map to $(Y/\Delta(f))$.

(This trivial argument
is in fact
the motivation for the definition \ref{orbd} given above).

Notice that for the reverse inequality, (at least when $f$ is a special fibration) 
one would need to approximate
arbitrarily closely holomorphic orbifold maps from $\Bbb D$ to $Y$ by
liftable (to $X$) ones. The definition of multiplicities (based on $inf$) shows that 
there is no local (on $Y$) obstruction for this.

We can now formulate more precisely Conjecture IV$_H$:

\begin{conjecture}{\bf Conjecture IV'$_H$:} Let $c_X$$:X\ra C(X)$ be the core of
$X$, with $c_X$ admissible. Then:

(1) $d_X=(c_{X})^*$$(d_{(X/\Delta(c_{X})})$. (In other terms: $d_X$
should be nothing but the lift of the orbifold pseudometric of
$(X/\Delta(c_X))$ by $c_X$).

(2) Let $(X/\Delta)$ be a klt orbifold of general type. Then its
Kobayashi pseudometric is a metric outside some proper algebraic
subset $S\subset X$.
\end{conjecture}

Remark that Conjecture IV'$_H$ obviously implies Conjecture IV$_H$.

\section{\bf Arithmetic Aspects}

 \subsection{\bf Fields of Definition}

We consider here the situation in which $X$ is a projective manifold
defined over the field $K\subset \Bbb C$, supposed to be
finitely generated over $\Bbb Q$ (this can always be achieved).

\begin{proposition} Let $f:X\ra Y$ be an almost
holomorphic fibration. Then:

(a) $f$ is defined over some finite extension $K'$ of $K$.

(b) If $f$ is either $J_X$ the Iitaka fibration, or $r_X$ the
rational quotient, or $c_X$ the core of $X$, then $f$ is defined
over $K$ itself.
\end{proposition}

{\bf Proof:} I would like F.Bogomolov for explaining me the proof and
indicating me some of the examples below.

{\bf (a)} The family of fibres of $f$ (with their reduced
structure for the generic one) is then an irreducible component
of the Hilbert Scheme $Hilb(X)$ of $X$. Let $\overline K$ be the
algebraic closure of $K$. Then $G:=Gal(\overline{K}/K)$ acts
on the set of irreducible components of $Hilb(X)$ with finite orbits,
since this action preserves the Hilbert polynomial of any component.
This shows the claim.

{\bf (b)}  When $f=J_X$, the claim is clear,
since $J_X$ is defined by the linear system $\mid m.K_X\mid$ for some
suitable integer $m>0$,
and this system is invariant under the Galois action.

In the same way, because of the maximality properties of the general
fibres of $r_X$ and $c_X$, we see that $G$ preserves the component of
$Hilb(X)$
parametrising the generic fibres of either $r_X$ or $c_X$.
More precisely: the conjugate of a rational curve (resp. a special
variety) is still rational (resp. special).
(Notice that the same argument in fact applies to any fibration
enjoying a similar maximality property) $\square$

\subsection{\bf Arithmetically Special Manifolds}

\begin{definition} \label{aspecdef}The projective manifold $X$ is said to be
{\bf Arithmetically special} ($A$-special, for short), if its set $X(K')$ of
$K'$-rational points is Zariski dense in $X$, for some suitable 
finite extension
$K'/K$. (This notion is usually called ``potential density").
\end{definition}

In analogy with the notion of $H$-speciality, we formulate the:

{\bf Conjecture III$_A$:} $X$ is special if and only if
it is $A$-special.

(In the first version of the present article, I hesitated to include this
conjecture; but the independent kind suggestion from B.Totaro
to include it decided me. I thank him for this suggestion).

One can strengthen the above conjecture by replacing the Zariski 
topology by the
metric topology. Examples in [C-T-S-SD] show that further
field extensions may be needed (ie: $X(K')$ may be Zariski
dense without being metrically dense).

Now guided by Lang's arithmetic conjecture that a variety of general
type should be ``Mordellic", and the connections he found between
hyperbolicity and
``Mordellicity", we infer from the observations above in \ref{komet}, the following
definitions and conjectures:

{\bf Conjecture IV$_A$:} Let $X$ be a projective manifold defined
over the field $K\subset \Bbb C$, assumed to be finitely generated
over $\Bbb Q$.

   Let $c_X:X\ra C(X)$ be its core (which is defined over $K$, by the
above remark). Let $K'/K$ be any finite extension.
  Then the $K'$-rational points of $X$ are mapped by
$c_X$ to a proper
algebraic subset which is finite outside of
the proper algebraic subset $S$ of $C(X)$ defined in Conjecture IV'$_H$.\\

{\bf Known cases that Special implies $A$-special:}

{\bf (a) Curves.} For genus 0 curves, the existence of a $K$-rational 
point implies already the
density of $X(K)$. For genus 1 curves, this is a special case of 
[L-N], for example.

{\bf (b) Abelian varieties.} This is shown in [L-N], but is older 
(going back to A.Weil at least).

{\bf (c) Unirational manifolds.} The conjecture holds, almost by 
definition. For rational surfaces,
see [CT], especially th.5 there, as well as the references in 
[CT-S-SD] and [CT], for example.
  (The works devoted to this
class of surfaces lie much deeper than what is needed for our purpose 
here, because the
authors usually work with a {\bf fixed} base field $K$).

{\bf (d) Special elliptic surfaces with two sections.}(due to [B-T])The base is a 
rational or elliptic curve. So
(finitely) extending the base field, we can assume everything to be 
defined over $K$, and the two sections
having a dense set of $K$-points. They induce then, using addition in 
the fibres of the elliptic fibration a
dense set of $K$-rational points.

 This is the principle of 
construction of [Bo-T] for a $K3$-surface
over a number field, and a dense set of $K$-rational points. By 
taking a suitable base change by an
elliptic double cover of the base $\Bbb P^1(\Bbb C)$, one can even 
get a surface $X$ with $\kappa(X)=1$, and $X(K)$ dense. Of course, 
$X$ is special.

(Colliot-Th\'el\`ene informed me that $K3$ examples with this 
properties are already in [SD], and with $\kappa=1$
in [Cox]).

{\bf (e) Quartics in $\Bbb P^4(\Bbb C)$.} These are $A$-special if 
defined over a number field (by [H-T]).

{\bf (f) Fano Manifols, Rationally connected Manifolds.}
The answer is not known (although the analogous property property is 
known over function fields  (see \ref{rcffspec} below).

Let us list some properties common to special and $A$-special manifolds.

First, in complete analogy with \ref{albsurj} and \ref{hspecalbsurj}:

\begin{theorem}\label{aspecalbsurj} Let $X$ be $A$-special. Then:

{\bf (a)} The Albanese map of $X$ is surjective.

{\bf (b)} Any finite \'etale cover $u:\widetilde X\ra X$ of $X$ is 
$A$-special (if $X$ is defined over a number field $K$).
\end{theorem}

{\bf Proof:} (a) If not, we have a surjective map $\beta:X\ra Z$, 
where $Z$ is a subvariety of general type
of an Abelian variety. Everything is defined over the base field $K$, 
after a finite extension.
After [Fa], $Z(K)$ is however not Zariski dense. This contradicts our 
assumption that $X$ is $A$-special.

(b) This is an immediate consequence of the Chevalley-Weil Theorem 
([H-S],Ex.C7, p.292), which states that
in this situation, $\widetilde X(K')$ contains $u^{-1}(X(K))$, for a 
suitable finite extension $K'/K$. (The restriction
that $K$ is a number field can probably be removed, here).

{\bf Known cases of conjecture IV$_A$:}

The list here is very short.

{\bf (a) Curves.} This is known by the solution of Mordell's 
conjecture in [Fa] for
curves of general type and previous papers by the same author.

{\bf (b) Subvarieties of Abelian varieties.} Such a variety is not 
$A$-special if of general type, by [Fa].
The distribution of $K$-rational points is as predicted by Lang's conjectures.

{\bf (c) Non-special surfaces not of general type.}

\begin{theorem} Let X be a projective surface defined over 
the number field $K$. Assume $X$ is not
special and not of general type. Then conjecture IV$_A$ holds for $X$.
\end{theorem}

{\bf Proof:} Because of the hypothesis, either $X$ is birationally 
ruled over a curve of genus 2 or more, and
the conclusion follows from [Fa], or $X$ has an elliptic fibration 
which is of general type (in the sense of
(1.9)). But then $X$ has a finite \'etale cover $\widetilde X$ which 
is elliptic over a curve of general type.
The conclusion now follows from combining [Fa] and \ref{aspecalbsurj} (the 
Chevalley-Weil theorem).$\square$

(Notice that the preceding argument in fact solves, more generally, 
conjecture IV$_A$ when $ess'(X)=1$ over number fields:

\begin{proposition} Assume $ess'(X)=1$, and $X$ is defined 
over a number field. Then Conjecture IV$_A$ holds for $X$.
\end{proposition}

The hypothesis here means that the core of $X$ is a fibration of $gcd$-general type 
 (not just of general type). See \ref{gcdmult} for this notion. 
It is an interesting question to know 
if the arguments can be adapted to the case where $ess(X)=1$.

\begin{remark}
\end{remark}

{\bf (1)} The preceding result is known in various versions (see 
[Da], [Da-G], [CT-S-SD]. I thank 
D. Abramovitch and J.L. Colliot-Th\'el\`ene for learning me these references).

{\bf (2)} Presumably, the conclusion extends if $K$ is finitely 
generated over $\Bbb Q$.

{\bf (3)} Lang's connection between hyperbolicity and arithmetic properties
and the reformulation above of Conjecture IV$_H$ as Conjecture IV'$_H$  also suggest strongly 
that one should be
able to attach to any projective orbifold $(Y/\Delta)$ defined over 
$K$ the notion of its
$K$-rational points $(Y/\Delta)(K)$. 

Such a definition can actually be given quite naturally, in 
complete analogy with the hyperbolicity and function field cases considered here, as explained to me by 
E.Peyre ([Pey]).

An important property is that if $f:X\ra Y$ is a fibration defined over $K$, then 
$f(X(K))\subset (Y/\Delta(f))(K)$, the latter being the set of $K$-rational 
points of the base orbifold $(Y/\Delta(f)$ of $f$. 

Then Conjecture IV$_A$ above can be reformulated as the assertion that
an orbifold of general type defined over $K$ does not have a Zariski
dense set of $K-$rational points.

{\bf Conjecture IV'$_A$}: Let $c_X: X\ra C(X)$ be an admissible 
representative of the core of $X$, defined over $K$, finitely generated over $\mathbb Q$. 

If $K'/K$ is any finite extension $S\subset C(X)$ is defined as in 
conjecture IV'$_H$, then $(C(X)/\Delta(c_X))(K')$ 
has only finitely many points outside of $S$ (and contains $c_X(X(K'))$).

\subsection{\bf A Related Conjecture of Abramovitch and Colliot-Thelene}

A conjecture related to conjecture III$_A$ is formulated in [H-T],
and attributed to D. Abramovitch and J.L. Colliot-Th\'el\`ene 
(I thank J.L. Colliot-Th\'el\`ene for learning me this reference)
: it claims that if $X$ is what we called $w$-special and 
defined over
a number field, then $X$ is $A$-special (``potentially dense", 
there), and conversely.

This conjecture thus coincides with conjecture III$_A$ if the 
properties of being special and $w$-special coincide.

Although no example of a $w$-special, but not special manifold is known today,
it seems plausible that the two properties are distinct (special is 
anyway more restrictive).

The difference in the two conjectures is that if there exists an $X$ 
(defined over a number field) which is $w$-special, but not special, it
should be $A$-special by [H-T], but not by IV$_A$, which would imply 
that its $K'$-rational points should lie in finitely many fibres of $c_X$, 
outside of $c_X^{-1}(S)$.

\section{\bf The Function Field Case}

\subsection{\bf  The (Complex) Function Field Case}

The following extremely interesting function field variant of
Conjectures III$_A$ and IV$_A$  above was
suggested to me by P. Eyssidieux. They are included here with his 
kind authorisation.

(One could of course consider other fields than $\Bbb C$ as well).

We need first some conventions.

{\bf Conventions.} Let $f:X\ra C$ be a fibration onto a 
complex projective
curve $C$, with $X$ complex projective, smooth connected. Let 
$S^0(f)\subset Chow(X)$ be the
Zariski open subset consisting of holomorphic sections of $f$, and 
let $S(f)$ be
its closure in $Chow(X)$, the Chow variety of $X$.

We also denote by $S^0(X)\subset X$ the union of all images of 
sections of $f$, and
by $S(X)$ the Zariski closure of $S^0(X)$ in $X$. We say that an 
irreducible component
$\Sigma$ of $S^0(f)$ is $X$-covering if the union of all images of 
sections in $\Sigma$
is Zariski dense in $X$ (and so covers a nonempty Zariski open subset of $X$).

Let finally $b:C'\ra C$ be any finite covering (possibly ramified), 
with $C'$ a smooth
connected projective curve. We then let $f':X'\ra C'$ be deduced from 
$f$ by the base change $b$.
We define also $S(f')$, $S(X')$,... as before, with $f'$ in place of $f$.

The arithmetic interpretation is that if $K$ is the (non algebraically closed)
function field of $C$ (ie: its field
of meromorphic functions), then $X$ is a model of the $K$-variety 
$X_K$ defined over $K$
by the function field of $X$. Moreover, the sections of $f$ are exactly the
$K$-rational points $X(K)$ of $X_K$; and $X(K)$ is Zariski dense in 
$X_K$ if and
only if so is $S^0(X)$, or equivalently: $S(X)=X$.

The finite extensions $K'/K$ are exactly the ones given by the 
coverings $b:C'\ra C$ above.

We shall say that $f:X\ra C$ is {\bf split} if there exists a 
manifold $F$ and a birational map
$\sigma:F\times C'\ra X$ over $C'$, after some suitable base change 
$b$. This is (essentially)
Lang's terminology in [La 2].

We need to introduce the orbifold version of the preceding notions: 

An orbifold $(X/\Delta)$ over $C$ will be just a fibration $f:X\ra C$ as above, together 
with an orbifold structure on $X$. The notions of Kodaira dimension, canonical bundle, 
general type,... are defined by restriction to the general fibre. 

If $f:X\ra C$ and $g:Y\ra C$ are fibrations, a fibration $h:X\mero Y$ over $C$ 
is the usual $h$ such that $g\circ h=f$.

When $h$ is holomorphic, we can naturally define its base orbifold (an orbifold 
structure on $Y$) in the expected way. This leads to the notion of Kodaira 
dimension of $h$, denoted $\kappa(Y,f)/C$.

A special fibration $f:X\ra C$ is thus defined in the same way as we already did. 
By the existence of a relative core, this is the same as having a special general fibre.

In this context, it is maybe useful to define the sections of an orbifold $(X/\Delta)$ 
over $C$. These are simply the sections of $f$ which meet each component $\Delta_j$ of 
$\Delta$ at order at least $m_j$ at each point where the section meets $\Delta_j$, if 
$\Delta=\sum_{j\in J}(1-1/m_j). \Delta_j$.

We say that $(X/\Delta)$ is {\bf split} if it becomes birationally 
trivial after some base change $b$ as in the non-orbifold version. (Birationally trivial 
means that $(X/\Delta)$ is birational to $(F/D)\times C$ over $C$, for some ``constant" 
orbifold $(F/D)$).

\begin{definition}\label{ffspec} We say that the $K$-variety $X_K$ is {\bf 
special} (resp. {\bf rationally
connected}; resp. an {\bf abelian variety}; resp. {\bf of general 
type}...) if so is the general fibre of $f$.
\end{definition}

Let $f=g\circ h$ be factorisation of $f$ by fibrations $h:X\ra Y$ and 
$g:Y\ra C$. We say that this
factorisation is the $K$-core (resp. the $K$-rational quotient; resp. 
the $K$-Albanese map) if it induces
the corresponding reduction on the general fibre of $f$. These 
factorisations thus exist
(by \ref{relquot}) and analogous results for the other cases), are uniquely defined,
and moreover compatible with the base changes $b$.

Because of the split case, the formulation of the function field version
of conjectures III$_A$ and IV$_A$ is not exactly the same.

To formulate them we need to introduce
the $K$-core $f=g_f\circ c_f$ of $f$ (see \ref{relquot}). We thus have a 
factorisation of $f$ by $c_f:X\ra C(f)$ and
$g_f:C(f)\ra C$.

\begin{conjecture}\label{conj3ff}{\bf (Conjecture III$_{FF}$)}  The $K'$-rational points of 
$X_K$ are Zariski dense in $X$, for some
finite extension $K'/K$ if and only if $g_f:(C(f)/\Delta(c(f)))\ra C$ is split. In 
particular, if $X_K$ is special, $X(K')$ is
Zariski dense in $X$, for a suitable finite extension $K'/K$. (One might say, in analogy with the 
terminology introduced in \ref{hspecdef} and \ref{aspecdef}, that they are $FF$-special ($FF$ standing for: function field).
\end{conjecture}

Geometrically, this just means that $S(X')=X'$ if and only if $C(f)$ 
is birationally a product over $C$ (after a
preliminary covering of $C$, maybe).

\begin{conjecture}\label{conj4ff}{\bf (Conjecture IV$_{FF}$)} Assume that $g_f:(C(f)/\Delta(c(f)))\ra C$ is 
not split. There exists a Zariski closed
subset $T\subset C(f)$, $T\neq C(f)$ such that for any finite 
extension $K'/K$, $c_f(S(X'))$ is contained in $T$
and a finite union of images of sections of $g_f$. 
\end{conjecture}

\begin{remark} If $g_f$ is split, one can consider also only
the sections which are not constant to then recover exactly the same 
formulations as in the arithmetic case.
\end{remark}

{\bf Known cases when $X_K$ is special.} One wants to know if 
$X(K')$ becomes Zariski dense for a suitable
finite extension $K'/K$.

{\bf (a) Curves.} For rational curves, this is clear,
because such surfaces $X$ are birationally ruled. For elliptic
fibrations, it is a special case of the next case
(alternatively, by a suitable base change, construct 2 sections; the first one
is used as the zero section to define a relative group law;
then choose the second section so as to have infinite order.
This is an immediate construction).

{\bf (b) Abelian varieties.} This is shown in a more general form in [L-N].

{\bf (c) Rationally connected varieties.}\label{rcffspec}

They are $FF$-special (for ``function-field-special).
This immediately follows from
the glueing lemma of [Ko-Mi-Mo], which
implies that if the generic fibre of $f:X\ra C$ is rationally connected,
and if $f$ has a section, then any finite set of $X$ mapped injectively to
$C$ is contained in a section of $f$.

Notice that [G-H-S] announces that $f$ itself has a section, so that no
field extension is needed in that case.

{\bf (d) Surfaces.} The answer is known, except for $K3$-surfaces as fibres:

\begin{theorem} Let $f:X\ra C$ be a fibration with its 
generic fibre a special, but not $K3$, surface.
Then $S(X')= X'$ for a suitable finite covering $b:C'\ra C$.
\end{theorem}

The proof will be given elsewhere. The case of elliptic $K3$ surfaces
is true and immediate. It is certainly
worth investigating the case of general $K3$'s.

{\bf Known cases when $X_K$ is not special.}

The list is then much more limited.
It is in fact restricted to cases when $X_K$ is of general type,
and where conjecture IV$_{FF}$ reduces to the function field case of 
Lang's conjectures.

{\bf (a) Curves.} This case has been settled by various approaches by 
H.Grauert [Gr"], Y.Manin [Man] and  A.N.Parshin [Par].

{\bf (b) Subvarieties of general type of Abelian varieties.} This the 
case where the general
fibre $F$ of $f$ is of general type and mapped birationally into its 
Albanese variety
by its Albanese map.

We now turn to another question, quite natural in our context, but 
for which we don't know any reference.

\subsection{\bf Split and ``Dominated" Fibrations: the Rational Quotient}

We shall make here some easy observations on the {\bf non-orbifold versions}. 

It would be interesting to extend them to the orbifold case as well.

\begin{definition}\label{domfib} The fibration $f:X\ra C$ is said to be 
{\bf dominated} if
there exists a projective manifold $F'$ and a dominant rational map 
$g:F'\times C'\ra X'$ ovec $C'$,
  for some suitable $b:C'\ra C$.
\end{definition}

If the above map $g$ can be chosen to be regular (resp. isomorphic 
above the generic point of $C$),
  we say that $f$ is {\bf regularly dominated} (resp. {\bf isotrivial}).

Remark that if $f$ is dominated, the manifold $F'$
can be chosen to be of the same dimension as $F$, the generic fibre 
of $f$, so that $g$ is generically
finite. And remark also that a split fibration is dominated.

Our objective now is to introduce Conjecture III$_{SD}$ below, which says that
the obstruction to reverse implication lies in the rational quotient of $f$.

The relevance of notion of ``dominatedness" to our context is the following:

\begin{proposition} The following properties are 
equivalent, for a fibration $f:X\ra C$:

{\bf (a)} $f$ is dominated.

{\bf (b)} $S^0(X')= X'$, for a suitable $b:C'\ra C$. (See above for 
the definition of $S^0(X')$).

{\bf (c)} There exists an irreducible component $T'$ of $S(f')$ which 
is $X'$-covering, for some suitable
$b:C'\ra C$. (See \ref{geomquot}.)
\end{proposition}

{\bf Proof:} Obviously, (a) implies (b). To show that (b) implies 
(c), use the countability at infinity of the Chow-Scheme of $X$.
To show that (c) implies (a), take $T':=F'$, and use the (surjective) 
evaluation map of the incidence graph $G$, say,
  of the algebraic family of cycles parametrised by $T'$. The 
conclusion is implied by the fact that $G$ is bimeromorphic to 
$C\times T'$,
because the generic point of $T'$ parametrises a section of $f$.

\begin{conjecture}\label{sdomconj}{\bf (Conjecture III$_{SD}$)}

{\bf (i)} Let $r_K:X\ra R(X_K)$ be the rational quotient of the 
fibration $f:X\ra C$.
Then $X_K$ is dominated if and only if $R(X_K)$ is split.

A weaker form (shown below to be equivalent) is:

{\bf (ii)} Assume $X_K$ is dominated, but not split. Then $X_K$ is 
uniruled (ie: so is $F$).
\end{conjecture}

This means in particular that $X_K$ is dominated if rationally 
connected, and that it is split if it is
dominated and not uniruled.

We shall show several cases where the above conjecture holds. First, 
the ``if" part:

\begin{theorem} \label{splitdom}If $R(X_K)$ is split, then $X_K$ is dominated.
\end{theorem}

A special case (when $R(X_K)=C$) is:

{\begin{corollary} Assume $X_K$ is rationally connected. 
Then it is dominated.
\end{corollary}

It would be interesting to have concrete examples when $F$
is a three-dimensional cubic or quartic.

For example: the above
result asserts that one can uniformely dominate one-dimensional 
families of such
manifolds, although they are birationally distinct.

{\bf Proof (of \ref{splitdom}):} Assume $R(X)=C\times R$ (after finite base 
change, if needed). Let $r\in R$ be generic, and
let $r'_K:X_r:=r_K^{-1}(C\times \lbrace r \rbrace)\ra (C\times 
\lbrace r\rbrace )$ be the restriction of $r_K$ to $X_r$.
This is a rationally conected fibration, and so by [G-H-S] and 
[Ko-Mi-Mo], it is dominated. Thus so is $f$, using characterisation 
(b) above of dominatedness $\square$

\begin{theorem} If the weaker form (ii) of conjecture III$_{SD}$ 
 holds, the stronger form (i) also holds.
\end{theorem}

{\bf Proof:} The proof is actually a direct application of \ref{RCRG}, or of 
[G-H-S], which asserts the
equivalence of rational connectedness and rational generatedness. 
(This last notion is
defined as follows: $X$ is rationally generated if, for any fibration 
$h:X\ra Y$, $Y$ is
uniruled). Indeed: since $X_K$ is dominated, so is $R(X_K)$. If $R(X_K)$ is 
not split, it is uniruled by .
But [G-H-S] for details) in particular says that 
$R(X_K)$ is not uniruled.

If, conversely, $R(X_K)$ is split, we deduce from \ref{splitdom} above 
that $X_K$ is dominated $\square$

We now give some positive answers to the conjecture III$_{SD}$, (ii).

\begin{proposition} Assume that through the generic 
point of $X'$ there
exists a positive-dimensional family of sections of $f':X'\ra C'$. 
Then $X_K$ is uniruled.
\end{proposition}

This is a simple application of Mori's bend-and-break technique. Indeed:

{\bf Proof:} This is a simple application of Mori's bend-and-break technique. Indeed:
let $x\in X$ be generic with a 
positive-dimensional family of sections through it.
By Mori's techniques (no characteristic $p>0$ is needed, here), there 
is a rational curve $\Gamma$ in $X$ through $x$.
If the base $C$ is not rational, $\Gamma$ is vertical (ie: contained 
in the fibre of $f$ through $x$.
We can reduce to this case by a base change $b:C'\ra C$, with $C'$ 
non-rational, because sections of $f$ lift to sections of $f'$.$\square$

So the remaining case is when
$X$ is dominated, but the set of sections through the general point 
of any $X'$ is at most countable, not empty.

\begin{proposition} Assume $X_K$ is regularly 
dominated. Then $X_K$ is isotrivial.
\end{proposition}

This is essentially (up to the terminology) shown in [Mw].

\begin{corollary} Assume $F$ (the generic fibre of $f$) 
is an absolutely minimal
model (in the sense that any birational rational map $r:F_0\ra F$ , 
with $F_0$ smooth, is regular).
Then: $X_K$ is isotrivial if dominated.
\end{corollary}

\begin{example} The above condition is satisfied if $F$ is 
an abelian variety, or a
minimal surface with $\kappa(F)\geq 0$. 
\end{example}

Another important case is given by:

\begin{proposition} Assume $F$ is of general type. If 
$X_K$ is dominated, it is split.
\end{proposition}

This is shown in [Mae]. Observe that conjecture III$_{SD}$,(ii) 
consists in extending this to the case
where $F$ is not uniruled. Already the (conjecturally equivalent) 
case where $\kappa(F)\geq 0$ would be of
interest.}

\section{\bf Appendix: Geometric quotients} \label{geomquot}

We now turn to the exposition of various techniques of construction 
of quotients for equivalence relations of geometric origin.

These techniques have been introduced in [C80,81] and have found many 
subsequent applications ([C92], [Ko-Mi-Mo], [C94], [C01] among others).

They are based on fundamental properness properties of the Chow-Scheme 
$\mathcal C$$(X)$ of $X$ (see [Ba] if $X$ is not projective) if $X\in \mathcal C$, 
discovered in [Lieb] (see [Fu] for the Hilbert-Douady space case), further properties 
having been shown in [C80].

When $X$ is projective, these properties reduce to the well-known projectivity 
of the irreducible components of $\mathcal C$$(X)$ of $X$

The cycle of $X$
naturally parametrised by $t\in T$ is denoted by $V_t$. Conversely, if $V$ 
is a subvariety of $X$ (ie: irreducible reduced), we denote 
by $[V]$ (or sometimes simply $V$) the the correspopnding point $t\in\mathcal C$$(X)$.

We need some preliminary definitions and facts.

\medskip

\subsection{\bf  Covering Families} \label{covfam}

Let $X$ be a compact connected
normal analytic space (not necessarily in $\mathcal C$).

We shall consider finite or countable unions $T\subset \mathcal C$$(X)$
of compact irreducible subvarieties
   $T_{i}\subset \mathcal C$$(X)$, such that:

{\bf (a)} V$V_t$ is irreducible for $t$ generic in $T$.

{\bf (b)} $X$ is the union of all $V_t$'s, for $t\in T$.

(Equivalently: if $V\subset X\times T$ is the incidence graph of the
family $(V_{t})_{t\in T}$, and
$q:V\ra T$ and $p:V\ra X$ the restrictions to $V$ of the
projections, then $p$ is surjective and for any index $i$,
$q_i$ has generic fibres which are irreducible, where $V_i$, $q_i$, 
and $p_i$ denote the incidence graph of the family $T_i$ , and
the restrictions of
$p$ and $q$ to $V_i$ respectively).\\

We say then that $T\subset\mathcal C$$(X)$ is a {\bf covering family} of
$X$. Notice that by Baire's Category Theorem, at least one of the
families
$T_{i}\subset\mathcal C$$(X)$ has to be a covering family of $X$.\\

Because the maps $q_i$ are equidimensional, we have the following important:

\begin{proposition} \label{irred} If $T$ is a covering family of $X$, then:
\begin{enumerate}
\item[a.] $V_t$ is connected, for every $t\in T$.
\item[b.] $V_i$ is irreducible, for every $i$. More generally:
\item[c.] If $V'_i$ is deduced from $q_i:V_i\ra T_i$ by any surjective base 
change by $\beta:T'_i\ra T_i$, with $T'_i$ irreducible, then $V'_i$ is irreducible, too.

\end{enumerate}

\end{proposition}

Sometimes, we shall denote by $V$ not the incidence graph of the family $T$,
but its normalisation, or a smooth model of it, without further notice.

\begin{remark}: All results in this appendix hold, together with their proofs without any 
change if $X$ is arbitrary, irreducible (not necessarily compact), provided the maps $q_i,p_i$ are 
all assumed to be proper (so the $T_i$'s are not necessarily compact either). This was observed 
and used in [C94].
\end{remark}

{\bf (14.2) T-Equivalence (see also [C81])} \label{teq}

Let $T\subset
\mathcal C$$(X)$ be a covering family of $X$. Let $R(T)$ be the
equivalence relation on $X$ which is generated
by the symmetric reflexive binary relation $R^1(T)$ on $X$ for which
$x,x'\in X$ are related if and only if there exists $t\in T$ such that $V_t$
contains $x$ and $x'$.

Thus $x,x'\in X$ are $R(T)$-equivalent iff they can be joined by (ie: are
contained in) a {\bf $T$-chain}. This is by definition a {\bf
connected} union of
finitely many of the $V_t$'s.

We say that $X$ is {\bf $T$-connected} if $R(T)$ has a single
equivalence class (ie: any two points of $X$ can be connected by some
$T$-chain).

For these definitions, we do not assume that $T$ is irreducible. We
only assume that each irreducible  component of $T$ is compact, and
that $T$ has at most countably many such components.

Observe that $R^1(T)\subset X\times X$ has the following
description: $R^1(T)$=$(p\times p)(V\times _{T}V)$, with obvious
notations. From this and \ref{irred}, we conclude:

\begin{corollary}\label{irred'}
If $T$ is a covering family of $X$, then for any $i$, $R^1(T_i)$ is irreducible.
\end{corollary}

The fundamental result on $T$-equivalence is now the following:

\begin{theorem} \label{quot}([C81])
 Let $T\subset \mathcal C$$(X)$ be a covering family of $X$. 
 Assume that $X$ is normal.

 There exists a fibration $q_{T}:X\ra X_{T}$ such that its general fibre is an equivalence
class for $R(T)$. Of course, this map is unique, up to equivalence. Moreover, $q_T$ is almost 
holomorphic.

The map $q_T$ is called
the {\bf $T$-quotient} of $X$.
\end{theorem}

\begin{remark} If $T$ has only finitely many irreducible components, the word 
``general" can be replaced by `` generic" in \ref{quot}. See below.

It is not necessary to assume that $X\in \mathcal C$ in \ref{quot}.
\end{remark}

\begin{remark} 

\end{remark} The above result is shown in [C6] only when
$T$ is irreducible. The general case stated here can however be
easily reduced to
this special case by arguments sketched below in \ref{irredcase}.

See also
[C94] for a more general version (but the proof goes through) in which
$X$ is not assumed to be compact ,
but $p_i$ is assumed to be only proper. This version is used to construct
$\Gamma$-reductions (also called Shafarevitch maps). These
$T$-quotients appear in
many other problems.

Finally, the normality assumption is essential, as shown by examples
given in [C81].

Theorem \ref{quot} is exposed in [Ko1] in the algebraic context, via an
approch slightly different from (and no appropriate quotation of) the 
original one.

A carefull, corrected exposition along
this line is given in [De].

The proof of thorem \ref{quot} will be done in two steps: first the proof of the case where 
$T$ has finitely many components, each covering; second: the the easy reduction to this special case.

\subsection{\bf Construction of the Quotient for an Irreducible Covering Family}

We shall prove here Theorem \ref{quot} when $T$ is irreducible. The proof is exactly the same 
as in [C81], except for a slight simplification at one point, which shortens the exposition. 
This simplification consists in replacing $T$ by $T'$ in step 2. below.

Let us remark that the proof works also if $T$ is not necessarily irreducible, 
provided all components $T_i$ of $T$ are covering.

We thus start with a diagram:

\centerline{
\xymatrix{ 
V \ar[r]_-{p}\ar[d]_-{q} & X\\  
T  \\ 
}
}

in which $V,T,X$ are irreducible compact, $p,q$ surjective, and $q$ is equidimensional with generic 
fibre irreducible. Also $X$ is normal.

Because $T$ is fixed, we delete $T$ from the notations introduced in \ref{teq}, 
obtaining $R,R^1$ instead of $R(T),R^1(T)$.

We define an {\bf $n$-chain (of $T$)} as the connected union of $n$ of the $V_t$'s.

We say that $x,x'\in X$ are {\bf $n$-equivalent (for $T$)} if they are connected by (ie: contained in) 
some $n$-chain of $T$.

We denote by $R^n$ (or $R^n(T)$) the symetric set of $X\times X$ consisting of pairs $x,x'\in X$ which 
are $n$-equivalent. When $n=0$ (resp. $n=\infty$), we define $R^n$ as being $\Delta_X$, the diagonal of $X$
(resp. $R=\cup_{n\geq0} R^n$).

The sequence $R^n$ is of course increasing for the inclusion.

\begin{lemma}\label {R^n} For any $n\geq 1$, one has:
\begin{enumerate}
\item[1.] $R^{n+1}=(p)'((q')^{-1}((p)')^{-1}(R^n)$, 
where $(p)'$ and $(q'):V\times_{T}V\ra X=\Delta_X$  
denote the restrictions to $V\times_{T}V$ of 
$(p\times p):(V\times V)\ra (X\times X)$ and $(q\times q)$ respectively.

\item[2.] $R^n\subset (X\times X)$ is Zariski closed.

\item[3.] $R^1$ is irreducible.(see \ref{irred}).
\end{enumerate}
\end{lemma} 

{\bf Proof:} Easy. See [C81] for details $\square$

\begin{definition}\label{signR} 
We shall now say that an irreducible component of $R^n$ is {\bf significant} if it contains $\Delta_X$.
\end{definition}

We note $R_n\subset R^n$ the union of the significant irreducible components of $R^n$. 
Thus $R_n\subset X\times X$ is Zariski closed, and $R^n$ has only finitely many significant components, for any $n$. 
Notice that $R_1=R^1$ is irreducible (if $T$ is).

We also note $d_n(T):=d_n:=dim(R_n)$, for every $n\geq 0$. This is an increasing, bounded sequence of integers. 
So it stabilises at some $n_0$.

\begin{definition}\label{stat} We shall say that the irreducible covering family $T$ of $X$ is {\bf stationary} if:
\begin{enumerate}
\item[a.] $p$ is a modification
\item[b.] $d_2=d_1$
\end{enumerate}
\end{definition}

Theorem \ref{quot} is proved in two steps, corresponding to the following two propositions:

\begin{proposition}\label{stateq} Let $T$ be an irreducible covering family of $X$, normal. 
There exists a {\bf stationary} irreducible covering family $T'$ of $X$ such that $R(T)=R(T')$.
\end{proposition}

Propsition \ref{statquot} will be proved below.

\begin{proposition}\label{statquot} Let $T$ be an irreducible stationary covering family of $X$, normal. 
Then:
\begin{enumerate}
\item[a.] $h:=q\circ (p^{-1}):X\mero T$ is an almost holomorphic fibration.
\item[b.] $h:X\mero T$ is the $T$-quotient of $X$.
\end{enumerate}
\end{proposition}

{\bf Proof (of \ref{statquot}:)} The second assertion is an immediate consequence of the first, because 
$h$ is then holomorphic near each of its regular fibres (the ones not meeting the indeterminacy locus $I_h$ 
of $h$). So that each such fibre is then an equivalence class for $R(T)=R$.

Everything is clear from the definitions and assumptions in the first assertion, except for the fact that $h$ is almost holomorphic.
We check this by contradiction, assuming that each fibre $F=V_t$ of $h$ meets $I_h$. This implies that $F$ also meets $J$, 
the indeterminacy locus of $p^{-1}$. Let $x\in J\cap F$. Because $X$ is normal, $p^{-1}(x'):=W'$ is positive dimensional at each of 
its points; in particular at $x':=(t,x)$. Let $W$ be an irreducible component of $W'$ containing $x'$. 

Let $S:=p(q^{-1}(q(W)))\subset X$, and let $y\in F=V_t$ be generic. Then $S\times y\subset R^2(T)$ (because $(x,y)$ and $(s,x)$ are in 
$R^1(T)$, if $s\in S$, by construction. More precisely, if $F=V_t$ is irreducible reduced, then $S\times y\subset R_2(T)$, because $S$ is 
then irreducible and $t\in T$ generic. 

We now get that $d_2\geq d+m$, if $d:=dim(S)$, and $m:=dim(X)$, while $d_1=f+m$, if $f:=dim(F)=(dim(V)-m)$. 
To get a contradiction, it is thus sufficient to show that $d>f$. But this is obvious because:
 $dim(S)=dim(p(q(q^{-1}(W)))=dim(q(q^{-1}(W))=dim(W')+f>f$. The second equality comes from the fact that $p$ is a modification, 
and that $F$ is not contained in $J$, but is contained in $q(q^{-1}(W))$. 
$\square$.

We now start the:

{\bf Proof (of \ref{stateq}):} We shall denote, for any $A\subset X\times X$ by $r'$ (resp. $r"$) the restriction to $A$ of the first 
(resp. second) projection of $ X\times X$ to $X$, and by $A'(x)$ (resp. $A"(x)$) the fibre of $r'$ (resp. $r"$) above $x\in X$.

Let $s:X\times X\ra  X\times X$ be the involution exchanging the factors. We say that $A$ is {\bf symmetric} if $s(A)=A$. In this case, 
$A'(x)=A"(x),\forall x\in X$, and is simply denoted $A(x)$. 

The main step in the proof of \ref{stateq} is: 

\begin{lemma}\label{signirred} Let $n\geq 0$ be such that $d_n=d_{2n}$. Then:
\begin{enumerate}
\item[a.] $R_n$ is symmetric and irreducible.
\item[b.] $R_n(x)$ is irreducible and locally irreducible at $(x,x)$, for $x$ generic in $X$.
\end{enumerate}
\end{lemma}

Before we show \ref{signirred}, let us prove that it implies proposition \ref{stateq}:

Let $r':R_n\ra X$ be the first projection. By \ref{signirred}, it is a fibration. Let $g:X\mero \mathcal C$$(X)$ be 
the associated map, sending $x\in X$ generic to $R_n(x)$. Let $T'\subset \mathcal C$$(X)$ be the image of $g$. 
This is an irreducible covering family of $x$. 

It thus only remains to observe the following two properties:

\begin{enumerate}
\item[a.] $R(T)=R(T')$É

\item[b.] The family $T'$ is stationaryÉ
\end{enumerate}

Let us show [a.]: For $x$ generic in $X$, $R_n(x)=(R')^1(x)$, by definition, with $R':=R(T')$. Because $(R')^1$ 
is irreducible (by \ref{irred}), we conclude that $R^1(T')= R_n(T)\subset R(T)$. Thus $R(T')\subset R(T)$. 

In the opposite direction, because $R^1(T)$ is irreducible and contains the diagonal, $R^1(T)\subset R_n(T)=R^1(T')$. 
Thus: $R(T)\subset R(T')$.

(Observe this argument still works if $T$ consists of finitely (or infinitely!) many components $T_i$, each of which is covering 
(so that $R^1(T_i)$ contains the diagonal. But it does not work in general).\\

Let us next show [b.]: Because $R^1(T')=R_n(T)$, we have: $R^2(T')\subset R^{2n}(T)$. Thus also: 
$R_{2}(T')\subset R_{2n}(T)$. So: $d_2(T')\leq d_{2n}=d_n=d_1(T')$. Hence [b.].$\square$\\

To conclude the proof of Theorem \ref{quot}, it thus only remains to give the:

{\bf Proof (of Lemma \ref{signirred}):} We first prove by contradiction 
that $R_n(x)$ is irreducible at $(x,x)$, for $x\in X$ generic.

Otherwise,  
$R_n$ had two distinct local irreducible components 
$S_1,S_2$ at $(x,x)$, with $S_1$ say of dimension $d_n$;
 and $S_1(x),S_2(x)$ are irreducible at $(x',x')$, for $x'$ close to $X$,
 by Lemma \ref{irredfib} below (applied to each local irreducible component of $R_n$). 

We have a natural holomorphic map: $h:=W:=S_1(x) \times_{X} s(S_2)\ra R_{2n}(x)$ defined by sending 
$((x,y),(y,z))$ to $(x,z)$. Shrinking suitably the domain $W$ of this map, we can assume that it is irreducible, 
because the generic fibres of $W$ over $X$ are irreducible. But now $h(W)$ contains $S_1(x)\cap U$ (take $z=x$ and $S_2(x)\cap U$ 
(take $y=x$), if $U$ is a sufficiently small open neighborhood of $x$ in $X$. 

Because $S_2(x)$ is not contained in $S_1(x)$, by hypothesis, the rank of $h$ at the generic point of $W$ is thus 
strictly larger than: $dim(S_1)-m=d_n-m$, if $m:=dim(X)$. Because this rank is at most 
$d_{2n}-m$, this implies that $d_{2n}>d_n$, contradicting our assumptions.

The others assertions of \ref{signirred} are immediate consequences of lemma \ref{irredfib} below. 
$\square$.\\

To complete the proof of Theorem \ref{quot} for $T$ irreducible, it thus simply remain to show the:

\begin{lemma}\label{irredfib} Let $g:V\ra W$ be a surjective holomorphic map between 
complex spaces, with $\sigma:W\ra V$ a holomorphic section. Let $W':=\sigma(W)$.

If $V$ is lacally irreducible at the generic point of $W'$, then for $w\in W$ generic:
\begin{enumerate}
\item[a.]  The $g$-fibre $V_w$ over $w$ is locally irreducible at $w':=\sigma(w)$.
\item[b.]  $V_w$ is irreducible, if $g$ is proper.
\end{enumerate}
\end{lemma}

{\bf Proof:} We shall give, by lack of suitable reference a proof based on desingularisation. The only advantage is its shortness. 
Elementary algebro-geometric proofs were longer (see [C81] for an analytic proof).

Proof of [a.]: let $d:V'\ra V$ be a smooth model of $V$, and $g':=g\circ d:V'\ra W$. We can assume $W$ to be smooth connected. 
The fibres of $d$ over $W'$ are thus connected, by assumption (and Zariski's Main theorem). 
Let $W":=d^{-1}(W')\subset V'$. The rank of $g'$ is 
then maximal, equal to $m:=dim(W)$, at each $w"\in W^*\subset W"$, where $W^*$ 
is Zariski open and contains $W"\cap ((g')^{-1}(U))$, 
where $U$ is a dense Zariski open subset of $W$. This holds true by standard Zariski 
lower semi-continuity properties of the rank, and 
Bertini theorem.

For $w\in W$ generic, we conclude that the fibre $V'_w$ of $g'$ over $w$ is connected and smooth near $W_w"$. 
Thus $V_w=d(V'_w$ is irreducible near $w'$. Hence [a.].

Proof of [b.]: Let $V'_w$ be a smooth generic fibre of $g'$. We saw that $W_w":=d^{-1}(w')$ is connected. 
It is thus contained in a unique 
connected component $\overline{w}$ of $V'_w$. Let $g':=h\circ \bar{g}$ be 
the Stein factorisation of $g'$, with $\bar{g}:V'\ra \overline{W}$ connected and 
$h:\overline{W}\ra W$ finite. Now the previous observation defines 
a natural meromorphic section $\bar{h}$ of $h$, which sends $w\in W$ generic to $\overline{w}$.
This shows [b.].$\square$\\

The above proof of Theorem \ref{quot} easily gives an estimate of 
the length of the $T$-chains connecting two general $T$-equivalent points of $X$ 
(see [Ko-Mi-Mo], and [De] for similar results):

\begin{proposition}\label{length} Let $T$ be an {\bf irreducible} covering family of $X$, normal. 
Let $q_T:X\ra X_T$ be the $T$-quotient of $X$, and $F$ a regular 
(ie: contained in the regular locus of $q_T$) fibre of $q_T$. 
Then any two points of $F$ are joined by a $T$-chain of length at most $l:=2^{(f-r)}$, where $r:=dim(V_t)$, and 
$f:=dim(X)-dim(X_T)$ is the dimension of a generic fibre of $q_T$.

In particular:  $l\leq 2^{f-1}$, unless $T$ is the family of points of $X$. 
\end{proposition}

Actually, the statement holds for {\bf any} Chow-theoretic fibre $F$ of $q_T$, as the proof shows.

{\bf Proof:} Assume the sequence $d_1,d_2, d_4,..,d_{2^k}$ is strictly increasing, while $d_{2^k}=d_{2^{k+1}}$. 
Thus: $d_{2^k}\geq d_1+k$.
From the proof of \ref{quot}, we see that $d_1=m+r$, while $d_{2^k}=m+f$, if $m:=dim(X)$.
Thus: $m+f=d_{2^k}\geq d_1+k=m+r+k$, and $k\leq (f-r)$. Hence the claim, since by construction, 
two points of $F$ 
are connected by a $T$-chain of length $2^k$. 
$\square$.

\begin{remark}
\end{remark}
 The proofs given here work also when $T$ has all of its components $T_i$ covering, if $r$ denotes the infimum 
of the corresponding $r_i$'s. So, in this case, we always get: $l\leq 2^{f-1}$.
The estimate is then in this case slightly sharper than the one ($l\leq 2^f-1$) given in loc.cit. 

The estimate $l\leq 2^f-1$
 can easily be deduced but apparently not improved (see section \ref{irredcase} below)
in the general case from \ref{length}, and the reduction process described in the next section.

\subsection{\bf  Reduction to the case where $T$ is irreducible} \label{irredcase}

We show how to reduce the proof of Theorem \ref{quot} to the case where $T$ is irreducible. 

We shall consider the set $\mathcal F$$(X,T)$ of almost holomorphic fibrations $f:X\mero Y$ which are {\bf subordinate to $T$}. 
This means, by definition, that $R(f):=R(Y)\subset R(T)$. Here $Y$ is seen as the covering family parametrised by the family of fibres of $f$. 

Notice that the identity map of $X$ is in $\mathcal F$$(X,T)$. 

Let $T':=T_i$ be a component of $T$, and $V'$, $p'$,$q'$ the
data associated to $T'$, as in \ref{covfam}.

Let $f\in \mathcal F$$(X,T)$. We say that $T_i$ is {\bf $f$-covering} if 
$f\circ p':V'\ra Y$ is surjective.\\

The main step of the reduction process is contained in:

\begin{lemma} \label{subofib} Assume that $T'$ is $f$-covering for some $f$ subordinate to $T$.
 Then either:
\begin{enumerate}
\item[1.] $R(T')\subset R(f)$, or:
\item[2.] There exists a fibration $g:Y\mero Y'$ with $dim(Y')<dim(Y)$, and such that $f':=g\circ f:X\mero Y'$ is 
almost holomorphic and subordinate to $T$. (Moreover: $R(T')\subset R(f')$).
\end{enumerate}
\end{lemma} 

Before proving \ref{subofib}, we show how it implies the existence of $q_T$ in the general case:

$\square$ Take $f:X\mero Y$ be in $\mathcal F$$(X,T)$, and such that $dim(Y)$ is minimum. 
Let $J$ be the set of $i$'s such that $T_i$ is not 
$f$-covering. Let $T'$ be the union of all components $T_i$ of $T$ which are $f$-covering. 

From the preceding lemma, we see that $R(T_i)\subset R(f)$ for each $i\notin J$. For a general fibre $F$
of $f$, we thus see that $F$ is saturated for $R(T')$, and moreover that $R(T')=R(T)$ on this fibre, because 
the $T_i$'s are not $f$-covering if $i\in J$. Thus $F$ is saturated for $R(T)$. Finally, 
because $f$ is subordinate to $T$, we see that $F$ is a single $R(T)$ eaquivalence class $\square$

{\bf Proof of \ref{subofib}:} Let thus $T'$ be $f$-covering, and 
assume that $R(T')$ is not contained in $R(f)$. So $dim(V_{t'})>0$ if $t'$ is generic in $T'$.

 In this
case, define for $t'\in T'$ generic the new reduced cycle
$V'_{t'}\subset X$ by:
$V'_{t'}$:=$(f)^{-1}(f(V_{t'}))$. This cycle is irreducible,
because so are $V_{t'}$ and the generic fibre of $f$.
Let now $T"\subset \mathcal C$$(X)$ be the closure in $\mathcal C$$(X)$
of the points parametrising the cycles $V'_{t'}$ so obtained.

It is obvious to check that $T"$ is an irreducible covering
family of $X$, and that we have a natural factorisation map:
$g:Y\ra Y':=X_{T"}$
such that $f'=g\circ f$ if $f':X\ra Y'$ is the
$T"$-quotient of $X$. Obviously: $dim(Y')<dim(Y)$ because: $dim(V_{t'})>0$ if $t'$ is generic in $T'$. 

Moreover, we have an inclusion:
$R^1(T")\subset R^1(T)$, because $R(T")$ is generated by $R(f)$ and $R(T')$. 
This inclusion extends to $R(T")\subset R(T)$, 
so that $f'$ is subordinate to $T$, as claimed.

Let us show that $R(T")$ is generated by $R(f)$ and $R(T')$: The inclusions of $R^1(f)$ and $R^1(T')$ in 
$R^1(T")$ are obvious, by construction. So that $R(f)$ and $R(T')$ are contained in $R(f')$. 
The reverse direction is immediately given by the:

\begin{lemma}\label{R(T")} 
\begin{enumerate}
\item[1.] $R^1(T")$ is the main component of $(f\times f)^{-1}((f\times f)(R^1(T'))$.
\item[2.] $R^1(T)\subset R^2(S)$ if $S$ is the covering family with the two irreducible 
components $Y$ and $T'$. 
\end{enumerate}
\end{lemma}

{\bf Proof:} The first assertion follows immediately from the definition of the generic $V'_{t'}$ and 
the fact that $R^1(T")$ is irreducible, while the second follows from the first and Lemma 
\ref{R^n}.$\square$.\\

\begin{remark}\label{length'} ([De],[Ko-Mi-Mo]) Any two points of any Chow-theoretic fibre $F'$ of 
$q:X\mero X_T$ can be connected by a $T$-chain of length at most: $l'\leq 2^b-1$, if 
$b=dim(X)-dim(X_T)$ is the dimension of $F'$.
\end{remark}

{\bf Proof:} By the reduction process just described above, and using the same notations, we just need 
to show that if this holds for $f$ in place of $q_T$, this holds also for $f'$, so that 
$R(f')$ is generated by $R(f)$ and $R(T')$, $T'$ an irreducible $f$-covering family of $X$.

By \ref{length}, to points of any fibre (Chow-theoretic, as always) $F'$ of $f'$ are connected by a $T"$-chain of length 
$l"\leq 2^{b'-r'}\leq 2^{b'-b-1}$, where $b':=dim(F')$, $r':=dim(V'_{t'}\geq b+1$, with $b:=dim(F)$, 
$F$ any fibre of $f$. By the definition of $V'_{t'}$, we immediately see that any two of its points 
are connected by a $T$-chain of length at most $(2l+1)$, if two points of any fibre of can be connected 
by a $T$-chain of length $l$.

By the hypothesis on $f$, we can take $l:=(2^b-1)$. So that any two points of any fibre of $f'$ can be 
connected by a $T$-chain of length at most $l'\leq l".(2l+1)\leq 2^{b'-b-1}.(2.(2^b-1)+1)\leq (2^{b'}-1)$,
 as claimed $\square$.

\subsection{\bf Zariski Regularity}\label{Zreg}

\begin{definition}\label{Zregdef} Let $A\subset S$, where $S$ is a complex space. We say that $A$ is 
{\bf Z-regular (in S)} if for each irreducible Zariski closed subset 
$T$ of $S$, $A\cap T$ either contains the general point of $T$, 
or is contained in a countable union of Zariski closed proper 
subsets of $T$. (In which case we say that $A$ is of first Zariski category in $T$).
\end{definition}

There are many examples of this situation (actually it rather seems that counterexamples are unnatural in 
algebraic or analytic geometry):

\begin{theorem} \label{Zregex} 

Let $X\in \mathcal C$. Let $A_P(X)=A_P\subset \mathcal C$$(X)$ be the set of points $t$ parametrising the set of 
irreducible reduced subvarieties $V_t$ of $X$ possessing a certain property $P$. 
(By convention, points of $X$ are always supposed to posses $P$).

Then $A_P(X)=A_P$ is Z-regular in $\mathcal C$$(X)$ if $P$ is any one of the following properties.
\begin{enumerate}
\item[1.] $V_t$ is special;
\item[2.] $V_t$ is a (possibly singular) rational curve.
\item[3.] $V_t$ has algebraic dimension at least $d$, a given nonegative integer.
\item[4.] $\pi_1(\hat{V_t})_X\in \mathcal G$, where $\mathcal G$ is a class 
of (isomorphism classes) of groups, and the first term is the 
image in $\pi_1(X)$ of the fundamental group of the normalisation of $V_t$.
\end{enumerate}
\end{theorem}

{\bf Proof:} For 1 (resp. 3 ; 4) see corollary \ref{specZR} and 
proposition \ref{Cfib} below (resp. [C80]; [C01]). Assertion (2) is easy.

\begin{proposition}\label{Cfib}  Let $P$ be a property such that 
if $f:X\ra Y$ is any fibration such that $X\in \mathcal C$, 
then the set of irreducible fibres of $f$ having $P$ either contains the general point of $Y$, 
or is of first Zariski category in $Y$. 

Then: for any $X\in \mathcal C$, $A_P(X)$ is Z-regular.
\end{proposition}

{\bf Proof:} Let $T\subset \mathcal C$$(X)$ be irreducible Zariski closed. 
Let $V=V_T$ be the incidence graph of the family $T$, with 
projections $p:V\ra X$ and $q:V\ra T$. By [C80], $V\in \mathcal C$. Thus: $A_P(X)\cap T$
 either contains the general point of $T$, or is of first Zariski 
category in $T$, by the hypothesis made on $P$ $\square$.

\begin{proposition}\label{strat}
Let $A\subset S$ be Z-regular. There exists a countable or finite family of 
Zariski closed irreducible subsets $S_i$ of $S$ 
such that:
\begin{enumerate}
\item[1.] For each $i$, $A_i:=A\cap S_i$ contains the general point of $S_i$.
\item[2.] $A$ is the union of the $A_i$'s.
The family $S_i$ is moreover unique if it is {\bf irredundant}, that is: if 
there is no inclusion among the $S_i$'s. 
The $S_i$'s are called the {\bf components of $A$}.
\end{enumerate}
\end{proposition}

{\bf Proof:} Let $\overline{A}$ be the Zariski closure of $A$ in $S$, and 
let $A'$ be an irreducible component of $\overline{A}$.
Consider $A\cap A'$: either it contains the general point of $A'$, or it 
is contained in a countable union $B$ of (Zariski) closed subsets of $A'$.

In the first case, $A'$ is a component of $A$. Proceed the $\overline{A}$, after removing its component $A'$.

In the second case, remove $A'$ from $\overline{A}$, and replace it by $B$. 

Repeat this process, observing that for each component of $\overline{A}$, we let at each step decrease the 
diension of the remaining components. So the process ends in finitely many steps for each such component.

The uniqueness statement is easy to check: assume we have two irredundant families $S_i$ and 
$S'_j$ of components of $A$. It is sufficient to show that any $S_i$  is contained in some 
$S_j'$. If not, $(A\cap S_i)$ is contained in the countable union of the $(A\cap S_j')$. It is 
then of first Zariski category in $S_i$. This is a contradiction to the definition of a component of $A$.
$\square$

The notion of Z-regularity has application to the consruction of meromorphic quotients:

\begin{theorem}\label{regquot}

Let $X\in \mathcal C$ be normal. Let $A\subset \mathcal C$$(X)$ be Z-regular. Let $T:=T(A)$ be the family of components of $A$, as 
defined in \ref{strat}. (If $T$ is not covering, we add the family of points of $X$).
Let $q_A:X\mero X_{T(A)}=X_A$ be the $T$-quotient of $X$. 

Let $t\in A$ such that $V_t$ meets some general fibre $F$ of $q_T$. Then $V_t$ is contained in $F$.

(We shall call also $q_T$ the $A$-reduction of $X$).

\end{theorem} 

{\bf Proof:} By the defining property of $q_T$, its general fibre $F$ is an equivalence class for $R(T(A))$. In particular, 
any $V_t$, $t\in A$, is either disjoint from, or contained in $F$.

\begin{remark} Two points in such a fibre are joined by a $T(A)$-chain. Such a $T(A)$-chain turns out to be an $A$-chain 
if $A$ is stable under specialisation, in the sense defined just below.
So, it is important to know if the property $P$ defining $A$ is stable under specialisation (ie: 
if each component of $V_t$ has the property $P$ when $t\in S_i$, if $S_i$ is a component of $A_P$, and $A_P$ is 
Z-regular.
\end{remark}

 For example, this is easily seen to be the case for properties 2,3 in \ref{Zregex} above. The case of 1. is not presently known, 
but a positive answer to the following question is expected:

\begin{question} \label{specspec}
Let $f:X\ra Y$ be a special fibration, with $X\in \mathcal C$. Is each irreducible component of $X_y$ special, for 
each $y\in Y$?
\end{question}

An application of \ref{regquot} above (obtained explicitely in [C92], and in [Ko-Mi-Mo] in the algebraic case) is:

\begin{theorem}\label{ratquot} Let $X\in \mathcal C$ be normal. There exists an almost holomorphic map 
$r_X:X\ra R(X)$ called the {\bf rational 
quotient} of $X$ such that: 
\begin{enumerate} 
\item[1.] The fibres of $r_X$ are rationally connected 
(ie: any two of their points are contained in a connected union od finitely many 
rational curves)
\item[2.] If a rational curve $C$ inside $X$ meets 
some general fibre $F$ of $r_X$, then $C$ is contained in $F$. 
\end{enumerate}
\end{theorem}

Moreover, $A$-reductions can be constructed in a relative situation as well.

\begin{theorem}\label{relquot} Let $X\in \mathcal C$ be normal. Let $f:X\ra Y$ be a holomorphic fibration.
Let $A\subset \mathcal C$$(X)$ be Z-regular. 
Let $A_f\subset A$ be the set of $t\in A$ such that $V_t$ is contained in some fibre of $f$.
Then $A_f\subset \mathcal C$$(X)$ is also Z-regular.

If $q_{A_f}:X\mero X_{A_f}$ is the $A_f$-reduction of $X$, 
there exists a factorisation $h_{A_f}:X_{A_f} \mero Y$ such that $f=h_{A_f}\circ q_{A_f}$.

Moreover, for $y\in Y$ general, the restriction of 
$q_{A_f}$ to $X_y$ is the $A_y$-reduction of $X_y$, if $A_y:=A_f\cap \mathcal C$$(X_y)$.

The map $q_{A_f}$ is called the {\bf $A$-reduction} of $f$.
\end{theorem}

{\bf Proof:} The fact that $A_f$ is Z-regular is obvious, because $A$ is Z-regular, and $A_f$ is 
the intersection of $A$ with its Zariski closed subset $\mathcal C$$(X/Y)$ consisting of 
$t's$ such that $V_t$ is contained 
in some fibre of $f$. We thus have a quotient $q_{A_f}$ for $A_f$.

The natural map $f_*:\mathcal C$$(X/Y)\ra Y$ now defines the factoristaion $h_{A_f}$, 
because the equivalence classes for $T(A_f)$ 
are obviously contained in fibres of $f$. If one prefers: $R(f)$ contains obviously $R(T(A_f))$.

We shall now show that, for $y$ general in $Y$, $(T(A_f))_y=T(A_y)$.
This will immediately imply that $(X_y)_{A_y}=(X_y)_{(A_f)_y}$, and so the claim.

Notice first that, by definition, $A_y=(A_f)_y$, for any $y$. Next, if $T'(A_f)$ is the union 
of all components of $T(A_f)$ surjectively mapped to $Y$ by $f_*$, then $(T'(A_f))_y=(T(A_f))_y$ for 
$y$ general in $Y$. It is thus sufficient to show that $(T'(A_f))_y=T(A_y)$ for $y$ general in $Y$.

This we now show:first, we have, for $y$ general: $A_y=(A_f)_y\subset (T'(A_f))_y$. This gives the 
first inclusion:$T(A_y)\subset (T'(A_f))_y$. To show the reverse, let $T$ be a component of $T'(A_f)$. 
Then: $T_y\subset T(A_y)$ for general $y$. Indeed: $T/(A_f\cap T)$ is contained in a countable union of 
Zariski closed proper subsets $B_j$ of $T$. For each $j$, the set $C_j$ of $y\in Y$ such that 
$B_j\cap T_y$ has empty interior in $T_y$ is Zariski closed proper in $Y$. 

If $C$ is the union of the $C_j$'s, and $y$ is not in $C$, we conclude 
that $(T'(A_f))_y\subset T(A_y)$, as claimed $\square$

The main applications of \ref{relquot} are to construct relative versions of the rational quotient, of the core, and other 
natural fibrations. For example: applying \ref{relquot} to the family of rational curves on $X\in \mathcal C$, we get:

\begin{proposition}\label{relratquot} Let $f:X\ra Y$ be a fibration, $X\in \mathcal C$. Then $f$ has a {\bf rational quotient}, 
that is: a factorisation $f=R_f\circ r_f$, where $r_f:X\ra R(f)$ is a fibration which induces the rational quotient 
$r_{X_y}:X_y\ra R(X_y)=(R(f))_y$ for $y$ general in $Y$.
\end{proposition}

\subsection{\bf Stability}\label{stab}

\begin{definition}\label{defstab} Let $A\subset \mathcal C$$(X)$ be Z-regular, 
for $X\in \mathcal C$. We say that $A$ is {\bf stable} if it has the following two properties:
\begin{enumerate} 

 \item[(stab1)] The general fibre of $q_T:X\ra X_T$ is in $A$ if $T$ is 
any {\bf covering} irreducible component of $A$ (in the sense of Definition \ref{strat} above).

\item[(stab2)] $t\in A$ if $V_t=V\subset X$ is a subvariety and $f:V\mero W$ is a fibration with 
general fibre $V_w$ in $A$, such that a subvariety $W'\subset V$ exists, with $[W']\in A$, and $f(W')=W$.
(By abuse of language, we say that a suvariety $Z$ of $X$ is 
in $A$ if the point $[Z]\in \mathcal C$$(X)$ which parametrises this subvariety is in $A$.)
\end{enumerate} 
\end{definition}

Let us give some examples:

\begin{theorem}\label{exstab} Let $A_P(X)\subset \mathcal C$$(X)$ be Z-regular, 
for $X\in \mathcal C$ be as in theorem \ref{Zregex} above. Then $A_P$ is stable in the following two cases:
\begin{enumerate} 
\item[1.] $P$ is the property of being special, if $X$ is smooth.
\item[2.]  $A_P$ is the set of $t$'s such that $\pi_1(\hat{V_t})_X\in \mathcal G$, 
where $\mathcal G$ is a class of (isomorphism classes) of groups which is stable in the sense of [C98]. (Such classes are 
the class of finite, or solvable groups, for example).
\end{enumerate} 
\end{theorem}

The proof is given in section \ref{core} for [1.], and in [C98] for [2.]. The basic result about stability is the following:

\begin{theorem}\label{quotstab} Let $A\subset \mathcal C$$(X)$ be Z-regular and stable,
for $X\in \mathcal C$ normal. Let $q_A:X\mero X_A$ be the $T(A)$-quotient of $X$. Then:
\begin{enumerate} 
\item[1.]  The general fibre of $q_A$ is in $A$.
\item[2.]  The general fibre $F$ of $q_A$ contains $V$, If $V$ is a subvariety of $X$ which is in $A$ and meets $F$.
\end{enumerate} 
\end{theorem}

{\bf Proof:} The second property is just one of the defining properties of 
$q_A$. To show the first one, we observe that, because of [stab1.], it is 
obviously satisfied for the quotient of $X$ by any of the covering 
irreducible components of $T(A)$, an so for the first step of the construction of  $q_A$, as 
exposed in section \ref{irredcase}.

We now go to the induction step (``the general fibre of $q_{n+1}$ is in $A$ if 
$q_n$ has this same property"). 

We use here this same section \ref{irredcase} to which we refer, with the same notations (and $T:=T(A)$). 

Then the subvarieties $V'_{t'}$ constructed there from the varieties $V_{t'}$ are also in $A$, by [stab2.] applied to 
$q'$ and $V_{t'}\in A$, by the hypothesis $T_{n+1}\subset T(A)$, using the induction hypothesis (``the general fibre 
of $q'$ is in $A$").

 And so, by [stab1.] again, we see that the general 
fibre of $q_{n+1}$ is indeed in $A$ $\square$

{\bn F.Campana

D\'epartement de Math\'ematiques.

Universit\'e Nancy 1

BP 239.

F-54506.Vandoeuvre-les-Nancy.Cedex.}}

\end{document}